\documentclass[12pt,epsf,bezier]{article}
\usepackage{latexsym}
\usepackage{amsfonts}
\newcounter{ppp}
\newcounter{band}
\newcounter{annulus}
\newcounter{nonannulus}
\newcounter{ovals}
\newcounter{innermost}
\newcounter{touch}
\newcounter{bigons}
\newcounter{spokes}
\newcounter{triple}
\newcounter{derived}

\newcounter{spone}
\newcounter{sptwo}
\newcounter{figone}
\newcounter{ab}

\newcounter{theta}

\newcounter{pdfour}

\newcommand{\Lab}{{\rm Lab}}
\newcommand{\yyy}{{\cal Y}}
\newcommand{\ttt}{{\cal T}}

\newcommand{\bb}{{\cal B}}
\newcommand{\topp}{{\bf top}}
\newcommand{\bott}{{\bf bot}}
\newcommand{\vk}{van Kampen }

\newcommand{\iv}{^{-1}}

\newcommand{\bn}{{\bf N} }

\newcommand{\rr}{{\cal R} }
\newcommand{\xxx}{{\cal X} }

\newcommand{\ooo}{{\cal O} }

\newcommand{\ww}{{\cal W} }
\newcommand{\pp}{{\cal P} }
\newcommand{\qq}{{\cal Q} }
\newcommand{\sss}{{\cal S} }
\newtheorem{thm}{Theorem}[section]
\newtheorem{prop}[thm]{Proposition}
\newtheorem{lm}[thm]{Lemma}

\newtheorem{df}[thm]{Definition}

\begin{document}

\title{Isoperimetric Functions of Groups and\\ Computational Complexity of the Word Problem}

\author{J.-C. Birget, A. Yu. Ol'shanskii,\\ E. Rips, M. V. Sapir\thanks{The research of the first and the fourth authors was supported in part by an NSF grant. The research of the second author was supported in part by the Russian fund for fundamental research 96-01-420.}}
\date{}

\maketitle

\begin{abstract}
We prove that the word problem of a finitely
generated group $G$ is in NP (solvable in polynomial time by a non-deterministic Turing machine) if and only if this group is a subgroup of a finitely presented group $H$ with polynomial isoperimetric function. The embedding can be chosen in such a way that $G$ has bounded distortion in $H$. This completes the work started in \cite{SBR} and \cite{Ol1}. 
\end{abstract}

\bigskip


\tableofcontents
\section{Introduction}

A function $f: \bn\to\bn$ is called an {\em isoperimetric function}
of a finite
presentation $\pp=\langle  X\ |\ R\rangle$ of a group $G$
if for every number $n$ and
every word $w$ over $X$ which is equal to 1 in $G$, $|w|\le n$, there exists a
\vk diagram over $\pp$ whose boundary label is $w$ and area $\le f(n)$;
in other words, $w$ is a product in the free group of at most $f(n)$ conjugates of
the relators from $R$. One can also define an isoperimetric function of a presentation $\pp$ as 
an isoperimetric function of its Cayley complex (the Cayley graph equipped with 2-cells corresponding to relations of $\pp$), that is a function bounding the filling area of a loop of a given perimeter in the Cayley complex  (see \cite{gromov3}, \cite{gromov4}, \cite{MO}, \cite{Gersten}, \cite{gromov1},
\cite{short}). Isoperimetric
functions of manifolds and 2-complexes were studied since ancient Greeks. Dehn probably was
the first to notice a connection between them and the word problem in
groups (this is how he solved the word problem in surface groups). But the concept of an isoperimetric function for groups (at least the
name of this concept) probably first appeared in \cite{gromov3}, \cite{gromov4}.

The smallest isoperimetric function of a finite presentation $\pp$
is called the Dehn function of $\pp$ \cite{Gersten}.

Let $f, g:{\bf N}\to {\bf N}$ be two functions. We write $f \preceq g$ if there
exist non-negative constants $a, b, c, d$ such that $f(n)\leq ag(bn)+cn+d$.  All
functions $g(n)$ which we are considering in this paper grow at least as fast as
$n$. In this case $f(n)\preceq g(n)$ if and only if $f(n)\leq ag(bn)$ for some
positive constants $a, b$.  Two functions $f, g$ are called {\em equivalent} if
$f \preceq g$ and $g \preceq f$.

It is well known that Dehn functions 
corresponding to different  finite presentations of the same group are
equivalent (see \cite{MO} or \cite{Alo}, \cite{Gersten}).  This
allows us to speak about {\em the Dehn function} of a finitely presented group.  

It is also known \cite{Gersten} that the Dehn function of a finitely presented group is recursive if and only if the group has decidable word problem. Moreover, a subset of the set of authors of this paper has proved in \cite{SBR} (see the proof of Theorem 1.1 of \cite{SBR}) that for every finitely presented  group $G$ with Dehn function $T(n)$ there exists a nondeterministic Turing machine $M(G)$ which solves the word problem in $G$ and has time function equivalent to $T(n)$. Roughly speaking, this Turing machine takes a word $w$ over the generators of $G$ and just inserts relators of $G$. It stops and accepts $w$ when it gets 1. 
Clearly this machine solves the word problem in every finitely generated subgroup of $G$ as well. Therefore if a finitely generated group $G$ is a subgroup of a finitely presented group with polynomial isoperimetric function then the word problem in $G$ is in NP (i.e. it can be solved by a non-deterministic Turing machine with polynomial time function). 

The drawback is that the word problem in a group $G$ can be easy to solve but the Dehn  function of $G$ can be huge. A typical example is the Baumslag-Solitar group $G_{2,1}=\langle a,t\ | \ a^t=a^2\rangle$ where $a^t=t\iv at$. This group has exponential Dehn function \cite{gromov1}. But this group is a subgroup of $GL(2, {\bf Q})$ ($t=\left(\begin{array}{ll}1/2 &0\\ 0& 1\end{array}\right)$, $a=\left(\begin{array}{ll}1 &1\\ 0 & 1\end{array}\right)$), so the word problem in $G_{2,1}$ can be solved in at most quadratic time: it is easy to see that the word problem of every finitely generated matrix group over 
a field of rational numbers can be solved in at most quadratic time by a deterministic Turing machine. In fact it is possible to solve the word problem there in time $n(\log n)^2(\log\log n)$ (this follows from the fact that the product of two $n$-digit numbers can be computed in time $n\log n\log\log n$ \cite{Knuth}). Using the construction in \cite{SBR} Sapir showed (unpublished) that there exists a finitely presented group with word problem in NP and Dehn function not even primitive recursive. 

Nevertheless, in this paper we prove that every finitely generated group $G$  with word problem in NP can be embedded into a finitely presented group $H$ with polynomial isoperimetric function. Thus if we can solve the word problem in $G$ using a very clever and fast Turing machine, then we can use the simple-minded but almost as fast  Turing machine $M(H)$ to solve the word problem in $G$.  Moreover we prove that the group $G$ can have {\em bounded length distortion} in $H$, that is the word length of any element of $G$ in $G$ is within a constant factor of the word length of this element in $H$. Here is the main result of this paper.

\begin{thm} Let $G$ be a finitely generated group with word problem solvable by a non-deterministic Turing machine with time function $\le T(n)$ such that $T(n)^4$ is superadditive (that is $T(m+n)^4\ge T(m)^4+T(n)^4$ for every $m,n$). Then $G$ can be embedded into a finitely presented group $H$ with isoperimetric function equivalent to $n^2T(n^2)^4$ in such a way that $G$ 
has bounded distortion in $H$. In particular, the word problem of a finitely generated group is in NP if and only if  this group is a subgroup of a finitely presented group with polynomial isoperimetric function. 
\label{th1}
\end{thm}

Notice that this theorem strengthens the main result of \cite{Ol1} where it is proved that every finitely generated recursively presented group $G$ can be embedded without distortion into a finitely presented group. 

The class of finitely generated groups with word problem in NP is very large. We have mentioned that it includes all matrix groups over ${\bf Q}$. It also includes
\begin{itemize}
\item All finitely generated matrix groups over arbitrary fields: this follows from the fact that every finitely generated field is a finite extension of a purely transcendental extension of its simple subfield, and the fact that the word problem in the ring of polynomials over ${\bf  Q}$ or ${\bf Z}/p{\bf Z}$ is solvable in polynomial time, 
\item Polycyclic and finitely generated metabelian groups because they are representable by matrices \cite{Segal}, 
\item Automatic groups (in particular, hyperbolic groups) \cite{WPG}, 
\item Groups of piecewise linear transformations of a line with finitely many rational        singularities (including the R. Thompson group $F$), 
\item Every finitely generated subgroup of a diagram group \cite{GS}, 
\item Every free Burnside group $B_{m,n}$ for sufficiently large odd exponent $n$ (see, for example, Storozhev's 
argument in Section 28 of \cite{Ol89}).
\end{itemize}
This class is closed under free and direct products. It is easy to see using Magnus' embedding that for every normal subgroup $N$ of a free finitely generated group $F$ if $F/N$ has word problem in NP (resp. P) then $F/N'$ has word problem in NP (resp. P). Therefore every free group in the variety of all solvable groups of a given class has word problem in P. 
It is an interesting question whether this class also contains all one-relator groups. There are of course finitely generated groups with word problem not in NP, for example groups with undecidable word problem. But these groups and in some sense ``artificial". So perhaps the class 
of groups with word problem in NP (which by Theorem 1.1 is the class of all subgroups of finitely presented groups with polynomial Dehn functions) can be considered as the class of ``tame" groups.

An example of an embedding of one group into another where lengths are not distorted but areas are distorted can be found in \cite{Gersten3}. Some examples of groups with big Dehn functions embeddable into groups with small Dehn functions can be found in \cite{MO} and \cite{BBMS}.
Our theorem and Theorem 1.1 from \cite{SBR} show that any finitely generated group can be embedded
into a finitely presented group with bounded length distortion but with close to maximal possible area distortion. Indeed, Theorem 1.1 from \cite{SBR} shows that the isoperimetric function of a group $H$ containing a given group $G$ cannot be smaller than the time complexity $T(n)$ of the word problem for $G$, and we show that $G$ can be embedded into a finitely presented group with Dehn function at most $n^2T(n^2)^4$.  

For matrix groups our theorem implies that every such group is embedded with bounded length distortion into a finitely presented group with Dehn function at most $n^{10+\epsilon}$ for every $\epsilon>0$. 
It is interesting to know the smallest Dehn function of a finitely generated group containing,
for example,  the Baumslag-Solitar group $G_{2,1}$. Is it possible to lower this bound for $G_{2,1}$ to $n^3$
or even $n^2$? 

Using the proof of Theorem \ref{th1}, in order to embed a finitely generated group with word problem in NP into a finitely presented group with polynomial isoperimetric function, one needs first construct a Turing machine which solves the word problem, then convert it into an $S$-machine, then convert the $S$-machine into a group. As a result the group we construct will have a relatively complicated set of relations. In some particular cases like the Baumslag-Solitar group $G_{2,1}$, the free Burnside groups $B_{m,n}$ and others, the authors can modify this construction and get simple presentations of groups with polynomial isoperimetric functions where these groups embed. 

Theorem 1.1 is also interesting from the logic point of view. One can consider a group as a logical system where the defining relations
are axioms, and the inference rules are constructing step by step van
Kampen diagrams. The van Kampen diagrams are then {\em proofs} of their boundary
labels. Then the Dehn function becomes the syntactic complexity of proofs. The computational 
complexity of the word problem is the complexity in
the metaworld of the group theory. The embedding of groups becomes a conservative
extention of theories. 

With this vocabulary, the result means that the complexity of proofs in
the outer world of groups (that is, in its metamathematical semantics)
after appropriate conservative extention to a finitely axiomatized theory becomes the complexity of proofs
in the syntactical sense (up to a polynomial correction).

In this formulation, one can ask a similar question for any logical system! 
Several results of this type for general logical systems can be found in \cite{Kraj}.
Notice also that the semigroup analog of Theorem \ref{th1} was obtained in \cite{Bi}.

\section{Preliminaries}
\label{sect2}

\subsection{Turing Machines and $S$-machines}
\label{tm}

Let $T(n)$ be a function such that $T(n)^4$ is superadditive. Let $G$ be a group with 
finite symmetric set of generators $A=A^{-1}$, such that the
word problem of $G$ is decidable in time at most $T(n)$ by a nondeterministic 
Turing machine $M$. By a Turing Machine we always assume a non-deterministic multi-tape 
Turing machine (see \cite{HU} or \cite{SBR}). Each tape is finite but can be expanded by adding cells at the right end
of it. One of these tapes is called the {\em  input tape}. Each tape has its own alphabet, the set of letters which can appear on this tape. A {\em configuration} of a Turing machine is the sequence of words written on the tapes, position of the head relative to the tapes, and the state of the head. In the {\em input configuration} of this machine a word is written on the input tape and all other tapes are empty. We shall assume that $M$ has only one {\em stop configuration} $c_0$ and all tapes in this configuration 
are empty.  The fact that $M$ solves the word problem in $G$ means that the alphabet of the input tape is $A$ and for every word $u$ in the alphabet $A$ the machine can {\em accept}  the input configuration with $u$ written on the input tape (that is get to the stop state starting with this input configuration) if
and only if $u$ is equal to 1 in $G$.

Let $M$ be a {\bf nondeterministic Turing} machine with time complexity 
$\le T(n)$ that decides if a word in the alphabet $A$ is equal to 1 in $G$. 
In Proposition 4.1 of \cite{SBR}, it was proved that $M$ is polynomially equivalent to 
a so called $S$-machine $\sss(M)$. Let us repeat the definition of an $S$-machine from \cite{SBR} and then explain in what sense $\sss(M)$ is equivalent to $M$ (see Proposition \ref{prop1}).

Roughly speaking, the difference between $S$-machines and ordinary Turing machines is that $S$-machines are almost ``blind". They ``see" letters written on the tape only when these letters 
are between two heads of the machine and the heads are very close to each other. If the heads are far apart, the machine does not see any letters on the tape.

We define $S$-machines as rewriting systems \cite{KharSap}.
Let $k$ be a natural number. A {\em hardware} of an $S$-machine $\sss$
is a pair $(Y,Q)$ where
$Y$ is a $k$-vector of (not necessarily disjoint) sets
$Y_i$, $Q$ is a $(k+1)$-vector
of disjoint sets $Q_i$, $\overline Q=\bigcup Q_i$ and $\overline Y=\bigcup Y_i$ are also disjoint.
The elements of $\overline Y$ are called {\em tape letters}, the elements of
$\overline Q$ are called {\em state letters}. All these letters together form the
{\em alphabet} of the machine $\sss$.

With every hardware $\sss=(Y,Q)$ we associate the 
{\em language of admissible words} $L(\sss)=Q_1F(Y_1)Q_2...F(Y_k)Q_{k+1}$
where $F(Y_j)$ is the language of reduced group words in the alphabet
$Y_j$. Notice that in every admissible word, there is exactly one representative of each $Q_i$ and these representatives appear in this word in the order of the indices of $Q_i$. The $\overline Y$-letters which appear between
a $Q_i$-letter and a $Q_{i+1}$-letter in an admissible word belong to $Y_i$, $i=1,2,...,k$.

If $0\le i\le j\le k$ and $W=
q_1u_1q_2...u_kq_{k+1}$ is an admissible word then
the subword $q_iu_i...q_j$ of $W$ is called the $(Q_i,Q_j)$-subword of $W$
($i<j$).

An $S$-machine with hardware $\sss$ is a rewriting systems. The objects
of this rewriting system are all admissible words.

The rewriting rules, or {\em $S$-rules}, have the following form:
$$[U_1\to V_1,...,U_m\to V_m]$$
where the following conditions hold:
\begin{description}
\item Each $U_i$ is a subword of an admissible word starting with
a $Q_\ell$-letter and ending with a $Q_r$-letter (where $\ell=\ell(i)$
must not exceed $r=r(i)$, of
course).
\item If $i<j$ then $r(i)<\ell(j)$.
\item Each $V_i$ is also a subword of an admissible word whose $Q$-letters
belong to $Q_{\ell(i)}\cup...\cup Q_{r(i)}$ and which contains a $Q_\ell$-letter
and a $Q_r$-letter.
\item If $\ell(1)=1$ then $V_1$ must start with a $Q_1$-letter and if
$r(m)=k+1$ then $V_n$ must end with a $Q_{k+1}$-letter (so tape letters
are not inserted to the left of $Q_1$-letters and to the right of $Q_{k+1}$-letters).
\end{description}

To apply an $S$-rule to a word $W$ means to replace simultaneously subwords
$U_i$ by subwords $V_i$, $i=1,...,m$. In particular, this means that our rule is
not applicable if one of the $U_i$'s is not a subword of $W$. The following convention
is important:

{\bf After every application of a rewriting rule,
the word is automatically reduced. We do not consider reducing of an admissible
word a separate step of an $S$-machine.}

For example, if a word is $$q_1aaq_2bq_3ccq_4$$ and
$q_i\in Q_i$, $a\in Y_1$, $b\in Y_2$, $c\in Y_3$ and the $S$-rule is
\begin{equation}\label{rule1}
[q_1\to p_1a\iv,  q_2bq_3\to a\iv p_2b'q_3c],
\end{equation}
where $p_1\in Q_1$, $p_2\in Q_2$, $b'\in Y_2$, then the result of the
application of this rule is $$p_1p_2b'q_3cccq_4.$$

With every $S$-rule $\tau$ we associate the inverse $S$-rule $\tau\iv$
in the following way: if $$\tau=[U_1\to x_1V_i'y_1,\ U_2\to x_2V_2'y_2,...,
U_m\to x_mV_m'y_m]$$ where $V_i'$ starts with a $Q_{\ell(i)}$-letter and
ends with a $Q_{r(i)}$-letter, then $$\tau\iv=[V_1'\to x_1\iv U_1y_1\iv,\
V_2'\to x_2\iv U_2y_2\iv,...,V_m'\to x_m\iv U_my_m\iv].$$

For example, the inverse of the rule (\ref{rule1}) is

$$[p_1\to q_1a,  p_2b'q_3\to aq_2bq_3c\iv].$$

It is clear that $\tau\iv$ is an $S$-rule,  $(\tau\iv)\iv=\tau$, and that
rules $\tau$ and $\tau\iv$ cancel each other (meaning
that if we apply $\tau$ and then $\tau\iv$, we return to the
original word).

The following convention is also important:

{\bf We always assume that an $S$-machine is symmetric, that is if an
$S$-machine contains a rewriting rule $\tau$, it also contains the rule
$\tau\iv$.}

When an $S$-machine works, it applies its rules to admissible words. The sequence of these admissible words and the rules is called a {\em computation} of the $S$-machine.
We define the history of
a computation of an $S$-machine as the sequence (word)
of rules used in this computation. 
A computation is called {\em reduced} if the history of this computation is
reduced, that is if two mutually inverse rules are never
applied next to each other.

As usual the {\em length} of a computation $W_1,...,W_n$ is $n$, and the {\em area} is
$\sum_i|W_i|$. 

The following statement is a immediate corollary of Lemma 3.1 and Proposition 4.1 from \cite{SBR}. Notice that Proposition 4.1 from \cite{SBR} contains many properties of $\sss(M)$
which we are not going to explicitely use in this paper, so we do not present them here. But we are 
going to use lemmas from \cite{SBR} which are based on these properties. 

\begin{prop}\label{prop1} Let $M$ be a non-deterministic Turing machine with time function
$T(n)$ accepting language $L\subseteq A^+$. Then there exists an $S$-machine $\sss(M)$ and a 1-1 function 
$\sigma$ from the set of configurations of $M$ to the set of admissible words of $\sss(M)$
which satisfy the following properties. Let $c_0$ be the stop configuration of $M$ and denote $\sigma(c_0)$ by $W_0$. Then

\begin{enumerate}
\item A configuration $c$ of the machine $M$ is acceptable by $M$ if and only if $\sss(M)$ can take
$\sigma(c)$ to $W_0$. 
\item If $u$ is an input word for the machine $M$, $|u|=n$, and $c(u)$ is the corresponding input configuration of $M$ then $\sigma(c(u))=z_0\alpha^n z_1uz_2\delta^nz_3\omega^nz_4$ 
for some fixed letters $\alpha, \delta,\omega$ not occurring in $u$ and some $\overline Q$-words $z_0,z_1,z_2,z_3,z_4$ which do not depend on $u$, such that the contents of these words are disjoint and do not contain letters from $u$ or $\alpha, \delta, \omega$. 
\item If an admissible word $W$ contains the subwords $z_0$, $z_1$,$z_2$, $z_3$, $z_4$, and the $S$-machine $\sss$ takes $W$ to $W_0$ then $W=\sigma(c(u))$ for some group word $u$ over $A$. 
\end{enumerate}
\end{prop}

\subsection{Groups Simulating $S$-machines}

Now let $M$ be a Turing machine solving the word problem in our group $G$. We can assume that the set of defining relations of $G$ contains all relations of the form $aa^{-1}=1$, $a\in A$.  Let $\sss=\sss(M)$. 
Then $\sss$ also solves the word problem in $G$ in the following sense. For every group word $u$ in the alphabet $A$, the word $u$ is equal to $1$ in $G$ if and only if the $S$-machine $\sss$ takes the word $\sigma(c(u))$ to $W_0$ (in this case we say that $\sss$ accepts $\sigma(c(u))$ ). 

Let $\Theta$ be the set of all rules of $\sss$. In every pair of mutually inverse rules from $\Theta$ we choose 
one which we call {\em positive}. The set of all positive rules of $\sss$ is denoted by $\Theta_+$. 

Let $\overline Y$ be the set of all tape letters of the machine  which can appear in the admissible words of $S$
and let $\overline Q$ be the set of all state letters of $\sss$. Notice that $A\subseteq \overline Y$. Let us fix a large natural number $N$. Up to Section \ref{sect7}
it will be enough to assume that $N\ge 9$,  the exact value of $N$ will be given in Section \ref{sect7}. With every admissible word $W$ of the $S$-machine $\sss$ we associate the following word 

$$\kappa(W)=\kappa_1W \kappa_2 W^{-1} \kappa_3W...\kappa_{4N-1} W \kappa_{4N} W^{-1}.$$

{\bf Remark 1. } Notice that this word will turn into the word $K(W)$ from \cite{SBR} if we identify $\kappa_{2N+i}$
with $\kappa_i^{-1}$, $i=1,...,2N$. 

Now we can define a group $G_N(\sss)$ similar to the one defined in \cite{SBR}. 

The set of generators of $G_N(\sss)$ is $\overline Y\cup\overline Q\cup \Theta_+\cup \{\kappa_1,...,\kappa_{4N}\}$.
The set of relators of $G_N(\sss)$ is divided into three groups.

{\bf 1. Transition relations}. These relations correspond to
elements of $\Theta_+$.

Let $\tau\in \Theta_+$, $\tau=[U_1\to V_1,...,U_p\to V_p]$.
Then we include relations $U_1^\tau=V_1,...,U_p^\tau=V_p$
into $\pp_N(\sss)$. Here $x^y$ stands for $y\iv xy$.
If 
the letters of  one of the component of the vector $Q$, say, $Q_j$, do not appear in any of the $U_i$ then also include the relations
$q_j^\tau=q_j$ for every $q_j\in Q_j$.

{\bf 2. Auxiliary relations}.

These are all possible relations of the form $\tau x=x\tau$ where
$x\in \overline Y$, $\tau\in \Theta_+$
and all commutativity relations of the form
$\tau \kappa_i=\kappa_i\tau$, $i=1,...,4N$, $\tau\in \Theta_+$.

{\bf 3. The hub relation:}

$$\kappa(W_0)=1.$$

We also assume that all cyclic shifts and the inverses of these relations are in the set of defining relations 
of $G_N(\sss)$. 

The remark we made above shows that the group $G_N(\sss)$ constructed here has a natural homomorphism 
onto the group $G_N(\sss)$ in \cite{SBR}. This homomorphism is induced by the map $\kappa_{2N+i}\to \kappa_i^{-1}$. The following Proposition contains several properties of the group $G_N(\sss)$ from \cite{SBR}. It is easy to check that all these properties are true for the group we deal with here. The proofs are exactly the same (almost identical), so we do not present them here. In fact we could avoid changing 
the presentation of $G_N(\sss)$ in this paper, but it would make the paper longer and more complicated, so we decided to make the change.

{\bf Remark 2.} We can preserve these properties while changing the presentation of $G_N(\sss)$ even further: we can make different copies of the word $W$ in the definition of $\kappa(W)$
written in disjoint alphabets. This is interesting because this would make the hub $\kappa(W_0)$ linear (each letter occurs only once).  Then the group $G_N(\sss)$ would be an HNN extension of a free group with several stable letters, factorized by a linear relation. This may have nice geometric consequences. See also Remark 3 in Section \ref{he}. 

\begin{prop}\label{prop2} The following properties of the group $G_N(\sss)$ hold:
\begin{enumerate}
\item The machine $\sss$ takes an admissible word $W$ to $W_0$ if and only if the word $\kappa(W)$ is equal to 1 in $G_N(\sss)$. This and part 1 of Proposition \ref{prop1} imply that a 
word $u$ in the alphabet $A$ is equal to 1 in our original group $G$ if and only if the word 
$\kappa(\sigma(c(u)))$ is 1 in $G_N(\sss)$. 

\item The Dehn function of the group $G_N(\sss)$ is equivalent to $T(n)^4$.
\end{enumerate}
\end{prop}

\subsection{Van Kampen Diagrams, Bands and Annuli}

In order to analyze van Kampen diagrams over group presentations we mainly use notation and definitions from Ol'shanskii \cite{Ol89}, and from papers \cite{SBR}, \cite{Ol1}. 

In particular,
we use auxiliary $0$-edges, i.e. edges labeled by 1, and 0-cells, i.e. 
cells corresponding to the trivial relations of the form $a1^na\iv 1^m$ where $a$ is one of the generators,
$m$ and $n$ are integers or to the relation $1^n$ (see \cite{Ol89} for details).  Thus some edges in van Kampen diagrams can have label 1. One can insert 0-edges and 0-cells in a diagram in
order to make a path, which touches
itself, absolutely simple. 

We call a \vk diagram {\em reduced} if it does  contain 0-cells and 0-edges and does not have a pair of cells having a common edge and such that the label of the boundary of the subdiagram consisting of these two cells is equal to 1 in the free group (\setcounter{pdfour}{\value{ppp}}see Figure \thepdfour):

\bigskip

\unitlength=1mm
\linethickness{0.4pt}
\begin{picture}(98.00,25.00)
\put(74.00,10.00){\vector(0,1){14.33}}
\bezier{200}(74.00,10.33)(50.00,19.00)(74.00,25.00)
\bezier{200}(74.00,10.33)(98.00,17.00)(74.00,25.00)
\put(71.00,17.00){\makebox(0,0)[cc]{$u$}}
\put(59.00,17.67){\makebox(0,0)[cc]{$v$}}
\put(88.67,17.33){\makebox(0,0)[cc]{$v$}}
\put(68.67,23.33){\vector(3,1){2.67}}
\put(79.33,23.00){\vector(-2,1){2.33}}
\end{picture}
\begin{center}
\nopagebreak[4]
Fig. \theppp.
\end{center}
\addtocounter{ppp}{1}

If a diagram has such a pair of cells then these two cells can be removed and the corresponding subdiagram can be replaced by a diagram without cells. This operation is called {\em canceling} a pair of cells.

If $\Delta$ is an (ordinary) \vk diagram then $\partial(\Delta)$ denotes its
boundary.  

We always assume that boundaries of \vk diagrams and boundaries
of cells are oriented clockwise (this  is an insignificant difference
with \cite{Ol89}). The {\em contour} of a cell or a diagram is the union
of the boundary and its inverse.

A path in the 1-skeleton of a diagram $\Delta$ is called {\em simple} if it does not cross itself after an arbitrary small deformation.
A path is called {\em reduced} if it does not contain consecutive mutually
inverse edges.

The length of a path $p$, denoted by $|p|$
is the number of non-zero edges in it.

We shall use the following operations which can be performed on
arbitrary diagrams.

\vskip 0.1 in

{\bf Taking the inverse (mirror image).}
Let $\Delta$ be a \vk diagram over a symmetric
set of defining relations $\pp$.
Consider the mirror image $\Delta\iv$
of the graph $\Delta$.
Since $\pp$ is symmetric,
the graph $\Delta\iv$ is again a \vk diagram over $\pp$. We call $\Delta\iv$
the {\em inverse of $\Delta$}. 

{\bf Composition.} Let $\Delta_1$ and $\Delta_2$ be \vk diagrams over
a presentation $\pp$. Let $\partial(\Delta_1)=p_1p_1'$,
$\partial(\Delta_2)=p_2\iv p_2'$. Suppose that $\Lab(p_1)=\Lab(p_2)$ in the free
group. Here $\Lab(p)$ is the label of the path $p$. Then
by adding 0-cells we can make both paths $p_1$ and $p_2$
absolutely simple and the labels of $p_1$ and $p_2$ identical.
After that we can glue $\Delta_1$ and $\Delta_2$ by
identifying the corresponding edges of $p_1$ and $p_2$.
This operation will be called the {\em composition} of $\Delta_1$ and $\Delta_2$.

\vskip 0.1 in

Let $\Delta$ be any van Kampen diagram. The following definition is similar to
the definition of a dual graph of a diagram \cite{LS}, \cite{Ol89}.  Fix a point
in each of the cells in $\Delta$ and a point in the inside of each of the edges
of $\Delta$.
For each cell $\pi$ and each edge $e$ on the boundary of this cell,
fix a simple polygonal line $\ell(\pi,e)$ inside $\pi_i$ which connects the
point inside $\pi_i$ and the point inside $e$. We can choose these
lines in such
a way that $\ell(\pi, e)$ and $\ell(\pi, e')$ do not have common points except
for the fixed point inside $\pi$.

\vskip 0.1 in
The next definition of a band in a diagram is crucial for our paper as it was in 
\cite{SBR} and \cite{Ol1}.\vskip 0.1 in

Let $S$ be a set of letters and let $\Delta$ be a van Kampen diagram.  Fix 
pairs of $S$-edges on the boundaries of some cells from $\Delta$ (we assume that each of these
cells has exactly
two $S$-edges).

Suppose that
$\Delta$ contains a sequence of cells ($\pi_1,\ldots,\pi_n$)
such that for each $i=1,\ldots,n-1$ the cells $\pi_i$ and $\pi_{i+1}$ have a
common $S$-edge $e_i$ and this edge belongs to the pair of
$S$-edges fixed in $\pi_i$ and $\pi_{i+1}$.
Consider the line which is the union of the lines $\ell(\pi_i, e_i)$ and
$\ell(\pi_{i+1}, e_i)$, $i=1,\ldots,n-1$. This polygonal line
is called the {\em median} of this
sequence of cells. Then our sequence of cells ($\pi_1,\ldots,\pi_n$) with
common $S$-edges $e_1, \ldots, e_{n-1}$ is called an {\em $S$-band}
if the median is a absolutely simple curve or an absolutely simple closed
curve (\setcounter{band}{\value{ppp}} see Fig.  \theband).

\bigskip
\unitlength=0.90mm
\linethickness{0.4pt}
\begin{picture}(145.33,30.11)
\put(19.33,30.11){\line(1,0){67.00}}
\put(106.33,30.11){\line(1,0){36.00}}
\put(142.33,13.11){\line(-1,0){35.67}}
\put(86.33,13.11){\line(-1,0){67.00}}
\put(33.00,21.11){\line(1,0){50.00}}
\put(110.00,20.78){\line(1,0){19.33}}
\put(30.00,8.78){\vector(1,1){10.33}}
\put(52.33,8.44){\vector(0,1){10.00}}
\put(76.33,8.11){\vector(-1,1){10.33}}
\put(105.66,7.78){\vector(1,2){5.33}}
\put(132.66,8.11){\vector(-1,1){10.00}}
\put(19.33,30.11){\vector(0,-1){17.00}}
\put(46.66,29.78){\vector(0,-1){16.67}}
\put(73.66,30.11){\vector(0,-1){17.00}}
\put(142.33,30.11){\vector(0,-1){17.00}}
\put(117.00,30.11){\vector(0,-1){17.00}}
\put(16.00,21.11){\makebox(0,0)[cc]{$e$}}
\put(29.66,25.44){\makebox(0,0)[cc]{$\pi_1$}}
\put(60.33,25.44){\makebox(0,0)[cc]{$\pi_2$}}
\put(133.00,25.78){\makebox(0,0)[cc]{$\pi_n$}}
\put(145.33,21.44){\makebox(0,0)[cc]{$f$}}
\put(77.00,25.11){\makebox(0,0)[cc]{$e_2$}}
\put(122.00,25.44){\makebox(0,0)[cc]{$e_{n-1}$}}
\put(100.33,25.44){\makebox(0,0)[cc]{$S$}}
\put(96.33,30.11){\makebox(0,0)[cc]{$\dots$}}
\put(96.33,13.11){\makebox(0,0)[cc]{$\dots$}}
\put(26.66,4.78){\makebox(0,0)[cc]{$l(\pi_1,e_1)$}}
\put(52.66,4.11){\makebox(0,0)[cc]{$l(\pi_2,e_1)$}}
\put(78.33,4.11){\makebox(0,0)[cc]{$l(\pi_2,e_2)$}}
\put(104.66,3.78){\makebox(0,0)[cc]{$l(\pi_{n-1},e_{n-1})$}}
\put(139.00,3.11){\makebox(0,0)[cc]{$l(\pi_n,e_{n-1})$}}
\put(88.66,30.11){\makebox(0,0)[cc]{$p$}}
\put(88.66,13.11){\makebox(0,0)[cc]{$q$}}
\put(49.33,25.11){\makebox(0,0)[cc]{$e_1$}}
\put(33.33,21.11){\circle*{0.94}}
\put(46.66,21.11){\circle*{1.33}}
\put(60.33,21.11){\circle*{0.94}}
\put(129.66,20.78){\circle*{1.33}}
\put(73.66,21.11){\circle*{1.33}}
\put(117.00,20.78){\circle*{0.67}}
\put(95.00,20.56){\makebox(0,0)[cc]{...}}
\end{picture}

\begin{center}
\nopagebreak[4]
Fig. \theppp.
\end{center}
\addtocounter{ppp}{1}

Notice that an $S$-band can contain no cells. In this case the band is called {\em empty}.
An empty $S$-band consists of one $S$-edge and its median is empty.

We say that two bands {\em cross} if their medians cross. We say that a band is
an {\em annulus} if its median is a closed curve. In this case the first and the last cells of the band coincide 
\setcounter{annulus}{\value{ppp}}(see Fig. \theannulus a).

\bigskip
\begin{center}
\unitlength=1.5mm
\linethickness{0.4pt}
\begin{picture}(101.44,22.89)
\put(30.78,13.78){\oval(25.33,8.44)[]}
\put(30.78,13.89){\oval(34.67,18.00)[]}
\put(39.67,9.56){\line(0,-1){4.67}}
\put(25.67,9.56){\line(0,-1){4.67}}
\put(18.56,14.89){\line(-1,0){5.11}}
\put(32.56,7.11){\makebox(0,0)[cc]{$\pi_1=\pi_n$}}
\put(21.00,7.33){\makebox(0,0)[cc]{$\pi_2$}}
\put(15.89,11.78){\makebox(0,0)[cc]{$\pi_3$}}
\put(43.44,12.89){\line(1,0){4.67}}
\put(43.00,8.44){\makebox(0,0)[cc]{$\pi_{n-1}$}}
\put(19.22,10.89){\line(-1,-1){4.00}}
\put(25.89,20.44){\circle*{0.00}}
\put(30.78,20.44){\circle*{0.00}}
\put(30.33,1.56){\makebox(0,0)[cc]{a}}
\put(35.44,20.44){\circle*{0.00}}
\put(62.78,4.89){\line(0,1){4.67}}
\put(62.78,9.56){\line(1,0){38.67}}
\put(101.44,9.56){\line(0,-1){4.67}}
\put(101.44,4.89){\line(-1,0){38.67}}
\put(82.11,16.22){\oval(38.67,13.33)[t]}
\put(62.78,15.78){\line(0,-1){7.33}}
\put(101.44,16.44){\line(0,-1){8.22}}
\put(82.00,12.22){\oval(27.78,14.67)[t]}
\put(67.89,12.89){\line(0,-1){8.00}}
\put(95.89,13.11){\line(0,-1){8.22}}
\put(67.89,13.78){\line(-1,1){4.67}}
\put(95.67,14.00){\line(1,1){5.11}}
\put(73.67,9.56){\line(0,-1){4.67}}
\put(89.67,9.56){\line(0,-1){4.67}}
\put(76.56,21.33){\circle*{0.00}}
\put(81.44,21.33){\circle*{0.00}}
\put(86.11,21.33){\circle*{0.00}}
\put(65.22,12.22){\makebox(0,0)[cc]{$\pi_1$}}
\put(65.22,6.89){\makebox(0,0)[cc]{$\pi$}}
\put(71.00,6.89){\makebox(0,0)[cc]{$\gamma_1$}}
\put(92.78,7.11){\makebox(0,0)[cc]{$\gamma_m$}}
\put(98.56,7.11){\makebox(0,0)[cc]{$\pi'$}}
\put(98.33,12.67){\makebox(0,0)[cc]{$\pi_n$}}
\put(94.11,20.22){\makebox(0,0)[cc]{$S$}}
\put(85.89,6.89){\makebox(0,0)[cc]{$T$}}
\put(76.56,6.89){\circle*{0.00}}
\put(79.44,6.89){\circle*{0.00}}
\put(82.33,6.89){\circle*{0.00}}
\put(83.44,1.56){\makebox(0,0)[cc]{b}}
\end{picture}
\end{center}
\begin{center}
\nopagebreak[4]

Fig. \theppp.
\end{center}
\addtocounter{ppp}{1}

Let $\cal B$ be an $S$-band with common $S$-edges $e_1, e_2,\ldots, e_n$
which is not an annulus. Then the first cell
has an $S$-edge $e$ which forms a pair with $e_1$
and the last cell of $\cal B$ has an edge $f$ which forms a pair with
$e_n$. Then we shall say that $e$ is the {\em start edge} of $\cal B$ and
$f$ is the {\em end edge of $\cal B$}.  If $p$ is a path in $\Delta$ then
we shall say that a band
starts (ends) on the path $p$ if $e$ (resp. $f$) belongs to $p$.

Let $S$ and $T$ are two disjoint sets of letters, let ($\pi$, $\pi_1$, \ldots,
$\pi_n$, $\pi'$) be an $S$-band and let ($\pi$, $\gamma_1$, \ldots, $\gamma_m$,
$\pi'$) is a $T$-band. Suppose that:
\begin{itemize}
\item the medians of these bands form a simple closed
curve,
\item on the boundary  of $\pi$ and on the boundary of $\pi'$ the pairs of $S$-edges separate the pairs of $T$-edges,
\item the start and end edges of these bands are not contained in the
region bounded by the medians of the bands.
\end{itemize}
Then we say that these bands form an {\em $(S,T)$-annulus} and the curve formed by the medians of these bands is the {\em median} of this annulus (see Fig. \theannulus b). \setcounter{nonannulus}{\value{ppp}}
Fig. \thenonannulus\ shows that 
a multiple intersection of an $S$-band
and a $T$-band  does not necessarily produce an $(S,T)$-annulus.  

\begin{center}
\unitlength=1.50mm
\linethickness{0.4pt}
\begin{picture}(97.67,32.89)
\put(11.67,10.67){\framebox(24.00,3.56)[cc]{}}
\bezier{88}(35.67,14.22)(38.78,24.67)(42.33,14.22)
\bezier{200}(42.33,14.22)(45.22,-2.67)(15.22,10.67)
\bezier{168}(32.33,14.22)(39.22,32.89)(46.11,12.22)
\bezier{240}(46.11,12.22)(48.33,-6.89)(11.67,10.67)
\put(15.22,10.67){\line(0,1){3.56}}
\put(32.33,14.22){\line(0,-1){3.56}}
\put(13.44,12.22){\makebox(0,0)[cc]{$\pi$}}
\put(34.11,12.45){\makebox(0,0)[cc]{$\pi'$}}
\put(24.11,12.45){\makebox(0,0)[cc]{$\ldots$}}
\put(36.78,17.33){\line(-1,1){2.40}}
\put(43.44,16.00){\circle*{0.00}}
\put(44.33,12.67){\circle*{0.00}}
\put(43.89,9.11){\circle*{0.00}}
\put(18.56,14.22){\line(0,-1){3.56}}
\put(68.78,10.67){\framebox(20.44,3.56)[cc]{}}
\put(72.33,10.67){\line(0,1){3.56}}
\put(72.33,14.22){\line(0,0){0.00}}
\put(75.67,14.22){\line(0,-1){3.56}}
\put(85.00,14.22){\line(0,-1){3.56}}
\put(70.56,12.44){\makebox(0,0)[cc]{$\pi$}}
\put(87.22,12.22){\makebox(0,0)[cc]{$\pi'$}}
\bezier{80}(89.22,14.22)(92.33,22.00)(93.00,10.67)
\bezier{80}(93.00,10.67)(94.11,6.44)(78.78,6.44)
\bezier{100}(84.33,6.67)(65.00,4.89)(65.22,10.67)
\bezier{72}(65.22,10.89)(65.22,20.89)(68.78,14.22)
\bezier{168}(85.00,14.22)(95.00,32.00)(96.11,10.67)
\bezier{108}(96.11,10.67)(97.67,2.89)(78.78,2.89)
\bezier{128}(83.89,3.11)(62.11,0.89)(62.33,10.67)
\bezier{156}(62.33,10.67)(63.89,30.89)(72.33,14.22)
\put(90.11,16.22){\line(-1,1){2.67}}
\put(66.33,16.89){\line(1,1){2.89}}
\put(80.56,12.22){\makebox(0,0)[cc]{$\ldots$}}
\put(81.00,4.67){\makebox(0,0)[cc]{$\dots$}}
\end{picture}
\end{center}

\begin{center}
\nopagebreak[4]
Fig. \theppp.
\end{center}
\addtocounter{ppp}{1}

If  $\ell$ is the median of an $S$-annulus or an $(S,T)$-annulus then the maximal subdiagram of $\Delta$ contained in the region bounded by $\ell$ is called  the {\em inner diagram} of the annulus.

The union of cells of an $S$-band $\cal B$ forms a
subdiagram.
The boundary of this diagram, which we shall call the {\em boundary
of the band},
has the form $e^{\pm 1}pf^{\pm 1}q\iv$ (recall that we trace boundaries
of diagrams clockwise).
Then we say that $p$ is the
{\em top path} of $\cal B$, denoted by $\topp({\cal B})$,
and $q$ is the {\em bottom
path} of $\cal B$, denoted by $\bott({\cal B})$.

If ${\cal B}=(\pi_1,\ldots,\pi_n)$ is an $S$-band then
($\pi_n,\pi_{n-1},\ldots,\pi_1$)
is also an $S$-band which we shall call the {\em inverse of $\cal B$}.

We shall call an $S$-band {\em maximal} if it is not contained in any other
$S$-band. If an $S$-band $\ww$ starts on the contour of a cell $\pi$,
does not contain
$\pi$  and is not contained in any other $S$-band with these properties
then we call $\ww$ a {\em  maximal $S$-band starting on the contour
of $\pi$}.

\subsection{Known Facts About $G_N(\sss)$}

\label{known}

Now let us return to our presentation of the group $G_N(\sss)$.
Consider the following partition of the generating set $X$ of this group:

\begin{equation}
\label{part}
 \Theta_+ \cup Q_1\cup...\cup Q_k\cup \overline Y \cup\{\kappa_1\}\cup ...\cup\{\kappa_{4N}\}.
\end{equation}

Recall that $Q_i$ are the components of the vector of sets of state letters of the $S$-machine $\sss$ and $\overline Y$ is the set of tape letters of this $S$-machine.

The fact that these subsets are disjoint follows from the definition of $S$-machines.

Let $U$ be a block of this partition.
In order to consider $U$-bands and annuli, we need to
divide $U$-letters of some relation of $G_N(\sss)$ into pairs.

Every relation except for the hub contains a letter from
$\Theta_+$ and the inverse of this letter. These letters form a $\Theta$-pair.

Every auxiliary relation containing a $\overline Y$-letter,
contains a $\overline Y$-letter and the inverse of this letter. These
letters form a $\overline Y$-pair.

Every auxiliary relation containing $\kappa_i$ contains also $\kappa_i\iv$.
These two letters form a $\kappa_i$-pair.

If a transition relation has a letter from $Q_j$ then it has exactly two
letters from $Q_j^{\pm 1}$ (this follows from the definition of an $S$-rule).
They form a $Q_j$-pair.

For every $\kappa_j$-band the number $j$ is called the {\em index} of this band.

For convenience we repeat here some auxiliary results from \cite{SBR} which we shall 
use in this paper. Recall that the {\em area} of a diagram is the number of cells in this diagram.

\begin{lm}\label{cent} ({see \rm Lemma 7.1 in \cite{SBR}.}) If $\Delta$ contains a $\kappa$-annulus
then one can get a new diagram with the same boundary label as $\Delta$ by removing 
all cells of the annulus, identifying the outer and the inner boundaries of the annulus and adding no new cells. 
\end{lm}

\begin{lm} \label{gammapath} ({\rm Lemma 7.7 in \cite{SBR}.}) If $\Delta$ is reduced,
does not contain hubs and
$\partial(\Delta)$
has no $\Theta$-edges then $\Delta$ does not have cells.
\end{lm}

\begin{lm} \label{qr} ({\rm Lemma 7.4 in \cite{SBR}.}) If $\Delta$ is reduced and
does not contain hubs, then it does not have
$(\overline Q,\Theta)$-annuli.
\end{lm}

\begin{lm}\label{q} ({\rm Lemma 7.2 in \cite{SBR}.}) If $\Delta$ is reduced and does not have hubs then
it does not have $\overline Q$-annuli.
\end{lm}

\begin{lm} \label{r} ({\rm Lemma 7.6 in \cite{SBR}.}) If $\Delta$ is reduced and
does not contain hubs then it does not
contain
$\Theta$-annuli.
\end{lm}

\begin{df} \label{defdisc} A disc is any reduced van Kampen diagram over the presentation of $G_N(\sss)$  which has exactly one hub and no $\Theta$-edges on the contour.
\end{df}

\begin{lm} \label{prop13} ({\rm Proposition 10.1 in \cite{SBR}.}) 
The cyclically reduced boundary label of a disc is equal to $\kappa(w)$ for some admissible word $w$ accepted by the $S$-machine $\sss$. The area of the disc does not exceed big ``O" of the area of the corresponding computation. 
\end{lm}

\begin{lm} \label{disc} ({\rm see Lemma 11.1 in \cite{SBR}.}) If a diagram over the presentation of $G_N(\sss)$ contains discs then it has at least $4N-3$ $\kappa$-bands starting on the same hub and ending on the boundary of the diagram.
\end{lm}

\begin{lm} \label{bigons} ({\rm See the proof of Lemma 11.1 in \cite{SBR}.}) Let $\Delta$ be a reduced van Kampen diagram over the presentation of $G_N(\sss)$ containing two hubs connected by two $\kappa$-bands. Then the number of hubs in $\Delta$ can be reduced by 2 without changing the boundary label of $\Delta$.
\end{lm}

\begin{lm}\label{ar} ({\rm see Lemma 7.8 in \cite{SBR}.}) If $\Delta$ is reduced and
contains no $\overline Q$-edges then it does not
have $(\overline Y, \Theta)$-annuli.
\end{lm}

\subsection{Hyperbolic Graphs}
\label{hg}

Here we present some graph theoretic concepts and results from \cite{Ol1} and \cite{SBR}.
Let $\Gamma$ be a finite plane
2-complex with vertices $v_0,...,v_n$, $n\ge 1$. The vertex $v_0$ is called {\em external} and other vertices are called {\em internal}. We assume that there exists a curve from $v_0$ to the infinity without intersecting the edges of $\Gamma$ (that is $v_0$ is indeed an external vertex of the plane graph $\Gamma$).  We also assume that $\Gamma$ does not have faces with 1 vertex $v_i$ or bigons with vertices $v_i, v_j$, $i,j\ge 1$; there are also no 1-gons at the vertex $v_0$ 
but $v_0$ can be connected by several edges with an internal  vertex $v_i$. Such a graph $\Gamma$ will be called an $\ell$-graph if the degree of every internal vertex is $\ge \ell\ge 6$. The next lemmas show that for $\ell\ge 7$ $\ell$-graphs have {\em hyperbolic} properties. In particular the next Lemma shows that the ``curvature" of $\ell$-graphs is negative. 
It shows that the number of external edges of a hyperbolic graph is proportional to the number of internal edges.
  
\begin{lm} ({\rm See Lemma 3.2 in \cite{Ol1} or Lemma 11.5 in \cite{SBR}.}) In an $\ell$-graph, there exists a vertex $o$ 
of degree $d\ge \ell$ and $d-3$ consecutive (clockwise) edges $e_1,...,e_{d-3}$ connecting $o$ with the external vertex $v_0$, such that there are no vertices of $\Gamma$ between $e_i$ and $e_{i+1}$, $i=1,...,d-4$. 
\label{3.2}
\end{lm}

\begin{lm} ({\rm See Lemma 3.4 in \cite{Ol1}}.) Let $d_0, d_1,...,d_n$ be degrees of vertices $v_0, v_1,...,v_n$ of an $\ell$-graph $\Gamma$. Then $$d_0\ge 6-6n+\sum_{i=1}^n d_i.$$
\label{3.4}
\end{lm}

Let two edges $e$ and $e'$ go out of some internal vertex $v$ of an $\ell$-graph $\Gamma$. Then 
all edges going out of $v$ can be listed in the clockwise order as follows: $e, f_1,...,f_{n_1}e', f_1',...,f_{n_2}'$. We say that the edges $e$ and $e'$ form a {\em big angle}
if $n_1\le 2n_2$ and $n_2\le 2n_1$. 

We call any path $e_1,...,e_\alpha$ in $\Gamma$ {\em weakly curved} if the vertices  $(e_1)_+$,...,$(e_{\alpha-1})_+$ are internal and the edges $e_i\iv$ and $e_{i+1}$ form a big angle at the vertex $(e_i)_+=(e_{i+1})_-$ for every $i=1,...,\alpha-1$. 

\begin{lm} ({\rm See Lemma 3.5 in \cite{Ol1}.}). If a weakly curved path in an $\ell$-graph starts in an internal vertex then it cannot be closed (if $\ell\ge 6$; \setcounter{ovals}{\value{ppp}}
see Fig. \theovals a).
\label{3.5}
\end{lm}

This lemma implies that any maximal weakly curved path in an $\ell$-graph starts and ends in the external vertex.

Let $x$ be a Jordan arc on the plane which has exactly two common points with $\Gamma$: vertices $x_-$ and $x_+$. We shall call such an arc {\em $\Gamma$-normal}. Suppose that a simple path $p$ in $\Gamma$ passes through a vertex $v$ which is equal to $x_-$ or $x_+$. Then we can cyclically (clockwise) order all edges going out of $v$: $f_1,...,f_d$ so that $x$ comes to $v$ between $f_d$ and $f_1$. We say that the arc $x$ is $k$-{\em separated} from $p$ if the edges $f_1,...,f_k$, $f_{d-k+1},...,f_d$ and their inverses do not occur in $p$. 

Two maximal weakly curved paths $p$ and $p'$ in an $\ell$-graph $\Gamma$ are called {\em diverging} in their internal vertices $v$ and $v'$ if $v\ne v'$ and $v$ and $v'$ can be connected by a $\Gamma$-normal arc $x$ which is 2-separated from $p$ and $p'$.

\begin{lm} ({\rm See Lemma 3.6 in \cite{Ol1}.}) Any two diverging paths in an $\ell$-graph $\Gamma$ do not have common internal vertices (see Fig. \theovals a).
\label{3.6}
\end{lm} 

By Lemma \ref{3.5} every maximal weakly curved path $p$ in an $\ell$-graph $\Gamma$ is absolutely simple, so it divides the plane into two regions $\ooo_1$ and $\ooo_2$ (the inner and 
the outer regions). If for every vertex $v\ne v_0$ on $p$ the number of edges going from $v$ to $\ooo_1$ (resp. $\ooo_2$) is bigger by at least 8 than the number of edges going from $v$ to $\ooo_2$ (resp. $\ooo_1$) then the path $p$ is called an {\em oval} and we denote
$\ooo_1$ (resp. $\ooo_2$) by $\ooo(p)$. It is easy to understand that if $\ell\ge 28$ then for every internal edge $e$ of an $\ell$-graph there exists an oval passing through $e$.

We call a simple Jordan curve $\Gamma$-{\em transversal} if it does not pass through the vertices of $\Gamma$ and crosses the edges of $\Gamma$ transversally. We say that  a $\Gamma$-transversal curve {\em envelopes} an internal vertex $v$ of $\Gamma$ if it intersects consecutively edges $f_1,...,f_s$ going out of $v$ where $s-2$ is at least a half of the degree of $v$, and in the sectors between consecutive edges and the curve, there are no other vertices of $\Gamma$. 

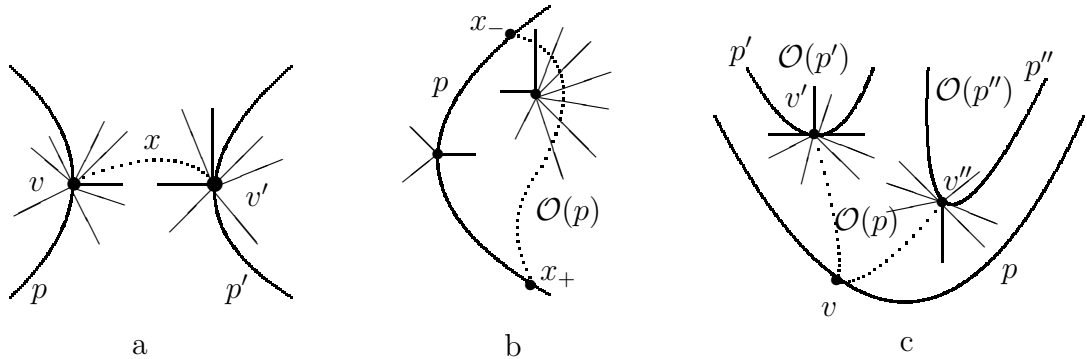
\begin{center}
\unitlength=1.50mm
\linethickness{0.8pt}
\begin{picture}(96.34,31.11)
\bezier{120}(1.23,25.78)(12.34,15.78)(1.23,5.33)
\bezier{136}(26.12,25.78)(12.56,13.78)(26.12,5.33)
\bezier{20}(6.78,15.55)(16.12,19.55)(19.45,15.33)
\put(6.78,15.11){\line(1,2){2.89}}
\put(7.00,15.33){\circle*{1.26}}
\put(19.45,15.33){\circle*{1.33}}
\put(7.45,15.78){\line(1,1){4.22}}
\put(7.45,15.33){\line(1,0){3.78}}
\put(7.23,14.89){\line(3,-1){4.22}}
\put(7.00,14.44){\line(1,-2){2.22}}
\put(6.78,15.55){\line(-1,2){2.00}}
\put(6.78,15.55){\line(-6,5){3.56}}
\put(6.78,15.11){\line(-2,-1){5.11}}
\put(19.23,16.00){\line(0,1){6.00}}
\put(19.23,15.33){\line(-2,5){2.67}}
\put(19.23,15.33){\line(-1,0){4.89}}
\put(19.45,15.11){\line(-1,-1){4.00}}
\put(19.67,15.55){\line(1,1){4.00}}
\put(19.45,15.11){\line(5,2){6.00}}
\put(19.23,14.89){\line(5,-6){4.00}}
\put(13.89,19.11){\makebox(0,0)[cc]{$x$}}
\put(3.67,15.33){\makebox(0,0)[cc]{$v$}}
\put(23.23,14.44){\makebox(0,0)[cc]{$v'$}}
\bezier{188}(49.00,5.55)(29.22,16.88)(49.00,31.11)
\put(45.67,28.66){\circle*{0.89}}
\put(47.45,6.44){\circle*{0.89}}
\bezier{25}(45.67,28.66)(53.22,26.66)(48.34,16.88)
\bezier{17}(48.34,16.88)(44.34,12.00)(47.67,6.44)
\put(47.89,23.33){\circle*{0.89}}
\put(47.89,23.33){\line(0,1){5.78}}
\put(47.89,23.11){\line(2,5){2.44}}
\put(47.89,23.11){\line(1,1){4.89}}
\put(48.11,23.11){\line(3,1){6.44}}
\put(47.67,23.11){\line(6,-1){6.67}}
\put(47.89,23.11){\line(1,-1){5.33}}
\put(47.89,23.11){\line(1,-3){2.44}}
\put(47.89,23.55){\line(-1,0){3.11}}
\put(48.11,23.33){\line(-1,-2){2.00}}
\put(39.22,18.00){\circle*{0.99}}
\put(39.22,18.00){\line(3,2){2.00}}
\put(39.00,18.00){\line(1,0){3.56}}
\put(39.22,18.00){\line(-1,1){2.44}}
\put(39.22,18.00){\line(-5,-4){3.11}}
\put(3.89,5.55){\makebox(0,0)[cc]{$p$}}
\put(21.45,6.00){\makebox(0,0)[cc]{$p'$}}
\put(12.78,0.89){\makebox(0,0)[cc]{a}}
\put(45.89,1.11){\makebox(0,0)[cc]{b}}
\put(43.67,29.33){\makebox(0,0)[cc]{$x_-$}}
\put(49.89,7.11){\makebox(0,0)[cc]{$x_+$}}
\put(50.78,12.88){\makebox(0,0)[cc]{$\ooo(p)$}}
\put(39.45,24.00){\makebox(0,0)[cc]{$p$}}
\bezier{292}(63.68,21.33)(81.01,-11.56)(96.34,21.33)
\bezier{104}(66.56,25.55)(73.23,13.78)(77.68,25.55)
\bezier{196}(82.56,25.55)(81.45,2.00)(93.01,24.66)
\put(74.56,6.89){\circle*{0.89}}
\put(72.56,19.78){\circle*{0.89}}
\put(83.90,13.78){\circle*{0.99}}
\bezier{17}(72.56,19.55)(75.23,9.78)(74.56,6.89)
\bezier{17}(83.90,13.78)(77.90,6.44)(74.79,6.66)
\put(72.56,19.78){\line(-1,0){4.00}}
\put(72.79,19.78){\line(-2,-1){4.22}}
\put(72.56,19.78){\line(-1,-4){1.11}}
\put(72.79,19.78){\line(1,0){4.22}}
\put(72.56,19.78){\line(5,-2){5.33}}
\put(72.56,19.78){\line(1,-1){4.00}}
\put(73.01,20.00){\line(1,3){1.11}}
\put(72.56,20.00){\line(0,1){4.00}}
\put(71.01,23.22){\makebox(0,0)[cc]{$v'$}}
\put(72.56,26.44){\makebox(0,0)[cc]{$\ooo(p')$}}
\put(66.12,27.11){\makebox(0,0)[cc]{$p'$}}
\put(73.90,4.22){\makebox(0,0)[cc]{$v$}}
\put(89.90,7.33){\makebox(0,0)[cc]{$p$}}
\put(77.23,12.00){\makebox(0,0)[cc]{$\ooo(p)$}}
\put(83.90,13.78){\line(5,-2){4.89}}
\put(83.90,13.55){\line(1,-1){3.78}}
\put(83.90,13.78){\line(0,-1){5.33}}
\put(83.90,13.55){\line(-1,1){3.78}}
\put(83.90,13.78){\line(-3,1){4.89}}
\put(83.45,13.78){\line(-5,-1){4.44}}
\put(83.90,13.78){\line(5,4){3.33}}
\put(85.34,15.89){\makebox(0,0)[cc]{$v''$}}
\put(86.79,23.55){\makebox(0,0)[cc]{$\ooo(p'')$}}
\put(92.56,25.78){\makebox(0,0)[cc]{$p''$}}
\put(80.79,1.11){\makebox(0,0)[cc]{c}}
\end{picture}
\end{center}
\begin{center}
\nopagebreak[4]
Fig. \theppp.
\end{center}
\addtocounter{ppp}{1}

\begin{lm} \label{3.7} ({\em See Lemma 3.7 in \cite{Ol1}.}). Let $\Gamma$ be an $\ell$-graph, and $p$ be an oval in $\Gamma$. Suppose that  a $\Gamma$-transversal curve $x$ is contained in $\ooo(p)$ except for the two points $x_-$ and $x_+$ which belong to $p$. Then either (1) $x$ intersects the same edge of $\Gamma$ twice 
and in the region bounded by this edge and $x$ there are no other vertices of $\Gamma$ or (2) $x$ envelopes an internal  vertex of $\Gamma$ (see Fig. \theovals b).
\end{lm}

Consider an arbitrary oval $p$ in an $\ell$-graph $\Gamma$. Let $v$ be an internal vertex on $p$. Let $v'$ be another internal vertex belonging to the closure of the region $\ooo(p)$ and let $p'$ be an oval passing through $v'$. If the ovals $p$ and $p'$ are diverging in $v$ and $v'$ and the corresponding $\Gamma$-normal arc $x$ is $3$-separated from $p'$ and $x$ is not in $\ooo(p')$ then we call the oval $p'$ 
{\em derived} from the oval $p$ in the vertices $v$ and $v'$. By Lemma \ref{3.6} $v'$ cannot belong to $p$.

\begin{lm} ({\rm See Lemma 3.8 in \cite{Ol1}.} (1) In the above notation, the region $\ooo(p)$ contains the region $\ooo(p')$ defined by the derived oval $p'$ .

(2) Let $p''$ be another oval derived  from $p$ in vertices $v$ and $v''$, $v'\ne v''$. Then the regions $\ooo(p')$ and $\ooo(p'')$ do not intersect (see Fig. \theovals c). 
\label{3.8}
\end{lm}

\section{Higman Embedding}
\label{he}

Now we are going to give the presentation of  a  finitely presented group $H_N(\sss)$ containing $G$ which will have isopermetric function $n^2T(n^2)^4$. 

This group will be constructed as a sequence of HNN extensions (similar to but simpler than) the well known Aanderaa construction in \cite{Ro}).

The set of {\em generators} of $H_N(\sss)$ 
consists of all the generators of $G_N(\sss)$, together 
with the following new letters: 
\, $\rho$, $d$, and the set $B$ which is in one-to-one correspondence
with the set $A$ (we assume that this correspondence preserves inverses). 

The set of {\em relations} of $H_N(\sss)$ consists of relations of the previously
given finite presentation of $G_N(\sss)$, together with the following new 
relations:

\medskip 

$\rho x = x \rho$ \ \ \ for all \  $x \in  
\{ \alpha, \delta, \omega\} \cup A \cup \{z_j \ |\ 0 \leq j \leq 4 \} \cup
\{ \kappa_i \ |\ 3 \leq i \leq 4N \}$; \ 

$\rho^{-1} \kappa_1 \rho = \kappa_1 d^{-1}$;

$\rho^{-1} \kappa_2 \rho = d \kappa_2$;

\medskip 

$d x = x d$ \ \ \ for all \ $x \in  
\{ \alpha, \delta, \omega\} \cup \{z_j \ | \ 0 \leq j \leq 4\}$;

$d^{-1} a_i d = a_ib_i$ \ \ \ for all $a_i \in A$ 
(where $b_i$ is the letter in $B$ that corresponds to $a_i$); 

\medskip 

$a_ib_j = b_ja_i$ \ \ \ for all $a_i \in A$, $b_j \in B$, 
that is every element from $A$ commutes with every element from $B$.


Let $G_b$ be the copy of $G$ over the generating set $B$ (recall that $B$ is in 
one-to-one correspondence with the original generating set $A$ of $G$). From 
now on will simply write $H$ for $H_N(\sss)$.
We will show later (see Lemma \ref{embedd} below) that the identity map on $B$ induces an isomorphism between $G_b$ and the subgroup 
$\langle B \rangle$ of $H$. 

{\bf Remark 3.} It can be easily deduced from lemmas from Section 4 that that $H$ is obtained from $G_N(\sss)$ by a sequence of 3 HNN extensions (in the first extension the stable letters are letters from $B$, in the second extension $d$ is the stable letter and in the third extension $\rho$ is the stable letter; the reader can easily recover from the presentation of $H$ the vertex groups and the edge groups of these extensions). Notice that if we used the version of $G_N(\sss)$ with linear word $\kappa(w)=\kappa_1w_1...\kappa_{4N}w_{4N}$ where $w_i$ are copies of $w$ written in disjoint copies of the alphabet $\overline Y\cup\overline Q$ (as we discussed in Remark 2), then in the definition of $H_N(\sss)$ we could eliminate the letter $d$
passing its responsibilities to $\rho$. That is $\rho$ would still commute with all letters which could appear in $w_i$, $i\ge 2$, and with all $k_i$, $i=1,...,4N$. It will also commute with copies of $z_0,...,z_4$, $\alpha, \delta, \omega$ in $w_1$, but for every $A$-letter $a$ which can appear in $w_1$ we would have $a^\rho=ab$.  This presentation would reduce the number of HNN extensions needed to built $H_N(\sss)$ to 2 and it would slightly simplify the proofs. The disadvantage is that in this case we would have to explain in details how to rewrite lemmas in \cite{SBR} using this new presentation of $G_N(\sss)$ and it would take too long. So we decided to use the above presentation of $H_N(\sss)$. 

\bigskip

The following three lemma shows that the subgroup generated by $B$ in $H$ is a homomorphic image of $G_b$ (the homomorphism is induced by the identity map on $B$).
 
\begin{lm}  \label{Main diagram}
The relations of $G_b$ follow from the presentation of $H_N(\sss)$ (i.e. the identity map on $B$ induces
a homomorphism from $G_b$ onto the subgroup $\langle B \rangle$ of $H_N(\sss)$).
\end{lm}

\noindent {\em Proof.} Let $u_b = 1$ be any relation in $G_b$.
We are going to describe a van Kampen diagram with boundary $u_b$. Let $u=u_a$ be the copy of the word $u_b$ rewritten in the alphabet $A$.

By  Proposition \ref{prop2} there exists a van Kampen diagram $\Delta$ with boundary label $\kappa(\sigma(c(u)))$.
This word is a product of  elements of the set $A\cup \{\alpha, \delta, \omega, z_0, z_1,z_2, z_3, z_4, \kappa_1,...,\kappa_{4N}\}$. For each element $x$ of this set there exists a relation of $H$ of the form 
$x^\rho=y$ (where $y$ is a word). Therefore there exists a $\rho$-annulus $R$ whose inner boundary is also bounded by $\kappa(\sigma(c(u)))$. Let us glue $\Delta$ inside this annulus. Denote the resulting van Kampen diagram by $\Delta_1$. Notice that all cells of $R$ are commutativity cells (corresponding to relations $x^\rho=x$) except for two cells corresponding to the relations 
$\kappa_1^\rho=\kappa_1d^{-1}$ and $\kappa_2^\rho=d\kappa_2$. 

The cyclic word written on the boundary of $\Delta_1$ has the form $Vd^{-1}\sigma(c(u))dW$ for some words 
$V$ and $W$ such that $V\sigma(c(u))W$ is a cyclic shift of $\kappa(\sigma(c(u)))$. The word $\sigma(c(u))$ is a product of words from $A\cup \{\alpha, \delta, \omega, z_0, z_1,z_2, z_3, z_4\}$. For each element $x$ of this set we have a relation $x^d=y$ where $y=x$ if $x\in \{ \alpha, \delta, \omega, z_0, z_1,z_2, z_3, z_4\}$
and $y=ab$ if $x=a\in A$. Since every element of $A$ commutes with every element of $B$, we can form an
annular diagram $\cal A$ with outer boundary labeled by the word $d^{-1}\sigma(c(u))d\sigma(c(u))^{-1}$ and the inner boundary labeled by $u_b$.  Let us glue the diagrams $\Delta_1$ and $\cal A$ along the path labeled by $d^{-1}\sigma(c(u))d$. The resulting annular diagram $\Delta_2$ will have outer boundary labeled by $\kappa(\sigma(c(u)))$ and the inner boundary labeled by $u_B$. 

Now draw diagram $\Delta_2$ on the top half of a sphere with the north pole inside the disc bounded by the inner boundary of 
$\Delta_2$ and the outer boundary  of $\Delta$ occupying the equator of the sphere. Since the outer boundary of $\Delta_2$ has the same label as the boundary of $\Delta$, we can draw $\Delta\iv$ (the mirror image of $\Delta$) on the 
bottom half of the sphere, and then glue two halves together. The stereographic projection of the resulting 
spherical graph on the plane gives us a van Kampen diagram over the presentation of $H$ 
with boundary label $u_b$.
The lemma is proved (\setcounter{figone}{\value{ppp}}see Fig. \thefigone) $\Box$

\unitlength=1 mm
\linethickness{0.6pt}
\begin{center}
\begin{picture}(68.50,74.83)
\put(38.00,36.50){\oval(36.67,42.33)[l]}
\put(50.00,36.50){\oval(24.00,42.33)[r]}
\put(33.33,57.67){\line(1,0){18.67}}
\put(31.00,15.33){\line(1,0){22.33}}
\put(30.00,57.67){\vector(1,0){5.00}}
\put(49.67,57.67){\vector(1,0){3.67}}
\put(29.33,57.67){\circle*{0.67}}
\put(47.33,57.67){\circle*{0.67}}
\put(29.33,57.67){\line(0,1){8.00}}
\put(29.33,65.67){\line(1,0){5.67}}
\put(35.00,65.67){\line(0,-1){7.67}}
\put(47.50,57.67){\line(0,1){8.00}}
\put(47.50,65.67){\line(1,0){5.83}}
\put(53.33,65.67){\line(0,-1){8.00}}
\put(29.33,61.33){\vector(0,1){3.83}}
\put(31.00,65.67){\vector(1,0){3.17}}
\put(35.00,58.33){\vector(0,1){2.17}}
\put(35.00,60.50){\vector(0,1){4.00}}
\put(48.50,65.67){\vector(1,0){4.00}}
\put(47.50,57.67){\vector(0,1){2.83}}
\put(47.50,60.50){\vector(0,1){4.00}}
\put(21.42,36.67){\oval(15.83,58.00)[l]}
\put(20.00,65.67){\line(1,0){9.33}}
\put(60.92,36.67){\oval(15.17,58.00)[r]}
\put(20.83,7.67){\line(1,0){41.33}}
\put(50.50,65.67){\line(1,0){11.33}}
\put(21.50,56.17){\vector(-3,4){6.00}}
\put(19.67,46.00){\vector(-1,0){6.17}}
\put(19.67,32.33){\vector(-1,0){6.17}}
\put(21.33,17.00){\vector(-3,-4){5.83}}
\put(29.33,15.33){\vector(0,-1){7.67}}
\put(35.00,15.33){\vector(0,-1){7.67}}
\put(47.50,15.33){\vector(0,-1){7.67}}
\put(53.33,15.33){\vector(0,-1){7.67}}
\put(60.00,16.67){\vector(1,-1){6.83}}
\put(62.00,32.33){\vector(1,0){6.50}}
\put(62.00,46.00){\vector(1,0){6.50}}
\put(60.33,56.00){\vector(3,4){6.17}}
\put(31.33,56.00){\makebox(0,0)[cc]{$^{\kappa_1}$}}
\put(50.33,56.00){\makebox(0,0)[cc]{$^{\kappa_2}$}}
\put(50.17,17.00){\makebox(0,0)[cc]{$^{\kappa_{2N+1}}$}}
\put(31.83,17.00){\makebox(0,0)[cc]{$^{\kappa_{2N+2}}$}}
\put(34.67,15.33){\vector(-1,0){4.67}}
\put(52.83,15.33){\vector(-1,0){4.50}}
\put(35.00,60.50){\line(1,0){12.50}}
\put(28.00,61.33){\makebox(0,0)[cc]{$^{\rho}$}}
\put(16.50,47.17){\makebox(0,0)[cc]{$^{\rho}$}}
\put(27.50,11.50){\makebox(0,0)[cc]{$^{\rho}$}}
\put(54.50,61.83){\makebox(0,0)[cc]{$^{\rho}$}}
\put(64.83,47.33){\makebox(0,0)[cc]{$^{\rho}$}}
\put(54.67,11.83){\makebox(0,0)[cc]{$^{\rho}$}}
\put(46.00,11.83){\makebox(0,0)[cc]{$^{\rho}$}}
\put(36.17,11.83){\makebox(0,0)[cc]{$^{\rho}$}}
\put(33.67,59.17){\makebox(0,0)[cc]{$^{\rho}$}}
\put(49.00,59.17){\makebox(0,0)[cc]{$^{\rho}$}}
\put(49.00,62.83){\makebox(0,0)[cc]{$^d$}}
\put(33.50,62.83){\makebox(0,0)[cc]{$^d$}}
\put(34.67,65.67){\line(1,0){4.33}}
\put(48.67,65.67){\line(-1,0){5.33}}
\put(38.83,60.50){\line(0,1){5.17}}
\put(43.33,60.50){\line(0,1){5.17}}
\put(37.67,63.17){\makebox(0,0)[cc]{$^d$}}
\put(44.50,63.17){\makebox(0,0)[cc]{$^d$}}
\put(41.00,62.33){\makebox(0,0)[cc]{$^u$}}
\bezier{28}(38.83,65.67)(38.33,69.67)(41.67,69.17)
\bezier{28}(38.83,65.67)(42.50,65.83)(41.67,69.17)
\bezier{72}(38.83,65.67)(34.17,74.17)(42.50,73.67)
\bezier{60}(42.33,73.67)(48.83,72.50)(43.33,65.67)
\put(40.00,35.83){\makebox(0,0)[cc]{Disc}}
\put(42.00,74.83){\makebox(0,0)[cc]{$^u$}}
\put(40.33,70.33){\makebox(0,0)[cc]{$^{u_b}$}}
\put(33.67,7.67){\vector(-1,0){3.00}}
\put(52.00,7.67){\vector(-1,0){3.67}}
\end{picture}
\end{center}

\begin{center}
\nopagebreak[4]
Fig. \theppp.
\end{center}
\addtocounter{ppp}{1}


\section{Forbidden Subdiagrams} \label{sect4}

In order to analyze the isoperimetric function of $H=H_N(\sss)$ we first use a new
presentation, which we will call {\bf the disc-based presentation}
of $H$. Although $H$ is finitely presented the
disc-based presentation will be infinite (this help prove important lemmas from Section 4).
The generators of the disc-based presentation are the same as those of
the original presentation. 
The relators of the disc-based presentation consist of the finite set of
original relators of $H$, together with the following new relators:

\medskip

$\kappa(w)=1$ for all admissible words $w$ accepted 
by $\sss$;

\smallskip

$u_b = 1$ for all $u_b$ such that $u_b = 1$ in $G_b$ and $u_b$ is a cyclically reduced word.

\medskip

Since by Proposition \ref{prop2} and Lemma \ref{Main diagram} these relations hold in $H$, the disc-based 
presentation is indeed a presentation of $H$. 
Notice that if a \vk diagram over the original finite presentation of $G_N(\sss)$ contains a disc (that is by Definition \ref{defdisc} a reduced diagram over the presentation of $G_N(\sss)$  which has exactly one hub and no $\Theta$-edges on the contour) then by Lemma \ref{prop13}
there exists a cell of the disc-based presentation of $H$ with boundary label equal in the free group to a cyclic shift of the boundary label of the disc. Therefore we can cut off the disc
from our diagram and replace it by a diagram with one cell.
All the lemmas in this section are about van Kampen diagrams 
relative to the disc-based presentation of $H$; and when we say 
``van Kampen diagram'', we
mean a van Kampen diagram over this presentation. 
Such a van Kampen diagram consists of {\it discs} with labels of the form
$\kappa(w)$ (= 1 in $G_N(\sss)$), $G_b$-cells with labels of the form $u_b$
(= 1 in $G_b$), and cells of the original finite presentation of $H$. 

\bigskip

With every van Kampen diagram over the disc-based presentation we associate a 
4-vector, called the {\em type} of this diagram. The first coordinate is the number of disc cells in this diagram, the second coordinate is the number of $(\rho,\kappa)$-cells, the third  coordinate is the total number of $\kappa$-, $\Theta$-, $(d,B)$-, $G_b$- and $\rho$-cells, excluding the discs
and $(\rho,\kappa)$-cells, and
the fourth coordinate is the number of $(A, B)$-cells (corresponding to the relations $a_ib_j=b_ja_i$). We order the types lexicographically. When we say that a diagram $\Delta$ is {\em minimal}, we always mean that the type of diagram is minimal among the types of all diagrams with the same boundary label.
We shall say that the discs have the highest rank among all cells of a diagram, the $(\rho,\kappa)$-cells have smaller rank, the cells counted in the third coordinate of the type have rank smaller than the rank of $(\rho,\kappa)$-cells, and the $(A,B)$-cells have the smallest rank. 
Notice that we can view \vk diagrams over the disc based presentation as graded diagrams \cite{Ol89}. Similar grading is used in \cite{Ol1} and \cite{SBR}. Notice also that by Proposition \ref{prop2} and Lemma \ref{prop13} we can add discs to the presentation of $G_N(\sss)$ without changing the group. This presentation will be called the {\em disc-based presentation} of $G_N(\sss)$.

The following Lemma shows that $G_N(\sss)$ is embedded into $H$.

\begin{lm} \label{gn} Let $\Delta$ be a minimal diagram over the disc based presentation of $H$ with boundary over the alphabet of generators of $G_N(\sss)$ then $\Delta$ is a diagram over the 
disc-based presentation of $G_N(\sss)$. 
\end{lm}

\noindent {\em Proof.} Indeed, if we kill letters $B$, $\rho$ and $d$ in any relation of the presentation of 
$H$ containing these letters, we get a trivial relation. Therefore if $\Delta$ contains a $B$-, $\rho$- or $d$-cell, then by replacing in $\Delta$
all $B$-, $\rho$-, $d$-edges with 0-edges, and then removing 0-cells, we can lower the type of the diagram. 
$\Box$ 

\medskip

Notice that the proof of Lemma \ref{gn} shows that $G_N(\sss)$ is embedded in $H$ and is a retract of $H$.

\medskip 

The rest of this section is similar to Section 2.4: we show that a minimal diagram over the disc-based presentation of $H$ does not have certain bands and annuli. 

In order to consider bands in diagrams over $H$ we need to choose pairs of letters in every relation 
over the disc-based presentation of $H$. We shall preserve the pairs in the relations over 
the 
presentation of $G_N(\sss)$. In every $\kappa$-cell except for the discs, there are two $\kappa$-edges, they form a pair. So we can consider $\kappa$-bands. In $(A, B)$-, $(\rho,A)$- and $(d,A)$-relations, there are two $\overline Y$-letters, they form a pair. Thus we can include $(B,A)$-, $(\rho,A)$- and $(d,A)$-cells in $\overline Y$-bands. Every $\overline Q$-relation except for the discs contains 
two $\overline Q$-letters, they form a pair. So we can include $(d,\overline Q)$- and $(\rho,\overline Q)$-cells in $\overline Q$-bands.
Every $B$-relation except $G_b$-relations and $(d,B)$-relations contain two $B$-letters, so we can consider $B$-bands (consisting of $(A,B)$-commutation cells). Every $\rho$-relation contain $\rho$ and $\rho\iv$, these letters form a pair.  So we can consider $\rho$-bands. Finally every $d$-relation except for $(d,\rho)$-relations contains $d$ and $d\iv$, these letters form a pair. so we can consider $d$-bands. 
 
In the rest of this section we shall often use the following obvious fact without references.

\begin{lm}\label{fact} Let $X$ denote $\kappa$, $\overline Y$, $\Theta$, $\overline Q$, $B$, $\rho$ or $d$. Then in every $X$-band in a reduced \vk diagram over the disc-based presentation of $H$, two neighbor cells can have only $X$-edges and $\overline Y$-edges in common. Therefore 
if an $X$-band does not have a $Z$-edge on its top or bottom path where $Z\ne X$ denotes $\kappa$,  $\Theta$, $\overline Q$, $B$, $\rho$ or $d$, then the $X$-band does not have $Z$-cells.
\end{lm} 

\begin{lm} \label{d}  A minimal diagram over the disc-based presentation of $H$ 
has no $d$-annuli.
\end{lm}
 
\noindent {\em Proof.} Suppose by contradiction that there is a $d$-annulus.
Consider any $d$-annulus $\Sigma$. The outer or the inner boundary $p$ of this annulus is a word 
over $\{ \alpha, \omega, \delta\} \cup A \cup \{ z_i \ |\  i = 0,
\ldots , 4\}$. Therefore the label of $p$ is a word over the generators of $G_N(\sss)$. 
By Lemma \ref{gn} the subdiagram $\Delta_1$ 
bounded by $p$ is over the disc-based presentation of $G_N(\sss)$. This immediately implies that $p$ is the inner boundary of $\Sigma$. 

If the subdiagram $\Delta_1$ 	contains discs then by Lemma \ref{disc} $\Delta_1$ contains a discs and $4N-3$ $\kappa$-bands starting on it which end on the boundary of $\Delta$. 
But a $d$-band cannot have $\kappa$-cells. Therefore $\Delta_1$ does not contain discs. Then by Lemma \ref{gammapath} $\Delta_1$ does not have cells. Therefore the label of $p$ is equal to 1 in the free group. This implies that $p$ contains two neighbor edges with mutually inverse labels. The corresponding neighboring cells of $\Sigma$ cancel which contradicts the minimality of $\Delta$. 
$\Box$
 
\begin{lm} \label{sigma-ann1}
In a minimal van Kampen diagram without $\kappa$-edges,
there are no $\rho$-annuli. 
\end{lm}
 
\noindent {\em Proof.} Suppose that a minimal diagram $\Delta$ without $\kappa$-edges contains a $\rho$-annulus $\Sigma$. Then $\Sigma$ consists of commutativity cells (corresponding to the relations $\rho x=x\rho$.  The inner and the outer boundaries of this annulus have the same labels. Therefore we can identify the inner and outer boundaries of $\Sigma$, removing all cells of $\Sigma$ from $\Delta$. This operation does not change the boundary label of $\Delta$, but decrease the type of the diagram. 
 $\Box$

\begin{lm} \label{sigma-kappa} 
In a reduced van Kampen diagram $\Delta$ over the disc-based presentation there are no $(\rho, \kappa)$-annuli, no $(\Theta, \kappa)$-annuli and no $(d,\kappa)$-annuli.
\end{lm}
 
\noindent {\em Proof.} By contradiction, assume $\Delta$ contains a 
$(\rho, \kappa_i)$-annulus, a $(\Theta, \kappa_i)$-annulus, or a $(d, \kappa_i)$-annulus for some $i$. 
Let $\Sigma$ be an innermost one. So $\Sigma$ is either a 
$(\rho, \kappa_i)$-annulus or a $(\Theta, \kappa_i)$-annulus or a $(d,\kappa_i)$-annulus, 
and it has no 
other $(\rho, \kappa_i)$-annuli, $(\Theta, \kappa_i)$-annuli, or $(d,\kappa)$-annuli in the subdiagram $\Delta_1$ bounded by the inner boundary of $\Sigma$.

Let us first consider the case where $\Sigma$ is a $(\Theta, \kappa_i)$-annulus.
Recall that the only cells which can occur in the $\kappa_i$-band are $\rho$-cells and $\Theta$-cells.
We claim that the $\kappa_i$-band of $\Sigma$ contains 
no $\rho$-cells. 
 
\bigskip
\begin{center}
\unitlength=1.50mm
\linethickness{0.4pt}
\begin{picture}(39.67,17.33)
\put(0.11,3.11){\line(1,0){38.22}}
\put(-1.00,1.78){\framebox(2.67,2.89)[cc]{}}
\put(37.22,1.78){\framebox(2.44,2.89)[cc]{}}
\put(19.33,10.22){\oval(38.44,14.22)[t]}
\put(-0.11,10.44){\line(0,-1){7.33}}
\put(38.56,10.67){\line(0,-1){7.56}}
\put(8.78,1.78){\framebox(2.89,2.89)[cc]{}}
\put(25.89,1.78){\framebox(2.89,2.89)[cc]{}}
\put(18.89,7.11){\oval(17.11,8.00)[t]}
\put(10.33,7.56){\line(0,-1){4.44}}
\put(27.44,7.33){\line(0,-1){4.22}}
\put(27.89,15.56){\makebox(0,0)[cc]{$\Theta$}}
\put(22.33,8.67){\makebox(0,0)[cc]{$\rho$}}
\put(19.89,1.11){\makebox(0,0)[cc]{$\kappa_i$}}
\end{picture}
\end{center}
\begin{center}
\nopagebreak[4]
\setcounter{innermost}{\value{ppp}}
Fig. \theppp.
\end{center}
\addtocounter{ppp}{1}

Indeed, if there were a $\rho$-cell on the $\kappa_i$-band, a $\rho$-band 
would start there, and it would form a $(\rho, \kappa_i)$-annulus with the 
$\kappa_i$-band of $\Sigma$ (since a $\rho$-band cannot cross a $\Theta$-band;  see Fig. \theinnermost).
But this 
would contradict the assumption that $\Sigma$ is innermost.

Hence, the $\kappa_i$-band of $\Sigma$ consists only of $\Theta$-cells. If it contains not only the intersection cells then we would get a $(\Theta, \kappa_i)$-annulus
inside $\Delta_1$ which again contradicts the assumption that $\Sigma$ is innermost.

It follows that the $\kappa_i$-band consists of the two intersection cells. Then these cells cancel: they have a common $\kappa_i$-edge, and their $\Theta$-edges have the same labels but point in the opposite directions (note that along one $\Theta$-band the labels of the 
$\Theta$-edges are the same everywhere; and the only $(\kappa_i, \Theta)$ cells 
are $(\kappa_i, \Theta)$-commutation cells).

Let us consider the case where $\Sigma$ is a $(\rho, \kappa_i)$-annulus.
The $\kappa_i$-band of $\Sigma$  contains no $\rho$-cells (except for the two intersection
cells between the $\rho$-band and the $\kappa_i$-band of $\Sigma$). The proof is quite similar
to the previous cases.

Indeed, if there were a $\rho$-cell on the $\kappa_i$-band, a $\rho$-band 
would start there, and it would form another $(\rho, \kappa_i)$-annulus with
the $\kappa_i$-band (since two $\rho$-bands cannot cross). 
But this would contradict the assumption that $\Delta$ is inner-most.
Similarly the $\kappa_i$-band cannot contain $\Theta$-cells because $\Theta$-bands and $\rho$-bands cannot intersect. Therefore the $\kappa_i$-band consists of two (adjacent) intersection $\rho$-cells. 
These two cells are adjacent and cancel (note that
for a given $i$ there exists only one kind of $(\rho, \kappa_i)$-cell). This contradicts the 
assumption that our diagram is reduced. 

Finally suppose that $\Sigma$ is a $(d,\kappa_i)$-annulus. The $d$-band of $\Sigma$ cannot contain $\kappa$-cells except for the intersection cells. A $\rho$-band in $\Delta_1$ starting on the boundary of the $\kappa_i$-band of $\Sigma$ cannot intersect the $d$-band. 
This as before implies that the $\kappa_i$-band of $\Sigma$ cannot contain $\rho$-cells except 
for the intersection cells. Similarly it cannot contain $\Theta$-cells (because a $\Theta$-band cannot intersect with a $d$-band). Therefore again the $\kappa_i$-band consists of two intersection cells that cancel. 
$\Box$

\begin{lm} \label{kappa-ann2} A reduced diagram over the disc-based presentation 
cannot contain $\kappa$-annuli.
\end{lm}

 \noindent {\em Proof.} Indeed, by Lemma \ref{sigma-kappa} a $\kappa$-annulus cannot contain $\Theta$-cells or $\rho$-cells, so it is empty. 
$\Box$

\medskip

\begin{lm} \label{B-ann} 
In a reduced diagram over the disc-based presentation there are no
$B$-annuli.
\end{lm}
 
\noindent {\em Proof.} A
$B$-annulus consists only of $(A, B)$- commutation cells, so
the labels of the inner and outer boundaries of the annulus are equal, and we can remove the annulus as in Lemma \ref{sigma-ann1}.
$\Box$
 
\begin{lm} \label{Q-sigma-ann1}
In a minimal diagram without $\kappa$-edges, there are no $(\overline Q, \rho)$-annuli, no
$(\overline Q, d)$-annuli and no $(\overline Q,\Theta)$-annuli.
\end{lm}
 
\noindent {\em Proof.} Let ${\Sigma}$ be an inner-most $(\overline Q, \rho)$-,
$(\overline Q, d)$-, or $(\overline Q, \Theta)$-annulus.  Let $\qq$ be the $\overline Q$-band and let $\bb$ be the $\rho$-, $d$, or $\Theta$-band in this annulus.

Suppose first that $\qq$ contains just two cells (the intersection cells with $\bb$).  Then these two $\rho$-, $d$-, or $\Theta$ cells have a common $\overline Q$-edge. Therefore these cells cancel. This contradicts the minimality of our diagram. 

Now suppose that $\qq$ contains more than two cells. The inner subdiagram of $\Sigma$  (the subdiagram bounded by the 
inner boundary of this annulus) cannot contain $\rho$-edges or $d$-edges by Lemmas \ref{d} and 
\ref{sigma-ann1}. 

Therefore the cells of  $\qq$ which are not the intersection cells with $\bb$ are $\Theta$-cells. A maximal $\Theta$-band $\bb'$  starting on the contour of $\qq$ cannot cross $\bb$ because there are no $(\Theta,d)$- or
$(\Theta, \rho)$-cells. Therefore we have a $(\overline Q, \Theta)$-annulus $\Sigma'$ contained inside the subdiagram bounded by the annulus $\rr$. This contradicts the assumption that $\Sigma$ is innermost.  
$\Box$
 
\begin{lm} \label{qann} In a minimal diagram without discs there are no $\overline Q$-annuli.
\end{lm}

\noindent {\em Proof.} Suppose that $\Sigma$ is a $\overline Q$-annulus in a minimal diagram $\Delta$.  Let $\Delta'$ be the subdiagram bounded by the external boundary of $\Sigma$. 
By Lemma \ref{Q-sigma-ann1}
$\Sigma$ does not contain $d$-cells or $\rho$-cells. Therefore the label of the boundary of $\Delta'$ is a word over the generators of $G_N(\sss)$. Thus $\Delta'$ is a diagram over the presentation of $G_N(\sss)$. But by Lemma \ref{q} a diagram over the presentation of $G_N(\sss)$
cannot contain $\overline Q$-annuli. $\Box$

\begin{lm}\label{theta} In a minimal diagram without discs there are no $\Theta$-annuli. 
\end{lm}

\noindent {\em Proof.} Let $\Sigma$ be a $\Theta$-annulus in a minimal diagram $\Delta$. Let $\Delta'$ be the subdiagram bounded by the outer boundary of $\Sigma$. Then the boundary label of $\Delta'$ is a word over the generators of $G_N(\sss)$, so $\Delta'$ is a diagram over the presentation of $G_N(\sss)$, and it remains to use Lemma \ref{r}
$\Box$

\begin{df} \label{def2} {\rm 
Let $\Phi$ be a $\Theta$- or $\rho$-band, and let $\Pi$ be a disc.
We say that $\Phi$ {\bf touches} $\Delta$ if the bottom path $\bott(\Phi)$ has at least 2 common $\kappa$-edges with the
contour of $\Pi$.  Let $\kappa_i$ (resp. $\kappa_j$) be the label of the first (resp. last) common $\kappa$-edge of $\Phi$ and $\partial(\Pi)$ (counted along the path $\bott(\Phi)$. 
Let $p_1$ be the shortest subpath of $\bott(\Phi)$ which contains these $\kappa$-edges and let $p_2$ be the subpath of the contour of $\Pi$ with the same ends as $p_1$ such that the subdiagram bounded by $p_1$ and $p_2$ does not contain the disc $\Pi$. Then the subdiagram bounded by $p_1$ and $p_2$ will be denoted by $\Delta(\Phi,\Pi)$ (\setcounter{touch}{\value{ppp}}
see Fig. \thetouch). }
\end{df}

\begin{center}
\unitlength=1.50mm
\linethickness{0.4pt}
\begin{picture}(54.78,33.11)
\put(29.33,5.00){\oval(33.11,10.44)[]}
\put(13.22,7.33){\line(-1,1){3.11}}
\put(10.11,10.44){\line(1,1){3.33}}
\put(13.44,13.78){\line(1,-1){3.56}}
\put(43.22,9.78){\line(1,1){3.78}}
\put(47.00,13.56){\line(1,-1){3.56}}
\put(50.56,10.00){\line(-1,-1){4.89}}
\put(10.11,10.22){\line(-1,-1){4.00}}
\put(13.00,7.56){\line(-1,-1){3.78}}
\put(46.56,6.00){\line(1,-1){3.56}}
\put(50.33,10.00){\line(1,-1){4.44}}
\bezier{204}(13.44,13.56)(30.11,33.11)(46.78,13.56)
\bezier{172}(16.78,10.22)(30.78,26.67)(43.22,9.78)
\put(30.33,20.67){\makebox(0,0)[cc]{$\Phi$}}
\put(30.33,4.00){\makebox(0,0)[cc]{$\Pi$}}
\put(29.22,13.56){\makebox(0,0)[cc]{$\Delta(\Phi,\Pi)$}}
\put(21.89,8.00){\makebox(0,0)[cc]{$p_2$}}
\put(35.22,14.89){\makebox(0,0)[cc]{$p_1$}}
\put(15.44,7.56){\makebox(0,0)[cc]{$\kappa_i$}}
\put(43.44,6.67){\makebox(0,0)[cc]{$\kappa_j$}}
\put(11.67,8.67){\line(0,1){2.89}}
\put(12.33,8.22){\line(0,1){4.00}}
\put(13.44,8.44){\line(0,1){4.67}}
\put(14.33,9.33){\line(0,1){2.89}}
\put(15.22,10.00){\line(0,1){1.33}}
\put(45.00,8.67){\line(0,1){2.22}}
\put(45.89,7.56){\line(0,1){4.44}}
\put(46.78,6.89){\line(0,1){5.78}}
\put(47.67,8.00){\line(0,1){4.22}}
\put(48.56,8.67){\line(0,1){2.67}}
\put(49.22,9.33){\line(0,1){1.33}}
\put(13.00,8.22){\line(0,1){4.44}}
\end{picture}
\end{center}
\begin{center}
\nopagebreak[4]
Fig. \theppp.
\end{center}
\addtocounter{ppp}{1}

\begin{lm} \label{emptydiagram} Suppose that a $\Theta$-band $\Phi$ touches a disc $\Pi$ in a minimal diagram $\Delta$ and $\Delta(\Phi,\Pi)$ does not contain discs. Then 
the subdiagram $\Delta(\Phi,\Pi)$ has no cells.
\end{lm}

 \noindent {\em Proof.} Let us use the notation from the Definition \ref{def2}.  Every maximal $\kappa$-band in $\Delta(\Phi,\Pi)$ starting on $p_2$ must end on $p_1$. By Lemma \ref{sigma-kappa} none of these $\kappa$-bands contain $\rho$-cells or $\Theta$-cells (otherwise there would be $(\rho, \kappa)$- or $(\Theta, \kappa)$-annuli. Every maximal $\kappa$-band starting on $p_1$ must end on $p_2$ because otherwise we would have a $(\Theta,\kappa)$-annulus. 
Therefore each of the $\kappa$-bands starting on the boundary of $\Delta(\Phi,\Pi)$ is empty. Therefore $\Delta(\Phi,\Pi)$ is a disjoint union of several subdiagrams $\Delta_i$, $i=1,...,k$, connected by 
$\kappa$-edges (bridges): for every $i=1,...,k-1$, a $\kappa$-edge connects $\Delta_i$ and $\Delta_{i+1}$.  The boundary of each $\Delta_i$ contains no $\kappa$-edges.

Consider one of the subdiagrams $\Delta_i$. The label of the boundary of this diagram is a word over the generators of $G_N(\sss)$. Therefore by Lemma \ref{gn} $\Delta_i$ is a diagram over $G_N(\sss)$. The contour of $\Delta_i$ does not contain $\Theta$-edges, and $\Delta_i$ does not contain discs by the assumption, so by Lemma \ref{gammapath} $\Delta_i$ does not contain cells. $\Box$

\medskip

The next construction plays a very important role both in \cite{SBR} and in \cite{Ol1}. We repeat it here.

\begin{df} \label{movingTheta} {\rm {\bf Moving $\Theta$-bands.} Suppose that in a minimal diagram $\Delta$ a $\Theta$-band $\rr$ touches a disc $\Pi$. By Lemma \ref{emptydiagram} then the bottom path of $\rr$ has a common subpath with the contour of $\Pi$ starting and ending with $\kappa$-edges. Let $p$ be the 
maximal common subpath with this property, so that $\bott(\Theta)=qpq'$, $\partial(\Pi)= p p_1$.
Without loss of generality we can assume that the label $\Lab(p)$ of the path $p$ has the form $\kappa_iw^{\pm 1}\kappa_{i+1}w^{\mp 1}\kappa_{i+2}..\kappa_j$ where $w$ is an admissible word accepted by the $S$-machine $\sss$ (it could be also the inverse of this word, but then we could take the mirror image of $\Delta$ instead of $\Delta$).  Then for some word $V$ we have that $\Lab(p)V$ is a cyclic shift of $\kappa(w)=\Lab(\partial(\Pi))$. There exists a $\Theta$-band $\rr'$ with the bottom path labeled by the word $V$ and the $\Theta$-edges having the same labels as in $\rr$.
Let $\rr''$ be the subband of $\rr$ with bottom path $p$, so $\rr=\rr_1\rr''\rr_2$. Let $e$ be the start edge of $\rr''$ and let $e'$ be the end edge of $\rr''$.
Cut the diagram $\Delta$ along the path $e\iv p_1 e'$. We can fill the resulting hole by gluing 
in the $\Theta$-band $\rr'$ and the mirror image $\overline\rr'$ of $\rr'$. The new diagram $\Delta'$ that we obtain this way will have two $\Theta$-bands instead of the old $\Theta$-band $\rr$. The first is $\rr_1(\overline\rr')\iv\rr_2$ (the inverse band $(\overline\rr')\iv$ differs from $\overline\rr'$ by the order of  cells) and the second one is $\rr''\rr'$. The second $\Theta$-band is an annulus which touches $\Pi$ along its inner boundary. If we replace the disc $\Pi$ be the corresponding \vk diagram over the presentation of $G_N(\sss)$, we see that the subdiagram $\Pi'$ bounded by the outer boundary of the annulus $\rr''\rr'$ is a diagram over the presentation of $G_N(\sss)$ with exactly one hub and no $\Theta$-edges on the boundary. By Definition \ref{defdisc} $\Pi'$ is a disc. We replace it by one cell of the disc-based presentation. Then we reduce the resulting diagram.
Notice that this construction amounts to enlarging the disc $\Pi$ and  moving the band $\rr$ through the disc (the band $\rr_1(\overline\rr')\iv\rr_2$ touches the disc $\Pi'$ along the $\kappa_{j+1}$-,
...,$\kappa_{i-1}$-edges) and then reducing the resulting diagram (\setcounter{theta}{\value{ppp}}see Fig. \thetheta). 
}
\end{df}
\begin{center}
\unitlength=1.50mm
\linethickness{0.6pt}
\begin{picture}(84.50,56.00)
\put(7.33,33.08){\oval(9.33,15.83)[l]}
\put(20.16,33.08){\oval(15.33,15.83)[r]}
\put(7.33,41.00){\line(1,0){13.33}}
\bezier{64}(9.33,51.67)(7.66,43.50)(14.83,41.00)
\bezier{52}(10.66,51.67)(9.83,44.67)(15.83,42.50)
\bezier{60}(16.00,42.50)(25.16,43.33)(28.66,39.17)
\bezier{72}(28.83,39.17)(32.66,30.67)(26.16,24.67)
\bezier{56}(26.16,24.50)(24.33,22.50)(13.00,23.00)
\bezier{60}(13.00,22.83)(5.50,21.17)(5.83,13.83)
\bezier{72}(3.66,13.83)(3.33,22.67)(12.00,25.17)
\put(3.66,13.83){\line(1,0){2.17}}
\put(9.33,51.67){\line(1,0){1.33}}
\put(35.33,32.17){\vector(1,0){6.67}}
\put(59.00,32.58){\oval(9.33,15.83)[l]}
\put(71.83,32.58){\oval(15.33,15.83)[r]}
\put(59.00,40.50){\line(1,0){13.33}}
\bezier{60}(67.66,42.00)(76.83,42.83)(80.33,38.66)
\bezier{72}(80.50,38.66)(84.33,30.16)(77.83,24.16)
\bezier{56}(77.83,24.00)(76.00,22.00)(64.66,22.50)
\bezier{84}(67.16,44.16)(68.66,56.00)(77.33,54.50)
\bezier{84}(67.66,42.16)(69.83,54.83)(77.50,52.83)
\bezier{80}(63.66,20.33)(65.33,9.50)(74.33,9.83)
\bezier{76}(64.66,22.33)(65.66,12.50)(74.33,11.66)
\bezier{92}(67.66,42.33)(52.50,44.33)(52.16,36.16)
\bezier{80}(52.16,36.16)(50.66,23.00)(57.16,22.16)
\bezier{32}(57.00,22.16)(63.16,22.16)(64.66,22.33)
\bezier{132}(67.16,44.16)(48.66,47.16)(49.83,32.50)
\bezier{116}(50.00,32.50)(49.33,17.66)(63.66,20.33)
\put(64.66,32.33){\makebox(0,0)[cc]{Disc}}
\put(15.50,32.67){\makebox(0,0)[cc]{Disc}}
\put(25.16,33.00){\makebox(0,0)[cc]{$p$}}
\put(4.16,33.67){\makebox(0,0)[cc]{$p_1$}}
\put(56.16,32.83){\makebox(0,0)[cc]{$p_1$}}
\put(76.83,32.50){\makebox(0,0)[cc]{$p$}}
\put(6.83,46.00){\makebox(0,0)[cc]{$q$}}
\put(2.83,20.00){\makebox(0,0)[cc]{$q'$}}
\put(13.16,48.00){\makebox(0,0)[cc]{$\rr_1$}}
\put(32.83,33.33){\makebox(0,0)[cc]{$\rr''$}}
\put(8.83,16.33){\makebox(0,0)[cc]{$\rr_2$}}
\put(73.50,49.16){\makebox(0,0)[cc]{$\rr_1$}}
\put(84.50,32.16){\makebox(0,0)[cc]{$\rr''$}}
\put(70.83,16.16){\makebox(0,0)[cc]{$\rr_2$}}
\put(46.33,31.83){\makebox(0,0)[cc]{$\overline\rr'$}}
\put(14.33,41.00){\line(0,1){2.33}}
\put(7.33,25.17){\line(1,0){15.00}}
\put(12.16,25.17){\line(0,-1){2.50}}
\put(59.16,24.66){\line(1,0){13.33}}
\put(67.66,40.50){\line(0,1){1.67}}
\put(64.66,22.33){\line(0,1){2.33}}
\put(64.66,22.33){\rule{0.17\unitlength}{2.33\unitlength}}
\put(74.33,9.66){\rule{0.20\unitlength}{2.00\unitlength}}
\put(64.66,22.33){\line(-1,-2){1.20}}
\put(67.67,42.16){\line(-1,2){1.20}}
\put(67.83,42.16){\line(-1,2){1.20}}
\put(67.67,42.16){\line(-1,2){1.20}}
\put(77.50,52.66){\line(0,1){2.20}}
\end{picture}

\nopagebreak[4]
Fig. \theppp.
\addtocounter{ppp}{1}
\end{center}

This construction has several important applications. For example, suppose that the band $\rr$ had more than $2N$ common $\kappa$-edges with the contour of the disc $\Pi$ and there are no other discs between $\rr$ and $\Pi$. Then the $\Theta$-band moving construction reduces the number of $\Theta$-cells in the diagram without increasing the number of cells of bigger ranks. Indeed, every maximal $(\overline Y, \overline Q)$-part of the boundary of a disc has the same length by Lemma \ref{prop13}. Therefore the $\Theta$-band moving construction reduces the number 
of $\Theta$-cells. On the other hand the $\Theta$-band moving construction replaces one disc $\Pi$ by another disc $\Pi'$ which does not increase the number of discs in the diagram, and it does not touch $(\kappa,\rho)$-cells and cells of the same rank as $\Theta$-cells). 
This implies the following statement. 

\begin{lm} \label{Theta-touch} In a minimal diagram there is no $\Theta$-band $\rr$ which touches a disc $\Pi$ in more than $2N$ edges and such that $\Delta(\rr,\Pi)$ contains no discs. 
\end{lm}

We will be using this construction more later.

Now we are going to define an analogous construction for $\rho$-bands. It is similar to a construction in \cite{Ol1}.

\begin{lm} \label{lm1}
Let $W=\sigma(c(u))$ where $u$ is a word in the alphabet $A$ such that $u=1$ in the group $G$. Then the equalities $W=\rho W\rho\iv$ and $W=dWd\iv$ can be deduced from the presentation of $H_N(\sss)$
without using $\kappa$-relations.
\end{lm}

\noindent {\em Proof.} By Proposition \ref{prop1} (4) $W\equiv z_0\alpha^n z_1uz_2\delta^nz_3\omega^nz_4$
and $u_b$ is equal to 1 in the group $G_b$.
Therefore using the non-$\kappa$-relations of $H_N(\sss)$ we can transform the word $dWd\iv$ into $z_0\alpha^n z_1uu_bz_2\delta^nz_3\omega^nz_4$. Now the subword $u_b$ can be deleted 
by using the $G_b$-relation $u_b=1$. 

The first statement of the lemma follows from $\rho$-relations of $H_N(\sss)$.
$\Box$

\begin{lm}\label{lm2} Suppose that a $\rho$-band $\Sigma$ touches a disc $\Pi$ in a minimal diagram $\Delta$ along two consecutive $\kappa$-edges whose indices differ from 1 and 2. Suppose further that the subdiagram $\Delta(\Sigma, \Pi)$ contains no discs. 
Then the label of $\partial(\Pi)$ has the form $\kappa(\sigma(c(u)))$ where $u=1$ in $G$. 
\end{lm}

\noindent {\em Proof.} We use the notation of Definition \ref{def2}.  Since $p_1$ and $p_2$ do not contain $\rho$- or $\Theta$-edges, by Lemmas \ref{sigma-ann1} and \ref{theta} the subdiagram $\Delta(\Sigma, \Pi)$ contains no $\rho$- or $\Theta$-cells and also (by assumption) contains no discs. Since the subdiagram $\Delta(\Sigma,\Pi)$ contains no discs, it contains no $\kappa$-edges except the two $\kappa$-edges where $\bott(\Sigma)$ touches $\Pi$ (these $\kappa$-edges are common for $p_1$ and $p_2$). Since the indices of these $\kappa$-edges differ from 1 and 2, we can conclude that the paths $p_1$ and $p_2$ contain no $d$-edges. Therefore by Lemma \ref{d} $\Delta(\Sigma,\Pi)$ contains no $d$-edges. 

Therefore this subdiagram contains 
no $\overline Q$-cells (because every $\overline Q$-cell is either a disc or a $\Theta$-cell, or a $d$-cell or a $\rho$-cell). Therefore each $\overline Q$-edge of $p_2$ also belongs to $p_1$. By Lemma \ref{Q-sigma-ann1}, every $\overline Q$-edge of $p_1$ belongs to $p_2$. Since every $(\rho,\overline Q)$-relation has the form $\rho z_i=z_i\rho$ where $i=0,1,2,3,4$, we can conclude that the word $V$ obtained from the label of $p_2$ by removing all non-$\overline Q$-letters, is a product of $z_i$, $i\in \{1,2,3,4\}$. 
In the case where $\Delta$ has only one hub, $\Delta$ is a disc
(see Definition \ref{defdisc} and Lemma \ref{prop13}).
By the definition of $S$-machines, the label of $p_2$ contains one representative from each component $Q_i$ of the vector $Q$, and these representatives are placed in $w$ in the order of indices of $Q_i$. Since words $z_j$ contain no common letters and $z_0z_1z_2z_3z_4$ contains representatives of all $Q_i$, each $Q_i$ represented precisely once, we can conclude that $V\equiv z_0\dots z_4$. Now the statement of our lemma follows from Proposition \ref{prop1} (part 3) and Proposition \ref{prop2} (part 1).
$\Box$

\begin{lm} \label{moverho} Let $\Delta$ be a van Kampen diagram over the disc-based presentation. Suppose that $\Delta$ consists of a disc $\Pi$ with boundary label $\kappa(\sigma(c(u)))$ for some word $u$ equal to 1 in $G$ and $\rho$-cells $\rho_1$,...,$\rho_n$ attached to $\Pi$ along $n\le 4N$ $\kappa$-edges $e_i$, $i\in \{1,...,n\}$
and attached to $\Pi$ in a similar manner (that is their $\rho$-edges either all point toward $\Pi$ or all point away from $\Pi$).
Let $S$ be the set of the indices of the $\kappa$-edges $e_i$,  $i=1,...,n$. Then there exists another diagram $\Delta'$ 
with the same boundary label as $\Delta$, consisting of a disc with the same boundary label as $\Pi$, $4N-n$ $(\rho, \kappa)$-cells $\pi_{i+1},...,\pi_{4N}$ attached to it and a number of cells of smaller ranks such that the indices of $\kappa$-edges of the cells $\pi_{n+1},...,\pi_{4N}$ form a complement to the set $S$ in the set $\{1,...,4N\}$.
\end{lm}

\noindent {\em Proof.} Let $W=\sigma(c(u))$. 
Consider an auxiliary diagram $\Delta_1$ obtained by attaching to $\Delta$ new $(\rho,\kappa)$-cells $\pi_{n+1}$, ...,$\pi_{4N}$ in the same manner as the other cells $\pi_i$, $i\le n$ (that is their $\rho$-edges point toward $\Pi$ or away from $\Pi$, as in $\pi_1, \pi_2,...,\pi_n$). 
The boundary of the diagram $\Delta_2$ has label $\prod\kappa_i U_i$ where $U_1\equiv (\rho d)^{\pm 1}W(\rho d)^{\mp 1}$ (depending on the orientation of $\rho$-edges), 
$U_j\equiv \rho^{\pm 1} W^{\pm 1}\rho^{\mp 1}$ for $j>1$. 

By Lemma \ref{lm1}, for every $i=1,...,4N$, to the paths labeled by $U_i$, we can attach a number of cells of ranks smaller than the rank of $\kappa$-cells, and obtain a diagram $\Delta_3$ whose boundary label is the same as the label of the disc $\Pi$. 

Thus the diagram $\Delta_3$ can be obtained by gluing $\Delta_1$ and a diagram $\Delta_0$ consisting of $4N-n$ $(\rho, \kappa)$-cells whose indices form a complement to $S$, and several cells of smaller ranks. This implies that, conversely, there exists a diagram $\Delta'$ with the same boundary label as $\Delta_1$, which can be obtained by gluing $\Pi$ with the mirror image $\Delta_0\iv$ of the diagram $\Delta_0$. $\Box$

The construction defined in the proof of Lemma \ref{moverho} will be called the {\em $\rho$-band
moving construction.}

The next Lemma gives the first application of this construction.

\begin{lm}\label{lm3} Suppose that a $\rho$-band $\Sigma$ touches a disc $\Pi$ along $n$
$\kappa$-edges in a minimal diagram $\Delta$. Then $n<4N-6$.
\end{lm}

\noindent {\em Proof.} By contradiction, suppose that $\Sigma$ touches $\Pi$ along $n\ge 4N-6$ $\kappa$-edges. Since $N\ge 6$, two of these edges have consecutive indices (mod $4N$) which differ from 1 and 2. By Lemma \ref{lm1} the boundary label of $\Pi$ is $\kappa(\sigma(c(u)))$ for some word $u$ which is equal to 1 in $G$. 

Consider a subdiagram $\Delta_1$ of $\Delta$ consisting of the disc $\Pi$ with $n\ge 4N-6$ $(\rho,\kappa)$-cells $\pi_i$ from $\Sigma$ attached to it in the same manner (all $\rho$-edges of these cells are pointing either toward the disc or away from it). 
By Lemma \ref{moverho} we can find a diagram $\Delta'$ which consists of the disc $\Pi$, $4N-n$ 
$(\rho, \kappa)$-cells and a number of cells of smaller ranks, and which has the same boundary label as $\Delta_1$. Since $4N-n<n$, the type of the diagram $\Delta'$ is smaller than the type of $\Delta_1$. Replacing subdiagram $\Delta_1$ in $\Delta$ with $\Delta'$ will make the type
of the whole diagram smaller which contradicts the minimality of $\Delta$. $\Box$

\bigskip

\bigskip

With a van Kampen diagram $\Delta$ over the disc-based presentation we can associate the following graph $\Gamma(\Delta)$, 
which we call the {\bf disc graph}.
The internal vertices are the discs, the external vertex is taken outside the diagram, the edges are the $\kappa$-bands with indices $\ge 4$ that
start on discs, any $\kappa$-band with index $\ge 4$ that starts on a disc and ends on the boundary will be considered as an edge connecting this disc with the external vertex.  
This graph is similar to the graphs considered in \cite{SBR} and \cite{Ol1} (see Section \ref{hg}).

\begin{df}{\rm A diagram $\Delta$ over the disc-based presentation of $H$ is called {\em normal}
if it does not contain two discs $\Pi_1$ and $\Pi_2$ connected by two $\kappa$-bands $\bb_1$, and $\bb_2$ with consecutive (mod $4N$) indices not equal 1 and 2, and the region between 
$\bb_1$ and $\bb_2$ does not contain discs.  
}\end{df}

We shall prove later that every minimal diagram is normal. Before that we need to show that
the disc graph of a normal diagram is hyperbolic in the sense of Section \ref{hg}.

\begin{lm} \label{one-three-gons}\ 
The disc graph of a minimal van Kampen diagram
has no one-gons.
\end{lm}
 
\noindent {\em Proof.}  Since a disc has only one $\kappa_i$-edge for any $i$,
there are obviously no one-gons in the disc graph.
$\Box$

\begin{lm} \label{bigons1} 
The disc graph corresponding to a minimal normal \vk diagram contains no bigons. 
\end{lm}

\noindent {\em Proof.} Indeed suppose that the disc graph of a minimal diagram $\Delta$ contains a bigon. The bigon consists of two discs $\Pi_1$ and $\Pi_2$ connected by two $\kappa$-bands. 
Let $i$ and $j$ be the indices of these $\kappa$-bands.
 Then either $i$ and $j$ are consecutive numbers $\ge 4$ or $i=4N, j=4$ (the case $j=4$, $i=4N$ is, of course, similar). In the first case Lemma \ref{bigon3} immediately gives a contradiction because our diagram is normal. 

In the second case notice that since the bigon must bound a face in the disc graph, there are no discs between the edges of the bigon. Consider the maximal $\kappa_3$-band starting on the contour of $\Pi_1$. This band cannot intersect the $\kappa_{4N}$-band or the $\kappa_4$-band.  Therefore it ends up on the contour of $\Pi_2$. Therefore the discs $\Pi_1$ and $\Pi_2$ are connected by a $\kappa_3$-band and by a $\kappa_4$-band, a contradiction with the assumption that our diagram is normal.
$\Box$ 

\medskip 
The following Lemma follows from Lemma \ref{3.2}.

\begin{lm} For every minimal normal diagram $\Delta$ with at least one disc there exists a disc $\Pi$  with $4N-6$ consecutive (modulo $4N$)
$\kappa$-bands
$\bb_1$,..., $\bb_{4N-6}$
starting on $\partial(\Pi)$ and ending on the boundary of $\Delta$. For every such disc let  $\Psi_\Delta(\Pi)$
be the subdiagram of $\Delta$ bounded by $\topp(\bb_1)$, $\bott(\bb_{{4N-6}})$,
$\partial(\Delta)$
and $\partial(\Pi)$, which contains $\bb_1$, $\bb_{{4N-6}}$ and does not contain
$\Pi$ (there is only one subdiagram in $\Delta$ satisfying these conditions).
Then there exists a disc $\Pi$ such that $\Psi_\Delta(\Pi)$
does not contain discs.
\label{lmm}
\end{lm}

\begin{lm} \label{bigon3} Suppose that a diagram $\Delta$ over the disc-based presentation contains two discs $\Pi_1$ and $\Pi_2$ connected by two $\kappa$-bands $\bb_1$ and $\bb_2$ (both start on the contour of $\Pi_1$ and end on the contour of $\Pi_2$) whose indices are consecutive numbers 
(mod $4N$) not equal to 1 or 2. Suppose also that the diagram $\Delta_1$ bounded by $\topp(\bb_1)$, 
$\bott(\bb_2)$ and the boundaries of $\Pi_1$ and $\Pi_2$ ($\Pi_1,\Pi_2\not\in\Delta_1$)
is reduced and normal (we do not assume that $\Delta$ is reduced or normal).
Then the number of  discs in the diagram can be reduced by 2
without changing the boundary label of $\Delta$. 
\end{lm}

\noindent {\em Proof.} Assume by contradiction that there exist two discs connected by two $\kappa$-bands $\bb_1$, $\bb_2$ and satisfying the conditions of the lemma. Let us denote the discs by $\Pi_1$ and $\Pi_2$. It is clear that the indices of the $\kappa$-bands must be consecutive numbers modulo $4N$.  We can assume that $\Pi_1$, $\Pi_2$, $\bb_1$ and $\bb_2$ are chosen in such a way that the diagram $\Delta_1$
is innermost. 

Let us first consider the case where there are no $\rho$-edges on the two
$\kappa$-bands $\bb_1$ and $\bb_2$. Then $\Delta_1$ has boundary label over $G_N(\sss)$.  Then by Lemma \ref{gn} $\Delta_1$ is a diagram over the disc based presentation of $G_N(\sss)$. 

By replacing the discs in $\Delta_1\cup\{\Pi_1,\Pi_2\}$ by diagrams over the (original finite) presentation of $G_N(\sss)$, and by applying Lemma \ref{bigons}, we can conclude that the number of discs in $\Delta$ can be reduced by 2 without changing the boundary label of $\Delta$. The new diagram has the same boundary label as $\Delta$ but the number of discs is smaller by 2. 

Now suppose that one of the $\kappa$-bands $\bb_1$ and $\bb_2$ contains $\rho$-cells. 
The diagram $\Delta_1$ does not contain discs. Indeed, since this diagram is normal, by Lemma \ref{lmm} there exists a disc inside $\Delta'$ with $4N-6$ consecutive external edges. At most two of the corresponding $\kappa$-bands can end on $\Pi_1$ or $\Pi_2$ (otherwise we would have two discs connected by two $\kappa$-bands and the diagram bounded by the discs and the bands smaller than $\Delta_1$). 
The other $4N-8>1$ $\kappa$-bands must intersect $\bb_1$ or $\bb_2$ which is impossible because two $\kappa$-bands cannot intersect.
 
Since by Lemma \ref{sigma-kappa}, $\Delta$ cannot contain $(\rho, \kappa)$- or $(\Theta,\kappa)$-annuli, 
every maximal $\Theta$-band and every maximal $\rho$-band crossing $\bb_1$ (resp. $\bb_2$) must cross $\bb_2$ (resp. $\bb_1$). 

Let $e_1$, $e_2$ be the start edges of the bands $\bb_1$ and $\bb_2$. Suppose that the
cell $\pi$ of $\bb_1$ attached to $e_1$ is a $(\Theta, \kappa)$-cell. Then the maximal $\Theta$-band $\ttt$ containing $\pi$ crosses $\bb_2$. The intersection cell of $\ttt$ and $\bb_2$ must be attached to $e_2$ because otherwise the $\Theta$- or $\rho$-band containing the cell of $\bb_2$ attached to $e_2$ forms a $(\Theta, \kappa)$- or $(\rho,\kappa)$-annulus with $\bb_2$.
By Lemma \ref{emptydiagram} the diagram between $\ttt$ and $\Pi_1$ is empty (because as we saw it does not contain discs), and we can move the $\Theta$-band $\ttt$ over the disc $\Pi_1$ enlarging the disc $\Pi_1$
(\setcounter{bigons}{\value{ppp}}see Fig. \thebigons). Let $\Delta'$ be the new diagram and $\Pi'_1$ be the enlarged disc. In the (reduced) diagram $\Delta'$, discs $\Pi'_1$ and $\Pi_2$ are still connected by two $\kappa$-bands $\bb_1'$ and $\bb_2'$ which are shorter than the bands $\bb_1$ and $\bb_2$.  

\bigskip
\begin{center}
\unitlength=1.50mm
\linethickness{0.4pt}
\begin{picture}(70.56,23.56)
\put(1.67,0.00){\dashbox{0.67}(29.78,23.56)[cc]{}}
\put(6.33,4.44){\framebox(64.22,14.22)[cc]{}}
\put(26.11,18.67){\line(0,-1){14.22}}
\put(55.44,18.67){\line(0,-1){14.22}}
\put(26.11,14.67){\line(1,0){29.33}}
\put(26.11,8.44){\line(1,0){29.33}}
\put(16.56,20.89){\makebox(0,0)[cc]{$\Pi_1'$}}
\put(12.11,10.89){\makebox(0,0)[cc]{$\Pi_1$}}
\put(23.67,16.67){\makebox(0,0)[cc]{$e_1$}}
\put(23.44,6.44){\makebox(0,0)[cc]{$e_2$}}
\put(28.78,16.89){\makebox(0,0)[cc]{$\pi$}}
\put(28.56,6.44){\makebox(0,0)[cc]{$\pi'$}}
\put(28.78,11.33){\makebox(0,0)[cc]{$\ttt$}}
\put(43.00,16.00){\makebox(0,0)[cc]{$\bb_1$}}
\put(42.78,6.22){\makebox(0,0)[cc]{$\bb_2$}}
\put(40.33,11.33){\makebox(0,0)[cc]{$\Delta_1$}}
\put(63.44,11.33){\makebox(0,0)[cc]{$\Pi_2$}}
\end{picture}
\end{center}
\begin{center}
\nopagebreak[4]
Fig. \theppp.
\end{center}
\addtocounter{ppp}{1}

Now let us assume that the cell $\pi$ is a $\rho$-cell. Then as in the previous paragraph, the intersection $\pi'$ of the $\rho$-band $\ttt$ containing $\pi$ with the $\kappa$-band $\bb_2$ must be 
a cell attached to $e_2$. Therefore the band $\ttt$ touches $\Pi_1$ along two consecutive $\kappa$-edges whose indices are not 1 or 2. Then by Lemma \ref{lm2} the label of the boundary of $\Pi_1$ is $\kappa(\sigma(c(u)))$ for some word $u$ which is equal to 1 in the group $G$. 
Consider the subdiagram of $\Delta$ consisting of the disc $\Pi_1$ and two cells $\pi$ and $\pi'$. By Lemma \ref{moverho} we can apply the $\rho$-band moving construction and replace $\Delta_1\cup\{\Pi_1,\Pi_2\}$ by a diagram $\Delta_1'$ with the same boundary label 
such that $\Delta_1'$ consists of the disc $\Pi_1$, $4N-2$
$(\rho, \kappa)$-cells attached to it along the $4N-2$ $\kappa$-edges not including $e_1$ and $e_2$ and a number of cells of smaller ranks. Let us replace $\Delta_1'$ for $\Delta_1\cup\{\Pi_1,\Pi_2\}$ in $\Delta$. Denote the new diagram by $\Delta_2$. The diagram $\Delta_1'$ contains only two discs $\Pi_1'$ and $\Pi_2$ connected by proper sub-bands $\bb_1'$ and $\bb_2'$ of the bands $\bb_1$ and $\bb_2$. After reducing the subdiagram $\Delta'$ bounded by the boundaries of $\Pi_1$ and $\Pi_2$, $\topp(\bb_1')$, $\bott(\bb_2')$ we still have two discs $\Pi_1'$ and $\Pi_2$ connected by two $\kappa$-bands with consecutive indices not equal to 1 and 2, which are shorter than $\bb_1$ and $\bb_2$. The subdiagram between $\bb_1'$ and $\bb_2'$ is normal because it is disc-free. 
Notice that the number of discs in the whole diagram does not change. Continuing in this manner, we shall make both $\kappa$-bands connecting our discs empty. But this situation has been considered before. $\Box$

Lemma \ref{bigon3} immediately implies

\begin{lm} \label{normal} Every minimal diagram is normal.
\end{lm}

\noindent{\em Proof.} By contradiction let $\Delta$ be a minimal non-normal diagram with the smallest possible number of cells. Then $\Delta$ contains two discs $\Pi_1$ and $\Pi_2$ connected by two $\kappa$-bands $\bb_1$ and $\bb_2$ with consecutive indices (distinct from 1 and 2). Since $\Delta$ contains minimal possible number of cells, the subdiagram of $\Delta$ bounded by $\Pi_1$, $\Pi_2$,  $\bb_1$, $\bb_2$ and not containing $\Pi_1$, $\Pi_2$ is normal. Then by Lemma \ref{bigon3} we get a contradiction with the minimality of $\Delta$. $\Box$

\begin{lm} \label{sigma-ann2} In a minimal \vk diagram $\Delta$ over the disc-based presentation
there are no
$\rho$-annuli and no $\Theta$-annuli.
\end{lm}
 
\noindent {\em Proof.}  Suppose by contradiction that there is a
$\rho$-annulus or a $\Theta$-annulus in our diagram. Let us pick an 
inner-most one and let us call it $\Sigma$. 

Assume first that there are no discs inside the diagram $\Delta_1$ bounded by the outer boundary of $\Sigma$.

Then $\Sigma$ cannot be a $\Theta$-annulus by Lemma \ref{theta}. 

If $\Sigma$ is a $\rho$-annulus then there are no 
$\kappa$-edges inside the diagram bounded by the annulus (otherwise there would be 
either $(\rho,\kappa)$-annuli ruled out by Lemma \ref{sigma-kappa} or $\kappa$-annuli ruled out by Lemma \ref{kappa-ann2}). Now we can get a contradiction by applying Lemma \ref{sigma-ann1}.
 
\begin{center}
\unitlength=1.50mm
\linethickness{0.4pt}
\begin{picture}(47.67,34.44)
\put(18.11,15.11){\oval(33.78,20.89)[]}
\put(18.56,15.11){\oval(28.89,17.33)[]}
\put(26.11,14.89){\circle{6.90}}
\put(26.33,14.89){\makebox(0,0)[cc]{$\Pi$}}
\put(7.22,9.78){\makebox(0,0)[cc]{$\Sigma$}}
\put(23.00,16.67){\line(-1,1){13.33}}
\put(11.00,30.89){\line(1,-1){13.11}}
\put(25.44,18.22){\line(-3,5){9.56}}
\put(17.22,34.44){\line(3,-5){9.78}}
\put(25.44,24.67){\makebox(0,0)[cc]{$\dots$}}
\put(28.33,17.78){\line(1,1){14.22}}
\put(43.44,30.89){\line(-1,-1){14.44}}
\put(29.44,15.33){\line(1,0){18.22}}
\put(47.44,14.00){\line(-1,0){18.00}}
\put(29.22,12.89){\line(1,-1){11.11}}
\put(39.44,1.11){\line(-1,1){11.11}}
\put(25.67,11.56){\line(0,-1){11.33}}
\put(6.56,30.89){\makebox(0,0)[cc]{$\bb_1$}}
\put(20.78,32.22){\makebox(0,0)[cc]{$\bb_2$}}
\put(28.33,0.89){\makebox(0,0)[lc]{$\bb_{4N-6}$}}
\put(27.44,0.00){\line(0,1){11.78}}
\end{picture}
\end{center}

\begin{center}
\nopagebreak[4]
\setcounter{band}{\value{ppp}}
Fig. \theppp.
\end{center}
\addtocounter{ppp}{1}

Let us now assume that the subdiagram $\Delta_1$ bounded by the outer boundary of the annulus $\Sigma$ contains discs.
The discs inside $\Delta_1$ form a disc graph, and by Lemma
\ref{lmm} at least one disc $\Pi$ in
this graph has  $4N - 6$ consecutive $\kappa$-bands 
$\bb_1,...,\bb_{4N-6}$ starting on the 
contour of $\Pi$ and ending on the boundary of $\Delta_1$ and such that the subdiagram $\Psi_{\Delta_1}(\Pi)$
is disc-free (see Fig. \thespokes). By Lemma \ref{sigma-kappa} every maximal $\rho$- or $\Theta$-band crossing $\bb_1$
must cross all $4N-6$ bands $\bb_1,...,\bb_{4N-6}$. Take the $\rho$- or $\Theta$-band $\Sigma'$ which touches $\Pi$ at the start edge of $\bb_1$. Then the intersection cell of $\Sigma'$ and $\bb_2$ must contain the start edge of $\bb_2$ because otherwise we would have a 
$(\rho,\kappa)$ or $(\Theta,\kappa)$ annulus ruled out by Lemma \ref{sigma-kappa}. Continuing in this manner, we conclude that the intersection cell of $\Sigma'$ and $\bb_i$, $i=1,...,4N-6$ 
contains the start edge of $\bb_i$. Therefore $\Sigma'$ touches $\Pi$ in $4N-6\ge 2N$ $\kappa$-edges. Since $\Psi_{\Delta_1}(\Pi)$ is disc-free, we can use Lemmas \ref{Theta-touch} and
\ref{lm3} and get a contradiction. $\Box$

\begin{lm} \label{Q-ann}  
In a minimal van Kampen diagram there are no $\overline Q$-annuli.
\end{lm}
 
\noindent {\em Proof.}  Assume by contradiction that there is a $\overline Q$-annulus $\qq$. Let $\Delta_1$ 
be the diagram bounded by the outer boundary of $\qq$. Then by Lemma \ref{lmm}, Lemma \ref{normal} and  Lemma \ref{kappa-ann2}  
 $\Delta_1$ cannot contain $\kappa$-edges (since a $\overline Q$-band cannot intersect a $\kappa$-band). It remains to apply Lemma \ref{qann}. $\Box$

\begin{lm} \label{Q-Th-ann} In a minimal van Kampen diagram 
there are no $(\overline Q, \Theta)$-annuli.
\end{lm}

\noindent {\em Proof.} Notice that if there were a disc inside the subdiagram $\Delta_1$ bounded by a $(\overline Q,\Theta)$-annulus, then the $\Theta$-band of this annulus would intersect $4N-6$ consecutive $\kappa$-bands starting on the contour of a disc (Lemma \ref{lmm} and Lemma \ref{normal}) which is impossible by Lemma \ref{Theta-touch} (this lemma applies to the $\Theta$-band which intersects the $4N-6$ $\kappa$-bands and is the closest to the disc). Therefore $\Delta_1$ does not contain discs.  Now the statement follows from Lemma \ref{Q-sigma-ann1}. $\Box$

\begin{lm} \label{Q-sigma-ann2}  
In a minimal van Kampen diagram there are no $(\overline Q,\rho)$-annuli and no $(\overline Q, d)$-annuli.
\end{lm}

\noindent {\em Proof.} By Lemma \ref{Q-sigma-ann1} it is enough to show that the subdiagram $\Delta_1$ bounded by the outer boundary of an innermost such annulus $\Sigma$ cannot have $\kappa$-edges.  

If the diagram $\Delta_1$  contains discs then by Lemmas \ref{lmm} and \ref{normal} there would exist
$4N-6$ $\kappa$-bands starting on the same disc and intersecting the $d$- or $\rho$-band of our annulus (the $\overline Q$-band cannot intersect with $\kappa$-bands), and the subdiagram $\Psi_{\Delta_1}(\Pi)$ is disc-free. But a $d$-band can intersect only with $\kappa_1$- and $\kappa_2$-bands and $4N-6>2$. So our annulus cannot be a $(\overline Q,d)$-annulus. But if $\Sigma$ is a
$(\overline Q,\rho)$-annulus, we can get a contradiction as in the proof of Lemma \ref{sigma-ann2}. 

Thus $\Delta_1$ contains no discs. Then by Lemma \ref{kappa-ann2} if $\Delta_1$ contains a $\kappa$-edge, the maximal $\kappa$-band of $\Delta_1$ containing this edge would form a $(\rho,\kappa)$- or $(d,\kappa)$-annulus with the $\rho$- or $d$-band of $\Sigma$. But this has been ruled out by Lemma \ref{sigma-kappa}. $\Box$

\begin{lm} \label{Y-sigma-ann}
Consider a minimal van Kampen diagram $\Delta$.
Then in $\Delta$, there are no $(\overline Y, \rho)$-annuli, no 
$(\overline Y, d)$-annuli, and no $(\overline Y, B)$-annuli.
\end{lm}

\noindent {\em Proof.} 
Let us consider a  
$({\overline Y}, \rho)$-, $({\overline Y}, d)$-, or
$({\overline Y}, B)$-annulus $\Sigma$ and suppose that $\Sigma$ is an innermost annulus of the corresponding form (that is if $\Sigma$ is a $(\overline Y, \rho)$-annulus then the diagram bounded by the outer boundary of this annulus does not contain other $(\overline Y, \rho)$-annuli, etc.). 
Let $\Delta_1$ be the subdiagram bounded by the outer boundary of the annulus. Let $\yyy$ be the
$\overline Y$-band of the annulus and let $\ttt$ be the $d$-, $\rho$- or $B$-band of the annulus.
Then as in the proof of Lemma \ref{Q-sigma-ann2} we prove that $\Delta_1$ contains no $\kappa$-edges (like $\overline Q$-bands, a $\overline Y$-band cannot cross a $\kappa$-band).

Suppose that $\ttt$ is a $\rho$-band. Then $\ttt$ cannot contain $\kappa$-cells, as we proved before. It also cannot contain $\overline Q$-cells other than the two intersection cells, otherwise there would be a $(\rho, \overline Q)$-annulus rulled out by Lemma \ref{Q-sigma-ann2} (recall that a $\overline Q$-band cannot cross a $\overline Y$-band). This and the fact that $\Delta$ does not contain $\overline Q$-annuli (Lemma \ref{Q-ann}) imply that $\Delta_1$ does not contain $\overline Q$-edges at all (otherwise $\Delta_1$ would contain a $\overline Q$-annulus ruled out by Lemma \ref{Q-ann}). If $\ttt$ contains $\overline Y$-cells distinct with the intersection cells then the maximal  $\overline Y$-bands starting on the boundary of $\ttt$ cannot end on a cell inside $\Delta_1$ (since as we just proved, $\Delta_1$ does not contain $\overline Q$-cells). This implies that these $\overline Y$-bands must form $(\overline Y, \rho)$-annuli with $\ttt$. This contradicts the assumption that our annulus is innermost. 
Thus $\ttt$ does not contain $\overline Y$-cells. Therefore $\ttt$ contains exactly two intersection cells. These cells must cancel which contradicts the minimality of $\Delta$.
This implies that $\Delta$ does not contain $(\overline Y, \rho)$-annuli.

Suppose that $\ttt$ is a $d$-annulus. If $\ttt$ contains $\overline Q$-cells, then $\Delta$ contains $(\overline Q, d)$-annulus, which is impossible by Lemma \ref{Q-sigma-ann2}. Therefore 
$\ttt$ does not contain $\overline Q$-cells. This and the absence of $\overline Q$-annuli (Lemma \ref{Q-ann}) implies that $\Delta_1$ does not contain $\overline Q$-cells. Therefore $\ttt$ cannot contain $\overline Y$-cells except for the intersection cells (otherwise as in the case when $\ttt$ is a $\rho$-band, we would get 
a contradiction with the assumption  that $\Sigma$ is innermost). Since every $d$-cell is either 
a $\overline Q$-cell or a $\overline Y$-cell, we can conclude that $\ttt$ consists of two intersection cells that cancel, a contradiction. This proves that $\Delta$ does not have $(\overline Y, d)$-annuli.

Finally assume that $\ttt$ is a $B$-band. Every cell in $\ttt$ is an $(A,B)$-commutation cell. 
A maximal $\overline Y$-band in $\Delta_1$ starting on the boundary of $\ttt$ cannot end in $\Delta_1$ because $\Delta_1$ cannot contain $\overline Q$-cells (otherwise we would have a
$\overline Q$-annulus which is impossible). Therefore every such $\overline Y$-band must form a 
$(\overline Y, B)$-annulus with $\ttt$ which contradicts the assumption that $\Sigma$ is innermost. Therefore $\ttt$ consists of two intersection cells that cancel.
This proves that $\Delta$ does not have $(\overline Y, B)$-annuli.
$\Box$

\begin{lm} \label{Y-bands}  
In a minimal van Kampen diagram, a
$\overline Y$-band cannot start and end on the same disc.
\end{lm}

\noindent {\em Proof.} Suppose there is a ${\overline Y}$-band $\yyy$ that starts 
and ends on the same disc $\Pi$ of $\Delta$. 

Let $\Delta_1$ be the subdiagram of $\Delta$ bounded by the top or the bottom path of $\yyy$ and the boundary of the disc, such that $\Delta_1$ contains the $\overline Y$-band but does not contain $\Pi$.  

The same argument as in the proof of the previous lemmas shows that $\Delta_1$ contains no 
$\kappa$-edges (recall that a $\overline Y$-band cannot cross a $\kappa$-band).
 
The subdiagram $\Delta_1$ contains no $\rho$-annuli (by Lemma \ref{sigma-ann1}), and no 
$d$-annuli (Lemma \ref{d}). Since (by Lemma \ref{Y-sigma-ann})
there are no $({\overline Y}, \rho)$-, $({\overline Y}, d)$-, or  $(\overline Y, B)$-annuli, $\Delta_1$ contains no $\rho$-, $d$- or $B$-edges (recall that the boundary of a disc does not have $\rho$-, $d$- or $B$-edges). 

Therefore $\Delta_1$ is a diagram over the presentation of $G_N(\sss)$. There are no $\overline Q$-edges inside $\Delta$. 
Indeed, there are no $\overline Q$-annuli by Lemma \ref{q}, a $\overline Q$-band cannot cross a
${\overline Y}$-band, and a $\overline Q$-band cannot start and end on the boundary of the same disc because 
no two $\overline Q$-edges on the path between any two consecutive  $\kappa$-edges on the boundary of a disc have mutually inverse labels). 

Now by Lemmas \ref{ar} and \ref{r}, $\Delta_1$ cannot have $\Theta$-edges either. Therefore $\Delta_1$ is empty. In particular, the $\overline Y$-band $\yyy$ is empty, which means that the label of the disc is not reduced, a contradiction with the definition of a disc (indeed the label of a disc is $\kappa(w)$ where $w$ is an admissible word for $\sss$.
$\Box$

\begin{lm} \label{Theta-Y} 
In a minimal van Kampen diagram
there are no $(\Theta, \overline Y)$-annuli.
\end{lm}
 
\noindent {\em Proof.}  Suppose, by contradiction, that there is a $(\Theta, \overline Y)$-annulus. Let  $\Sigma$ be the innermost one. Then the same argument as in the proof of Lemma \ref{Q-sigma-ann2} 
shows that the inner diagram of  $\Sigma$ does not contain discs. The $\Theta$-band of this annulus cannot
contain $\overline Q$-cells because the inner diagram $\Delta_1$ of $\Sigma$ does not contain terminal cells for $\overline Q$-bands (discs) and because Lemma \ref{Q-Th-ann} rules out $(\overline Q,\Theta)$-annuli. 
The subdiagram $\Delta_1$ does not contain $\overline Q$-cells because $\Sigma$ is innermost. This $\Theta$-band also cannot contain $\overline Y$-cells 
(except for the intersection cells) because our annulus is innermost and because the inner diagram of $\Sigma$ does not contain terminal cells for $\overline Y$-bands (discs and $\overline Q$-cells) . Therefore the $\Theta$-band consists of two intersection cells that cancel.
 $\Box$ 

\begin{lm} \label{Y-ann} 
In a minimal van Kampen diagram, 
there is no ${\overline Y}$-annuli.
\end{lm}
 
\noindent {\em Proof.}  Indeed, a $\overline Y$-annulus cannot contain discs in the inner diagram (Lemmas \ref{lmm} and \ref{normal}) and other $\kappa$-cells (by Lemma \ref{kappa-ann2}). Therefore this inner diagram cannot contain terminal cells for $d$-bands. Thus by Lemmas \ref{Y-sigma-ann} and \ref{Theta-Y} an $(\overline Y)$-annulus cannot have $(\overline Y, \rho)$-, $(\overline Y, d)$-, or $(\Theta, \overline Y)$-cells. 
Therefore it consists of $(A,B)$-cells. Therefore the inner boundary and the outer boundary of the annulus have equal labels. Thus we can remove the annulus by identifying its inner and outer boundaries. This contradicts the minimality of our diagram. 
$\Box$
 
The following Lemma is similar to Lemma 11.2 in \cite{Ol1} and Lemma 11.8 from \cite{SBR}. 

\begin{lm} Let $\rr_1$,...,$\rr_n$ be maximal $\overline Y$-bands starting on a path $p$ where $p$ is either a $\overline Y$-subpath of the boundary of a disc $\Pi$ or a subpath of 
$\partial(\Delta)$ containing no $\kappa$-edges and no $\overline Q$-edges. Suppose that the end edges of all $\rr_i$ are on the contours of $(\overline Q,\Theta)$-cells belonging to the same $\Theta$-band $\ttt$.
Then $n\le c$ where $c$ is two times the maximal number of $\overline Y$-letters in a $(\overline Q,\Theta)$-relation.
\label{11.2}
\end{lm}

\noindent {\bf Proof.} Indeed, if $n>c$ then there are three $\overline Y$-bands, say,  $\rr_1$, $\rr_2$, $\rr_3$ starting on $p$ and ending on three different cells $\pi_1$, $\pi_2$ and $\pi_3$ of $\ttt$. We can assume that $\pi_2$ is between $\pi_1$ and $\pi_3$. 
Consider the minimal subdiagram $\Delta_1$ of our diagram containing  $\overline Y$-bands
$\rr_1, \rr_2, \rr_3$, the minimal subpath of the path $p$ containing the starting edges of $\rr_1$, $\rr_2$, $\rr_3$, and the part of the band $\ttt$ between $\pi_1$ and $\pi_3$ (\setcounter{triple}{\value{ppp}}see Fig. \thetriple). 
Then $\Delta_1$ has no $\kappa$-edges on its contour. Therefore $\Delta_1$ does not contain discs by Lemmas \ref{lmm} and \ref{normal}. Therefore the maximal $\overline Q$-band $\qq$ in $\Delta_1$ containing $\pi_2$ divides $\Delta_2$ into two parts (that is if we delete the $\overline Q$-edges from $\qq$, the diagram $\Delta_2$ will fall into two pieces). 
The subpath of the path $p$ containing the start edges of $\rr_1,\rr_2,\rr_3$ is contained in one of these parts since it does not contain $\overline Q$-edges.  The cells $\pi_1$ and $\pi_3$ belong to different parts because $\qq$ cannot intersect $\ttt$ twice (Lemma \ref{Q-Th-ann}). Since the $\pi_1$ and $\pi_3$ are connected with the cells on $p$ by $\overline Y$-bands, one of these bands must intersect $\qq$ 
which is impossible (a $\overline Q$-band cannot cross a $\overline Y$-band).
 $\Box$

\medskip
\begin{center}
\unitlength=1.50mm
\linethickness{0.4pt}
\begin{picture}(45.67,25.78)
\put(0.11,0.22){\line(1,0){40.22}}
\put(8.11,15.11){\line(1,0){25.11}}
\put(11.22,15.11){\circle*{0.99}}
\put(29.67,15.11){\circle*{0.89}}
\bezier{136}(5.00,0.22)(1.22,23.56)(10.56,20.22)
\bezier{68}(10.33,20.22)(17.22,15.33)(10.33,10.22)
\bezier{64}(10.33,10.22)(3.89,8.67)(7.22,17.11)
\bezier{40}(7.22,17.11)(10.78,20.67)(11.22,15.33)
\bezier{200}(36.56,0.22)(45.67,25.78)(23.67,20.67)
\bezier{36}(23.67,20.67)(18.33,18.44)(19.67,15.11)
\put(19.44,15.11){\circle*{0.89}}
\bezier{84}(21.22,0.22)(31.44,5.11)(29.67,15.11)
\put(16.11,20.67){\line(3,-5){7.33}}
\put(9.22,13.11){\makebox(0,0)[cc]{$\pi_1$}}
\put(17.89,13.33){\makebox(0,0)[cc]{$\pi_2$}}
\put(29.22,17.11){\makebox(0,0)[cc]{$\pi_3$}}
\put(21.89,7.33){\makebox(0,0)[cc]{$\qq$}}
\put(34.78,15.11){\makebox(0,0)[cc]{$\ttt$}}
\put(5.89,5.56){\makebox(0,0)[cc]{$\rr_1$}}
\put(15.89,5.56){\makebox(0,0)[cc]{$\Delta_1$}}
\put(28.56,3.33){\makebox(0,0)[cc]{$\rr_2$}}
\put(36.78,7.56){\makebox(0,0)[cc]{$\rr_3$}}
\put(42.33,0.22){\makebox(0,0)[cc]{$p$}}
\end{picture}
\end{center}
\begin{center}
\nopagebreak[4]
Fig. \theppp.
\end{center}
\addtocounter{ppp}{1}
\medskip

\begin{lm}\label{11.2y} Let $\Delta$ be a minimal diagram and let $\partial(\Delta)=pq$ where $p$ is a $B$-path. Then at most 2 different maximal $B$-bands in $\Delta$ starting on $p$ end on 
contours of cells of the same $\overline Y$-band. 
\end{lm}

\noindent {\em Proof.} The proof is similar to the proof of Lemma \ref{11.2}: one needs only replace there $\overline Q$ by $d$ and $\overline Y$ by $B$ and notice that the minimal subdiagram $\Delta_1$ considered in the proof of Lemma \ref{11.2} does not contain cells which are terminal for $d$-bands (i.e. $(\kappa,\rho)$-cells). 
$\Box$

\begin{lm}\label{embedd} The homomorphism from $G_b$ onto the subgroup $\langle B\rangle$ induced 
by the identity map on $B$ (see Lemma \ref{Main diagram})
is an isomorphism.
\end{lm}

\noindent {\em Proof.} Consider a minimal diagram $\Delta$ over the disc-based presentation of $H$ with boundary label $u$ over $B$. By the \vk Lemma it is enough to proof that $u=1$
in $G_b$. 

Since the boundary label of $\Delta$ does not contain $\kappa$-edges, $\Delta$ does not contain discs (Lemmas \ref{lmm} and \ref{normal}). Since the word $u$ does not contain letter $\rho$, $\Delta$
does not contain $\rho$-edges (Lemma \ref{sigma-ann2}). Since $u$ does not contain letter $d$
and $\Delta$ does not contain $d$-annuli (Lemma \ref{d}), $\Delta$ does not contain $d$-edges
(maximal $d$-bands in $\Delta$ cannot start or end on cells in $\Delta$ because $\Delta$ does not contain $\rho$-cells).
Similarly we can exclude $\Theta$- and $\overline Q$-edges. Therefore $\overline Y$-bands in $\Delta$ cannot start (end) on cells from $\Delta$. Since $u$ does not contain letters from $\overline Y$, $\Delta$ cannot contain $\overline Y$-edges. Therefore all cells in $\Delta$ are
$G_b$-cells, so $u=1$ in $G_b$ by the \vk lemma. $\Box$

\begin{lm} \label{B-bridge} 
In a minimal van Kampen diagram,
no $B$-band can start and end on $G_b$-cells.
\end{lm}
 
\noindent {\em Proof.}  Assume, by contradiction, that there is a $B$-band $\bb$ that
starts and ends on $G_b$-cells. This $B$-band consists
of $(A,B)$-commutation cells.
There are two possibilities: the two $G_b$-cells at the ends of the $B$-band
could be different, or they could be the same. 

Let us first consider the case of a $B$-band $\bb$ that starts and ends on the same
$G_b$-cell $\pi$. The
inner diagram $\Delta_1$ bounded by $\pi$ and $\bb$ has boundary
label over $A \cup B$. Since $A\subset G_N(\sss)$, $G_N(\sss)$ is a retract of $H$ (see
Lemma \ref{gn} and a remark after the proof of this lemma), and $B$ is contained in the kernel of this retraction, we can conclude that the $A$-word $u$ and the $B$-word $v$ written on the boundary of $\Delta_1$ are equal to 1 in $H$. By Lemma \ref{gammapath} and Lemma \ref{gn}
$A$ is a set of free generators of $\langle A\rangle$ (we can apply these lemmas because 
a diagram with boundary label over $A$ cannot have discs by Lemma \ref{lmm} and \ref{normal}). 
Therefore $u=1$ in the free group. 
Therefore if $\bb$ is not empty it must contain a pair of $(A,B)$-cells that cancel. 
Thus $\bb$ is empty which means that $\pi$ touches itself. Then the subdiagram $\Delta_1$ of $\Delta$ has boundary label $v$. Let $\Delta'$ be the smallest subdiagram of $\Delta$ containing the $G_b$-cell $\pi$. Then $\Delta'$ contains $\Delta_1$ and $\pi$. The boundary label of $\Delta'$ is a word $w$ over $B$. By the \vk Lemma $w=1$ in $H$. Then by Lemma \ref{embedd}
$w=1$ in $G_b$, so we can extract $\Delta'$ from $\Delta$ and replace it by a diagram with one $G_b$-cell. If $\Delta_1$ contains cells, thius operation would reduce the type of $\Delta$ which is impossible since $\Delta$ is minimal. Therefore $\Delta_1$ contains no cells. But this means that the boundary label of $\pi$ is not a cyclically reduced word which contradicts the definition of a $G_b$-cell. Thus this case is impossible.

Let us now consider a $B$-band that starts and ends on two different 
$G_b$-cells; the $B$-band consists entirely of $(A, B)$-commutation cells, or is empty (consists of just one edge). 
Let $bu_b$ be the label (over $B^{\pm 1}$) of one of the $G_b$-cells, and 
let $b^{-1}v_b$ be the label of the other $G_b$-cell; $b$ is the common 
$B$-letter of all the $(A, B)$-cells of the $B$-band; let $w_a$ (over 
$A$) be the label of the top (or bottom) path of the $B$-band 
(it may be equal to 1, if the $B$-band is empty).  

The label of the boundary of the subdiagram consisting of the two $G_b$-cells and the
$B$-band is \, $u_b w_a v_b w_a^{-1}$ (and is = 1 in $H$). Modulo 
$(A,B)$-commutativity, $u_b w_a v_b w_a^{-1}$ is equal to $u_b v_b$ 
(hence this is also = 1 in $H$). 
Thus, in our van Kampen diagram we can remove the two $G_b$-cells and the
$B$-band, and replace them by the following subdiagram, whose boundary label
is still $u_b w_a v_b w_a^{-1}$:  
a $G_b$-cell with boundary label $u_b v_b$ and $(A,B)$-commutation cells, that
fill in the space between the $G_b$-cell and the closed
path labeled by $u_b w_a v_b w_a^{-1}$. 
\setcounter{ab}{\value{ppp}}See Fig. \theab. Thus we can lower the type of $\Delta$ which contradicts its minimality.
\\
\unitlength=1.50mm
\linethickness{0.4pt}
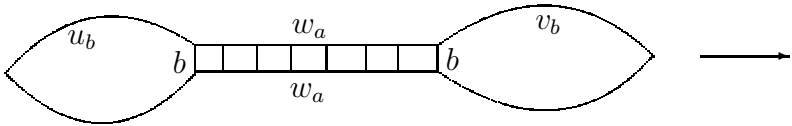
\begin{picture}(76.83,21.17)
\bezier{96}(24.17,11.17)(16.00,2.33)(7.33,11.17)
\bezier{92}(7.33,11.17)(14.83,19.83)(24.17,13.67)
\put(24.17,11.00){\line(0,1){2.33}}
\put(24.17,13.67){\line(1,0){21.50}}
\put(45.67,13.67){\line(0,-1){2.33}}
\put(45.67,11.33){\line(-1,0){21.50}}
\put(24.17,11.33){\line(0,1){2.33}}
\bezier{100}(45.67,13.67)(56.17,21.17)(64.83,12.67)
\bezier{100}(45.67,11.33)(56.83,3.83)(64.67,12.67)
\put(26.67,11.33){\line(0,1){2.33}}
\put(29.67,11.33){\line(0,1){2.33}}
\put(32.67,11.33){\line(0,1){2.33}}
\put(35.83,11.33){\line(0,1){2.33}}
\put(39.33,11.33){\line(0,1){2.33}}
\put(42.17,11.33){\line(0,1){2.33}}
\put(69.00,12.67){\vector(1,0){7.83}}
\put(22.83,12.17){\makebox(0,0)[cc]{$b$}}
\put(47.00,12.50){\makebox(0,0)[cc]{$b$}}
\put(55.50,15.50){\makebox(0,0)[cc]{$v_b$}}
\put(14.17,14.17){\makebox(0,0)[cc]{$u_b$}}
\put(34.33,15.17){\makebox(0,0)[cc]{$w_a$}}
\put(34.17,9.50){\makebox(0,0)[cc]{$w_a$}}
\end{picture}


\linethickness{0.4pt}
\begin{picture}(66.17,24.83)
\bezier{96}(36.00,11.33)(27.83,2.50)(19.17,11.33)
\bezier{92}(19.17,11.33)(26.67,20.00)(36.00,13.83)
\bezier{100}(36.00,13.67)(46.50,21.17)(55.17,12.67)
\bezier{100}(36.00,11.33)(47.17,3.83)(55.00,12.67)
\put(36.00,11.17){\line(0,1){2.50}}
\bezier{152}(36.00,1.67)(8.67,1.50)(8.33,12.33)
\bezier{140}(8.17,12.33)(14.67,24.83)(36.00,24.33)
\bezier{152}(35.83,24.33)(63.67,23.00)(66.17,12.67)
\bezier{160}(36.00,1.67)(64.50,1.50)(66.17,12.50)
\put(35.00,4.00){\makebox(0,0)[cc]{$(A,B)$-commutation cells}}
\put(35.17,12.50){\makebox(0,0)[cc]{$b$}}
\put(24.00,13.67){\makebox(0,0)[cc]{$u_b$}}
\put(45.33,15.33){\makebox(0,0)[cc]{$v_b$}}
\put(30.17,22.17){\makebox(0,0)[cc]{$u_bw_av_bw_a\iv$}}
\end{picture}

\begin{center}
\nopagebreak[4]
Fig. \theppp.
\end{center}
\addtocounter{ppp}{1}

\section{Estimating Areas in the Disc-based Presentation} \label{sec5}

Fix a minimal diagram $\Delta$ over the disc-based presentation of $H_N(\sss)$. Let $n$ be the length of the boundary of $\Delta$. We will assume that $n>3$. We are going to estimate the area of $\Delta$ in terms of $n$. Let $k$ be the number of components of the vector $Y$ of the $S$-machine $\sss$.
 
\begin{lm} \label{diskOn}
The number of discs in $\Delta$ is $O(n)$, 
and the number of maximal $\kappa$-bands is $O(n)$.
\end{lm}
 
\noindent {\em Proof.}  Let $\Gamma$ be the disc graph of $\Delta$, and let $v_0, v_1,...,v_m$ be the vertices of this graph, where $v_0$ is the external vertex. Without loss of generality we can assume that $m\ge 1$. Then as we know $\Gamma$ is an $\ell$-graph with $\ell=4N-3$. Let $d_0, d_1,...,d_m$ be the degrees of vertices $v_0,v_1,...,v_m$. Then by Lemma \ref{3.4} we have that 
$$d_0\ge 6-6m+\sum_{i=1}^m d_i.$$
Notice $d_i=4N-3$ for $i\ge 1$. Therefore 
$$d_0\ge 6-6m+(4N-3)m.$$
Since $N\ge 2$ we deduce that $d_0\ge m$. Notice that obviously $d_0\le n$. Therefore the number $m$ of discs in $\Delta$ does not exceed $n$.

It is well known (see for example Lemma 10.2 from \cite{Ol89}) that in any planar graph without 1-gons and bigons the number of vertices $n_0$ and the number of edges $n_1$ satisfy the following inequality: $n_1\le 3n_0$. Therefore the number of maximal $\kappa$-bands starting on discs is at most $6n$ (we consider a band and its inverse as different bands). Other $\kappa$-bands start and end on the boundary of $\Delta$, and their number is at most $n$. Since by Lemma \ref{kappa-ann2} $\Delta$ does not have $\kappa$-annuli, the total number of maximal $\kappa$-bands in $\Delta$ is at most $7n$.
$\Box$

In the rest of this section, we are going to estimate the number of non-disc and non-$G_b$-cells in $\Delta$ 
and the total perimeter of all discs and $G_b$-cells. The idea is the following. Every cell except for discs and $G_b$-cells is an intersection of two bands. For example, every $(\Theta,Q)$-cell is an intersection of a $\Theta$-band and a $Q$-band. So first we compute the 
number of different types of maximal bands in $\Delta$ using the absence of $Q$-annuli, $\Theta$-annuli, etc. proved in the previous section. Then we 
count the number of intersections of two bands. There we use the fact that there are no annuli formed by bands of different types. 

{\bf Remark 4.} 
\setcounter{spone}{\value{ppp}} Notice that a multiple intersection of a maximal $S$-band $\xxx$ and a maximal $T$-band $\yyy$ 
(where $S$ and $T$ are components of our subdivision of the set of generators of $H_N(\sss)$) does not formally  imply that $\yyy$ and $\ttt$ form a $(S,T)$-annulus. 
Indeed, $\yyy$ can start on the contour of a cell $\pi$, and $\ttt$ can go as a spiral around $\pi$ and intersect $\yyy$ many times as in Fig. \thespone. In this case $\yyy$ and $\ttt$ do not form a $(S,T)$-annulus because the start edge of $\yyy$ belongs to the region bounded by the medians of $\yyy$ and $\ttt$ (see the definition of an
$(S,T)$-annulus).
\medskip

\unitlength=1.00mm
\linethickness{0.4pt}
\begin{picture}(50.50,31.00)
\bezier{72}(36.83,21.83)(45.50,19.50)(44.83,10.17)
\bezier{76}(44.83,9.83)(42.17,0.17)(33.50,2.17)
\bezier{60}(33.17,2.17)(25.83,4.50)(25.50,12.17)
\bezier{84}(25.50,12.17)(25.50,23.17)(35.17,26.50)
\bezier{88}(35.17,26.67)(49.17,31.00)(50.50,23.33)
\put(38.00,14.17){\circle*{2.85}}
\put(38.83,15.00){\line(3,4){11.67}}
\put(53.33,30.33){\makebox(0,0)[cc]{$\yyy$}}
\put(23.00,15.67){\makebox(0,0)[cc]{$\ttt$}}
\put(35.33,12.00){\makebox(0,0)[cc]{$\pi$}}
\end{picture}
\begin{center}
\nopagebreak[4]
Fig. \theppp.
\end{center}
\addtocounter{ppp}{1}

But it is easy to see that a multiple intersection does imply the existence of $(S,T)$-annulus
if $\ttt$ does not start or end on a contour of a cell. 
Indeed let $\yyy=\yyy_1\yyy_2\yyy_3$, $\ttt=\ttt_1\ttt_2\ttt_3$ where $\yyy_2$ and $\ttt_2$ are parts of $\yyy$ and $\ttt$ between two consecutive intersections, and let $\pi_1$ and $\pi_2$ be these intersections (so these cells belong to both $\yyy_2$ and $\ttt_2$). Without loss of generality assume that $\pi_1$ precedes $\pi_2$ along both $\yyy$ and $\ttt$. Then the medians of $\yyy_2$ and $\ttt_2$ form a region $\Delta_1$ on the plane, and the cell $\pi$ is inside this region. Notice that the part $\ttt_1$ of $\ttt$ cannot escape from $\Delta_1$ without crossing $\yyy_2$ (it cannot cross itself). If $\ttt_1$ crosses $\yyy_2$ then the intersection cell must be between $\pi_1$ and $\pi_2$, so $\ttt_1$ and $\yyy_2$ form a $(S,T)$-annulus. But by our assumption, $\ttt$ does not start on a contour of a cell, so $\ttt_1$ must escape from $\Delta_1$, a contradiction. We are going to use this remark several times later.

\begin{lm} \label{sigmabandsOn} 
The number of maximal $\rho$- or $\Theta$-bands is $O(n)$.
The number of maximal $\overline Q$-bands is $O(n)$.
\end{lm}
 
\noindent {\em Proof.}  By Lemma \ref{sigma-ann2} every $\rho$-band and every
$\Theta$-band has to start and end on the contour of $\Delta$. This proves the first statement of the lemma. 
 
The $\overline Q$-bands all start and end on the outer boundary of the van Kampen diagram
or on discs. Since there are at most $n$ maximal $\overline Q$-bands starting on the boundary and 
there are $O(n)$ discs (by Lemma \ref{diskOn}), and each 
disc has a constant number of $\overline Q$-edges on the boundary (the number $k$ of components of vector $Q$ multiplied by $4N$), 
the second statement of the lemma follows. $\Box$
 
\begin{lm} \label{Num-s-k} 
The number of $(\rho, \kappa)$-cells in $\Delta$ is $O(n^2)$
The number of $(\Theta, \kappa)$-cells is also $O(n^2)$.
\end{lm}
 
\noindent {\em Proof.}  A $\rho$-band cannot cross the same
$\kappa_i$-band more
than once: this follows from the fact that there are no
$(\rho,\kappa)$-annuli (Lemma \ref{sigma-kappa}) and the fact that a maximal $\rho$-band does not start on the contour of a cell (see Remark 4).
Thus, each $\rho$-band has $O(n)$ intersections with $\kappa$-bands
(by Lemma \ref{diskOn}). Since there are $O(n)$ maximal $\rho$-bands
(Lemma \ref{sigmabandsOn}), the total number of intersections of
$\rho$-bands and $\kappa$-bands is $O(n^2)$.
 
The proof for $(\Theta, \kappa)$-cells is very similar: by Lemma \ref{sigma-kappa} there are no $(\Theta, \kappa)$-annuli and a maximal $\Theta$-band cannot start on the contour of a cell.
$\Box$
 
\begin{lm} \label{Num-k} 
The number of $\kappa$-cells in $\Delta$ is $O(n^2)$.
\end{lm}
 
\noindent {\em Proof.}  A $\kappa$-cell in our presentation is either a disc
(we have $O(n)$ of them, by Lemma \ref{diskOn}), or a $(\kappa, \Theta)$-cell
(there are $O(n^2) $ of them by Lemma \ref{Num-s-k}),
or a $(\kappa, \rho)$-cell (there are $O(n^2)$ of them by Lemma \ref{Num-s-k}).
Adding these bounds gives the result. $\Box$
 
\begin{lm} \label{Num-d-bands}
The number of maximal $d$-bands in $\Delta$ is $O(n^2)$.
\end{lm}
 
\noindent {\em Proof.}  Since $\Delta$ does not contain $d$-annuli (Lemma \ref{d}), a $d$-band starts either on the contour of $\Delta$ (which has $n$ edges) 
or on the contour of a $(\rho, \kappa)$-cell
(of which there are $O(n^2)$,  by Lemma \ref{Num-s-k}). $\Box$

\begin{lm} \label{Num-Q-Th} 
The number of
$(\overline Q, \Theta)$-cells in $\Delta$ is $O(n^2)$.
\end{lm}
 
\noindent {\em Proof.}
The proof is similar to the proof of Lemma \ref{Num-s-k}.
By Lemma \ref{Q-Th-ann} there are no $(\overline Q,\Theta)$-annuli and a maximal $\Theta$-band does not start on the boundary of a cell. Thus, each
$\overline Q$-band has $\leq n$
intersections with $\Theta$-bands (Lemma \ref{sigmabandsOn}). Since the number of maximal $\overline Q$-bands is $O(n)$ (Lemma  \ref{sigmabandsOn})
the total number of intersections of $\Theta$- and $\overline Q$-bands is $O(n^2)$. $\Box$
 
\begin{lm} \label{Num-Q-s}
The number of $(\overline Q,\rho)$-cells in $\Delta$ is $O(n^2)$.
\end{lm}
 
\noindent {\em Proof.}  There are no $\rho$-annuli
(Lemma \ref{sigma-ann2}), and no $\overline Q$-annuli (Lemma \ref{Q-ann}).
The number of maximal $\rho$-bands is  $O(n)$ and the number of maximal $\overline Q$-bands is $O(n)$, by Lemma \ref{sigmabandsOn}.
Since there are no $(\rho, \overline Q)$-annuli (Lemma \ref{Q-sigma-ann2}), 
each $\rho$-band intersects each $\overline Q$-band at most once. 
Since a maximal $\rho$-band cannot start on the contour of a cell, there are $O(n^2)$ intersections. $\Box$
 
\begin{lm} \label{Num-Q-d}
The number of $(\overline Q,d)$-cells is $O(n^3)$.
\end{lm}
 
\noindent {\em Proof.} The numbers of 
$\overline Q$-bands and $d$-bands are respectively $O(n)$ and $O(n^2)$ (Lemma \ref{sigmabandsOn} and Lemma \ref{Num-d-bands}). By Lemma \ref{Q-sigma-ann2} there are no $(d,\overline Q)$-annuli. We cannot apply Remark 4 now because a maximal $Q$-band and a maximal $d$-band can start on the boundary of a cell. So we should be more careful in proving that a maximal $\overline Q$-band $\yyy$ and a maximal $d$-band $\ttt$ cannot have a multiple intersection. 
Suppose that $\yyy$ starts on a disc $\Pi$, and $\ttt$ starts on a $\rho$-cell $\pi$ (\setcounter{sptwo}{\value{ppp}}see Fig.\thesptwo). 

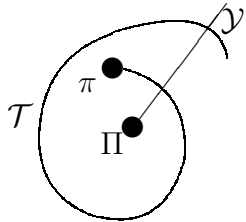
\unitlength=1.00mm
\linethickness{0.4pt}
\begin{picture}(51.33,31.00)
\bezier{72}(36.83,21.83)(45.50,19.50)(44.83,10.17)
\bezier{76}(44.83,9.83)(42.17,0.17)(33.50,2.17)
\bezier{60}(33.17,2.17)(25.83,4.50)(25.50,12.17)
\bezier{84}(25.50,12.17)(25.50,23.17)(35.17,26.50)
\bezier{88}(35.17,26.67)(49.17,31.00)(50.50,23.33)
\put(38.00,14.17){\circle*{2.85}}
\put(38.83,15.00){\line(3,4){11.67}}
\put(51.33,28.33){\makebox(0,0)[cc]{$\yyy$}}
\put(23.00,15.67){\makebox(0,0)[cc]{$\ttt$}}
\put(35.33,12.00){\makebox(0,0)[cc]{$\Pi$}}
\put(35.33,22.00){\circle*{2.98}}
\put(32.00,19.67){\makebox(0,0)[cc]{$\pi$}}
\end{picture}
\begin{center}
\nopagebreak[4]
Fig. \theppp.
\end{center}
\addtocounter{ppp}{1}

Let $\pi_1$ and $\pi_2$ 
be consecutive intersections of $\yyy$ and $\ttt$. 
Then as in Remark 4 let $\yyy=\yyy_1\yyy_2\yyy_2$ and
$\ttt=\ttt_1\ttt_2\ttt_3$ where $\yyy_2$ and $\ttt_2$ are parts of $\yyy$ and $\ttt$ between $\pi_1$ and $\pi_2$. Consider the subdiagram $\Delta_1$ of $\Delta$  bounded by the medians of $\yyy_2$ and $\ttt_2$. Then 
if $\yyy_2$ and $\ttt_2$ do not form an $(\overline Q, d)$-annulus, the start or end edge of $\ttt_2$ and $\yyy_2$ must belong to $\Delta_1$. Therefore $\yyy$ and $\ttt$ must start or end on contours of cells belonging to $\Delta_1$ (none of them can escape from $\Delta_1$ without forming a $(\overline Q,d)$-annulus). Without loss of generality we can assume that both the disc $\Pi$ and the $\rho$-cell $\pi$ are inside $\Delta_1$. But the contour of $\Delta_1$ cannot contain $\kappa$-edges because a $\overline Q$-band and
a $d$-band do not contain $\kappa$-cells. Therefore by Lemma \ref{lmm} $\Delta_1$ cannot contain discs, a contradiction. Therefore there are no multiple intersections of a $\overline Q$-band and a $d$-band. 
The statement of the lemma now follows immediately. $\Box$
 
\begin{lm} \label{Num-Q-cells}  
The number of $\overline Q$-cells is $O(n^3)$.
\end{lm}
 
\noindent {\em Proof.} In the presentation we use, every $\overline Q$-cell is either a disc (we counted 
$O(n)$ in Lemma \ref{diskOn}), or a $(\overline Q,\rho)$-cell (we have $O(n^2)$ 
of them by Lemma \ref{Num-Q-s}), or a $(\overline Q,d)$-cell (we have $O(n^3)\ )$,
by Lemma \ref{Num-Q-d}), or a $(\overline Q, \Theta)$-cell (we counted $O(n^2)$
in Lemma \ref{Num-Q-Th}). Adding these bounds yields the result. $\Box$

\begin{lm} \label{diskperim} 
The total number of  maximal ${\overline Y}$-bands and therefore the sum of the
perimeters of all the discs are $O(n^2)$.
\end{lm}
 
\noindent {\em Proof.} By Lemma \ref{Y-ann} $\Delta$ does not contain $\overline Y$-annuli. Therefore a maximal $\overline Y$-band can start on the contour of $\Delta$ or on a 
$(\Theta,\overline Q)$-cell or on a disc. The number of maximal $\overline Y$-bands starting on the contour of $\Delta$ is at most $n$, so we need to estimate only the number of other maximal $\overline Y$-band. Let $y(\Delta)$ be the number of maximal $\overline Y$-bands which start on discs in $\Delta$. Let $x(\Delta)$ be the number of maximal $\overline Y$-bands which start and end on discs of $\Delta$. Notice that we count every band which start and end on discs and the inverse of this band as two separate bands when we compute $x(\Delta)$ and $y(\Delta)$. 

Let us prove by induction on the number of discs in $\Delta$ that $x(\Delta)\le y(\Delta)/2$. 

By Lemma \ref{lmm} $\Delta$ contains a disc $\Pi$ and $4N-6$ consecutive $\kappa$-bands $\bb_1,...,\bb_{4N-6}$ starting on this disc and ending on the boundary of $\Delta$ and such that the diagram $\Psi_\Delta(\Pi)$ is disc-free. Create a diagram $\Delta_1$ by removing from $\Delta$ the subdiagram $\Psi_\Delta(\Pi)$ and the disc $\Pi$. The number of discs in $\Delta_1$ is fewer than in $\Delta$ by 1. Let $z$ be the number of $\overline Y$-edges on the path on the boundary of $\Pi$ between two consecutive $\kappa$-edges.

Let $p$ be the maximal subpath of the boundary of $\Pi$ which is not contained in $\Psi_\Delta(\Pi)$. Notice that $|p|$ contains $7z$ $\overline Y$-edges. 
The number $x(\Delta_1)$ of maximal $\overline Y$-bands in $\Delta_1$ which start and end on the contours of 
discs is $x(\Delta)$ minus the number of those maximal $\overline Y$-bands one of whose ends is on 
$p$ and the other end - on another disc. Therefore $$x(\Delta_1)\ge x(\Delta)-14z.$$

On the other hand, the number $y(\Delta_1)$ of the maximal $\overline Y$-bands  which start on discs of $\Delta_1$ is $y(\Delta)-4Nz$ (we subtract from $y(\Delta)$ the number of all maximal $\overline Y$-bands which start on $\Pi$). 

Since $\Delta_1$ has fewer discs than $\Delta$, we can conclude that $$x(\Delta_1)\le y(\Delta_1)/2.$$ Therefore $$x(\Delta)-14z\le (y(\Delta)-4Nz)/2.$$
Since $N\ge 7$, we conclude that $$x(\Delta)\le y(\Delta)/2.$$

This implies that $y(\Delta)$ does not exceed twice the number of maximal $\overline Y$-bands in $\Delta$ which start on a disc and end on the boundary of $\Delta$ or on a $(\overline Q,\Theta)$-cell. Since the number of $(\overline Q,\Theta)$-cells in $\Delta$ is $O(n^2)$ by Lemma \ref{Num-Q-Th}, the number of maximal $\overline Y$-bands starting on $(\overline Q,\Theta)$-cells is $O(n^2)$ (we have only finitely many $(\overline Q,\Theta)$-relations). Therefore $y(\Delta)=O(n^2)$ and the total number of maximal $\overline Y$-bands in $\Delta$ is also $O(n^2)$. $\Box$

\begin{lm} \label{NofSigmacells}  
The number of $\rho$-cells is $O(n^3)$.
\end{lm}

\noindent {\em Proof.} We already proved that the number of 
$(\rho, \kappa)$-cells and $(\rho, \overline Q)$-cells is $O(n^2)$ (Lemmas 
\ref{Num-s-k} and \ref{Num-Q-s}).  The only other $\rho$-cells are the 
$(\rho, \overline Y)$-cells.
The number of $\rho$-bands is $O(n)$ (Lemma \ref{sigmabandsOn}), and the 
number of ${\overline Y}$-bands is $O(n^2)$ (Lemma \ref{diskperim}). 
Since there are no $(\rho, \overline Y)$-annuli (Lemma \ref{Y-sigma-ann}), and a maximal $\rho$-band does not start or end on the contour of a cell,
it follows that each $(\rho, \overline Y)$-cell corresponds to a 
distinct pair of a maximal $\rho$-band and a maximal ${\overline Y}$-band. Thus, there are
at most $O(n^3)$ $\rho$-cells. $\Box$

\begin{lm} \label{Nofdcells} 
The number of $d$-cells is $O(n^4)$.
\end{lm}
 
\noindent {\em Proof.} We already proved that the number of
$(d, \overline Q)$-cells is $O(n^3)$ (Lemma \ref{Num-Q-d}), and the number of 
$\kappa$-cells is at most $O(n^2)$
(Lemma \ref{Num-k}). The only other $d$-cells are the
$(d, \overline Y)$-cells. 
The number of maximal ${\overline Y}$-bands is $O(n^2)$ (Lemma \ref{diskperim}), 
and the number of $d$-bands is $O(n^2)$ (Lemma \ref{Num-d-bands}). 
There are no $(d, \overline Y)$-annuli (Lemma \ref{Y-sigma-ann}). A similar argument as in Lemma \ref{Num-Q-d}
shows that a maximal $\overline Y$-band cannot intersect a maximal $d$-band more than once. 
Indeed, by Remark 4 the only ``bad" case is when both bands start or end on cells inside the region bounded by the medians of these bands. The cell where $d$-band starts or ends is a $\kappa$-cell. But the maximal $\kappa$-band containing this cell cannot cross the $\overline Y$-band or the $d$-band. 
This implies that the number of $(d, \overline Y)$-cells is $O(n^4)$. 
$\Box$

\begin{lm} \label{Num-Y-Th} 
The number of $(\overline Y, \Theta)$-cells is $O(n^3)$.
\end{lm}
 
\noindent {\em Proof.} We saw already that the number of $(\overline Q, \Theta)$-cells is $O(n^2)$ (Lemma \ref{Num-Q-Th}),
so we only have to count the $(\overline Y, \Theta)$-cells that have 
no $\overline Q$-edges.

We already proved that there are no $\Theta$-annuli
(Lemma \ref{sigma-ann2}), and no ${\overline Y}$-annuli (Lemma \ref{Y-ann}). 
We also saw that the number of $\Theta$-bands is
$O(n)$ (Lemma \ref{sigmabandsOn}), and the number of maximal  ${\overline Y}$-bands
is $O(n^2)$ (Lemma \ref{diskperim}). Since there are no 
$(\overline Y, \Theta)$-annuli (Lemma \ref{Theta-Y}),
and a maximal $\Theta$-band cannot start on the contour of a cell, every $\Theta$-band
intersects every maximal ${\overline Y}$-band at most once. Hence the number of 
$(\overline Y, \Theta)$-cells which do not contain $\overline Q$-edges (commutation cells) 
is $O(n^3)$. $\Box$ 

\begin{lm}  \label{Num-Th} 
The number of $\Theta$-cells is $O(n^3)$.
\end{lm}
 
\noindent {\em Proof.} A $\Theta$-cell is either a 
$(\overline Y, \Theta)$-cell (by Lemma \ref{Num-Y-Th}, there are $O(n^3)$
of these), or a $(\overline Q, \Theta)$-cell (by Lemma \ref{Num-Q-Th}, there are 
$O(n^2)$ of these), or a $(\kappa, \Theta)$-cell 
(by Lemma \ref{Num-s-k}, there are $O(n^2)$ of these).
Adding these bounds yields the Lemma.  $\Box$
 
\begin{lm} \label{Num-B-bands}
The number of maximal $B$-bands is $O(n^4)$ and the perimeter of each $G_b$-cell is $O(n^2)$.
\end{lm}
 
\noindent {\em Proof.} By Lemma \ref{B-bridge}, no $B$-band can both start and
end on a $G_b$-cell. Thus every maximal $B$-band must have at least one end either on 
the outer boundary of the van Kampen diagram (this accounts for
at most $n$ maximal $B$-bands), or on a $d$-cell. This accounts for at most $O(n^4)$ 
maximal $B$-bands, by Lemma \ref{Nofdcells}.  Thus there are at most $O(n^4)$ maximal $B$-bands. 

By Lemma \ref{11.2y} for every $G_b$-cell there are at most 3 maximal $B$-bands starting on the contour of this cell and ending on the contour of the same $\overline Y$-band. Notice that every edge on the boundary of a $G_b$-cell is the start edge of a $B$-band. Since by Lemma \ref{diskperim}
the total number of maximal $\overline Y$-bands is $O(n^2)$, the perimeter of each $G_b$-cell is $O(n^2)$.
 $\Box$
 
\begin{lm}  \label{Num-B-Y}  
The number of $(B, \overline Y)$-cells and the total number of 
$\overline Y$-cells are bounded by $O(n^6)$.
\end{lm}
 
\noindent {\em Proof.}
We saw (Lemmas \ref{Num-B-bands} and \ref{diskperim}) 
that there are $O(n^4)$ maximal $B$-bands and
$O(n^2)$ maximal ${\overline Y}$-bands.
Moreover, we cannot have $(\overline Y,B)$-annuli
by Lemma \ref{Y-sigma-ann}. Unfortunately, to exclude multiple intersections, we cannot use Remark 4 here because a $\overline Y$-band can start on a disc or a $\overline Q$-cell and a maximal $B$-band can start on a $G_b$-cell or on a $d$-cell. Suppose that a maximal $\overline Y$-band $\yyy$ intersects a maximal $B$-band $\ttt$ more than once. Then as in the proof of Lemma \ref{Num-Q-d} we can form a subdiagram $\Delta_1$ bounded by medians of a subband $\yyy_2$ of $\yyy$ and a subband $\ttt_2$ of $\ttt$. Then $\Delta_1$ must contain either a disc or a $\overline Q$-cell. The first option is impossible by Lemma \ref{lmm} because the boundary of $\Delta_1$ cannot contain $\kappa$-edges. The second option is also impossible because the boundary of $\Delta_1$ does not contain $\overline Q$-edges and there are no $\overline Q$-annuli in the diagram $\Delta$ (by Lemma \ref{qann}).

Therefore,  there are at most $O(n^6)$ $(B, \overline Y)$-cells.

Every ${\overline Y}$-cell is either a $(B, {\overline Y})$-cell (there are at
most $O(n^6)$ as we just proved), or a $\overline Q$-cell (there are at
most $O(n^3)$ of those, by Lemma \ref{Num-Q-cells}), or a $\rho$-cell
(there are at most $O(n^3)$ of those, by Lemma \ref{NofSigmacells}), or a 
$d$-cell (there are at most $O(n^4)$ of those, by Lemma \ref{Nofdcells}), or a
$\Theta$-${\overline Y}$-cell (there are at most $O(n^3)$ of those, by 
Lemma \ref{Num-Y-Th}). 
Summing these bounds yields the Lemma.
$\Box$

 
\section{Estimating the Isoperimetric Function of $H_N(\sss)$}
\label{sect6}

If a word $w$ over the generators of $H$ is equivalent to 1 in $H$ then there
exists a minimal van Kampen diagram over the disc-based presentation with boundary label $w$.  Let $n = |w|$ be the perimeter of this
van Kampen diagram. 

The number of $\kappa$-, $\overline Q$-, $\rho$-, $d$-, $\Theta$-, and 
${\overline Y}$-cells in this van Kampen diagram is at most $O(n^6)$ 
(by Lemmas \ref{Num-k}, \ref{Num-Q-cells}, \ref{NofSigmacells}, 
\ref{Nofdcells}, \ref{Num-Th}, \ref{Num-B-Y}). 

The only other cells in the diagram are the discs and the 
$G_b$-cells. By Proposition \ref{prop2} every disc of perimeter $p$ can be filled in by a $G_N(\sss)$-diagram with area at most
$O(T(O(p))^4)$. 
Moreover, by Lemma \ref {diskperim}, the sum of the perimeters of all the 
discs in our diagram is $O(n^2)$, for some positive constant $c$.   
Hence by superadditivity of $T^4$ we have: 
the total area of all the discs in the diagram, when replaced by $G_N(\sss)$-diagrams
is $O(T(O(n^2))^4)$. 

Every $G_b$-cell with perimeter $p$ can be replaced by a diagram with the same
boundary label, consisting of two discs (of perimeter $O(p)$ each), and 
$O(p^2)$ cells over the original presentation of $H$ ($O(p)$ $\rho$-cells and
$d$-cells, and $O(p^2)$ $(A, B)$ commutation cells); this follows from the
proof of Lemma \ref{Main diagram} and Fig. \thefigone. Let us fix one of these diagrams for each $G_b$-cell in $\Delta$. In Lemma \ref{Num-B-bands}, we
proved that the perimeter of each $G_b$-cell does not exceed $C_1n^2$ for some constant 
$C_1$. Since the total perimeter of all $G_b$-cells does not exceed the number of 
maximal $B$-bands in $\Delta$, Lemma \ref{Num-B-bands} shows that this total 
perimeter does not exceed $C_2n^4$ for some constant $C_2$. 

Then the set of all $G_b$-cells can be divided into subsets $M_1$,...,$M_t$ where $t\le C_3n^2$, such that the sum of perimeters of $G_b$-cells in each subset is $\le C_1n^2$. 

Since as we saw, every $G_b$-cell with perimeter $m$ can be filled by a diagram (as in Fig. \thefigone) with area at most $C_4T^4(C_5m)$ (see Proposition \ref{prop2}) for some constants $C_4$ and $C_5$, and since $T^4$ is a superadditive function, the sum of areas of diagrams filling the $G_b$-cells in each group $M_i$ does not exceed $C_4T(C_5C_1n^2)^4$. Therefore the total sum of areas of these diagrams for all $G_b$-cells in $\Delta$ does not exceed 
$C_3C_4n^2T(C_5C_1n^2)^4$ which is equivalent to $n^2T(n^2)^4$.

Thus, we can transform a minimal diagram $\Delta$ with boundary label $w$ into a diagram over the original presentation of $H_N(\sss)$ with the same boundary label and area bounded by a function equivalent to 
$n^2T(n^2)^4$. Since the total area of discs is equivalent to $T(n^2)^4$ and the number of other cells is $O(n^6)$, we conclude that $n^2T(n^2)^4$ is an isoperimetric function of the group $H_N(\sss)$ (since $T(n)$ growth at least linearly).

{\bf Remark 5.} Our construction of group $H_N(\sss)$ can be viewed as a substantial development of the Novikov-Boone-Higman-Aanderaa construction (with improvements by Britton, Collins and Miller)\cite{Ro}. The original construction (and all other known constructions for the Higman embedding theorem) contains Baumslag-Solitar relations $s\iv xs=x^2$ which lead to exponential lower bounds for the isoperimetric function of the groups. In \cite{SBR} and here we 
eliminated the Baumslag-Solitar relations by replacing arbitrary Turing machines with $S$-machines. There exists another approach based on replacing Baumslag-Solitar relations with relations which define groups with polynomial Dehn functions. One could use, for example,  embeddings of Baumslag-Solitar groups and similar groups into groups with polynomial Dehn functions. One could use also the R. Thompson group $F$ \cite{McKT} which has polynomial Dehn function \cite{Guba} and an injective but not surjective endomorphism of linear growth. This approach will most probably give  upper bounds for isoperimetric functions bigger than $n^2T(n^2)^4$.

\section{Length Distortion}
\label{sect7}

In this section we shall prove that the natural embedding of $G$ into $H$ (identifying $G$ with $G_b$) has bounded length distortion. This means that there exists a constant $C>0$ such that for every element $z\in G_b$ its length in $G_b$ is at most $C$ times its length in $H_N(\sss)$.

\subsection{Shadows}
\label{shadows}

Let us fix $N=9ck$ where $c$ is the constant from Lemma \ref{11.2} (notice that $c$ is twice the maximal number of $\overline Y$ letters in relations of the finite presentation of $G_N(\sss)$ and does not depend on $N$) and $k$ is the number of components of the vector $\overline Q$ of the $S$-machine $\sss$. Then the degree of every internal vertex of $\Gamma(\Delta)$ is $36ck-3> 33ck$. By Lemmas \ref{one-three-gons} and \ref{bigons1} $\Gamma(\Delta)$ is an $\ell$-graph for  $\ell=33ck$. Consider an oval $p$ passing through an internal vertex $\Pi$. By Lemma \ref{3.5} $p$ 
contains exactly two external edges $\bb$ and $\bb'$ and defines a region $\ooo(p)$. The external edges of $p$ cut a subpath $\overline p=\overline p(\Pi,p)\subset \ooo(p)$ of the boundary $\partial(\Delta)$. (The  end edges of the corresponding $\kappa$-bands $\bb$ and $\bb'$ do not belong to $\overline p$.) The path $\overline p$ is called the {\em shadow} of the disc $\Pi$ determined by the oval $p$.

Let $W$ be a word over the alphabet of generators of $H_N(\sss)$. Then $|W|_{y}$  denotes the number of $\overline Y$-letters  in $|W|$. Let $|W|_{y,\theta}$ (resp. $|W|_{y,b,\theta}$),
mean the number of $\overline Y$-letters (resp. $(\overline Y\cup B)$-letters) plus $4c$ times the number of $\Theta$-letters in $W$.
Thus we count every $\Theta$-letter with weight $4c$. 

The following Lemma is similar to Lemma 11.3 of \cite{Ol1}. 

\begin{lm} Let $W$ be a word written on the shadow $\overline p$ of a disc $\Pi$ of a minimal diagram $\Delta$. Let the label of $\Pi$ be $\kappa(u)$ for some word $u$. Then $2|u|_{y,\theta}=2|u|_y\le |W|_{y,\theta}$.
\label{11.3}
\end{lm}
 
\noindent{\em Proof.} Induction on the number of vertices from $\Gamma(\Delta)$ belonging to $\ooo(p)$. 

By the definition of oval at least    
$1/2(36ck-3)+4 > 18ck$ edges which are going out of the vertex $\Pi$, go in $\ooo(p)$. The word $u$ contains $k_u\le k$ maximal $\overline Y$-subwords (some of which may be empty). If $u$ does not contain $\overline Y$-subwords then there is nothing to prove. So let us assume that $u$ contains non-empty $\overline Y$-subwords. 

Thus the total number of maximal $\overline Y$-subpaths of $q=\partial(\Pi)\cap \ooo(p)$
is at least $18ckk_u$.  Ignoring the  subpaths between the first two and the last two edges going from $\Pi$ into $\ooo(p)$, we are left with at least $14ckk_u>14ck$ middle subpaths. Let us denote $K=14ck$ of these middle subpaths by $q_1,...,q_{K}$. Let $\Delta(p)$ be the subdiagram of $\Delta$ consisting of all cells intersected by the oval $p$ and those which are in $\ooo(p)$.  Consider 2 cases.

(1) The number $N_\theta$ of maximal $\Theta$-bands in the diagram $\Delta(p)$ is at least $\frac{1}{2c}|u|_y$. Notice that each of these bands must start or end on the shadow $\overline p$.
Indeed otherwise by Lemma \ref{3.7} we either have a $(\Theta,\kappa)$-bigon ruled out by Lemma \ref{sigma-kappa} or the median of the $\Theta$-band would envelope a disc $\Pi$ of $\Gamma(\Delta)$. This would lead to a contradiction with Lemma \ref{Theta-touch} as in the proof of Lemma \ref{sigma-ann2}. Therefore $|W|_{y,\theta}\ge 4c N_\theta \ge 2|u|_{y}$.

(2) Now suppose that $N_\theta< \frac{1}{2c}|u|_y$. Recall that  by Lemma \ref{sigma-ann2} $\Delta$ does not contain $\Theta$-annuli. Consider the set $\Omega$ of maximal $\overline Y$-bands starting on $q_1,...,q_{K}$. Notice that these bands cannot intersect $\kappa$-bands, so they are inside $\Delta(p)$. By Lemma \ref{11.2}, for every $i$ from 1 to $K=14ck$, the number of maximal $\overline Y$-bands starting on $q_i$ which end up on cells of the same $\Theta$-band is at most $c$. Therefore at most $14c^2k$ bands from $\Omega$ can have their end edges on cells belonging to the same $\Theta$-band. Hence the number of bands from $\Omega$ ending on $\Theta$-cells cannot be larger than $\frac{14c^2k}{2c}|u|_y=7ck|u|_y$ because of the restriction on the number $N_r$ of maximal $\Theta$-bands. Therefore at least 
$14ck|u|_y-7ck|u|_y=7ck|u|_y$ bands from $\Omega$ end up on discs inside the diagram $\Delta(p)$ or on $\overline p$. (Notice that these $\overline Y$-bands cannot end on $\partial(\Pi)$ by Lemma \ref{Y-bands}.) Let us denote this set of $\overline Y$-bands by $\Omega_1$. 

Consider two bands $\ttt_1$ and $\ttt_2$ from $\Omega_1$ which end on the same disc $\Pi'$
(if such a disc exists). In the area bounded by the medians of these bands, there are no discs
by Lemma \ref{lmm}. Therefore $\ttt_1$ and $\ttt_2$ start on $q$ and end on $\overline Y$-parts of the boundary of $\Pi'$ separated by at most 1 edge of the disc graph, that is by at most 4 $\kappa$-edges (recall that $\kappa_1$-, $\kappa_2$-, and $\kappa_3$-bands are not edges of the disc graph). This and the fact that the degree of every internal vertex of a disc graph is $\ge 28$ allows us do define a derived oval $p'$ passing through $\Pi'$ (using the median of any $\overline Y$-band connecting $\Pi$ and $\Pi'$; \setcounter{derived}{\value{ppp}} see Fig. \thederived a). 

\bigskip
\begin{center}
\unitlength=0.80mm
\linethickness{0.4pt}
\begin{picture}(151.91,53.44)
\put(16.92,50.44){\line(1,0){56.67}}
\put(35.25,32.10){\circle*{4.00}}
\put(34.58,15.77){\line(-2,1){11.00}}
\put(23.58,21.44){\line(-1,2){6.67}}
\put(16.92,34.77){\line(0,1){15.67}}
\put(25.58,50.44){\line(0,-1){11.33}}
\put(25.58,39.10){\line(1,-1){8.00}}
\put(37.25,31.10){\line(1,0){6.33}}
\put(43.58,31.10){\line(1,1){8.33}}
\put(51.92,39.44){\line(0,1){11.33}}
\put(47.08,15.94){\oval(18.33,6.33)[]}
\put(34.92,15.44){\line(1,0){18.33}}
\put(53.25,15.44){\line(2,1){10.33}}
\put(63.58,20.44){\line(3,5){9.67}}
\put(73.25,36.77){\line(0,1){13.67}}
\put(36.58,30.44){\line(2,-1){3.67}}
\put(35.92,30.10){\line(1,-1){3.00}}
\put(33.25,31.10){\line(-1,0){4.00}}
\put(33.58,30.44){\line(-5,-2){4.00}}
\put(38.58,17.10){\line(-1,1){7.00}}
\put(31.58,24.10){\line(1,2){3.00}}
\put(34.92,30.10){\line(0,-1){4.67}}
\put(39.58,18.44){\line(-5,6){2.33}}
\put(37.92,16.44){\line(-2,1){3.67}}
\put(44.58,19.10){\line(0,1){3.00}}
\put(47.58,19.10){\line(0,1){3.00}}
\put(50.92,19.10){\line(0,1){3.33}}
\put(53.58,19.10){\line(1,2){1.67}}
\put(55.25,18.10){\line(1,1){3.00}}
\put(63.58,20.77){\circle*{1.33}}
\put(73.25,36.44){\circle*{2.00}}
\put(63.58,20.44){\circle*{2.00}}
\put(24.25,20.77){\circle*{2.00}}
\put(16.92,34.77){\circle*{1.89}}
\put(73.58,36.44){\line(-1,1){3.00}}
\put(73.25,36.77){\line(-1,0){3.00}}
\put(73.25,36.44){\line(-1,-1){2.67}}
\put(63.58,20.10){\line(0,1){3.67}}
\put(63.58,20.44){\line(-1,1){2.33}}
\put(63.58,20.44){\line(-1,0){3.00}}
\put(24.25,20.44){\line(6,1){3.33}}
\put(24.25,20.77){\line(1,1){2.67}}
\put(24.25,21.10){\line(0,1){2.67}}
\put(17.25,34.77){\line(5,-1){3.00}}
\put(17.25,34.77){\line(1,1){2.67}}
\put(16.58,34.10){\line(-1,-1){2.00}}
\put(16.25,34.44){\line(-1,0){2.67}}
\put(23.92,20.77){\line(-1,0){2.67}}
\put(24.58,20.44){\line(-1,-6){0.33}}
\put(63.58,20.10){\line(0,-1){2.67}}
\put(63.58,20.44){\line(6,-1){2.67}}
\put(73.58,36.44){\line(6,-1){2.67}}
\put(73.58,36.44){\line(-1,-6){0.33}}
\put(51.92,39.10){\circle*{1.89}}
\put(25.58,39.44){\circle*{2.00}}
\put(43.92,31.10){\circle*{1.89}}
\put(36.25,33.10){\line(6,5){3.33}}
\put(33.92,33.77){\line(-1,3){1.33}}
\put(118.75,16.27){\oval(19.67,7.00)[]}
\put(118.58,15.77){\line(-1,0){15.33}}
\put(103.25,15.77){\line(-5,3){11.67}}
\put(91.58,22.77){\line(-1,2){6.67}}
\put(84.92,36.44){\line(0,1){14.00}}
\put(84.92,50.44){\line(1,0){66.33}}
\put(151.25,50.44){\line(0,-1){14.33}}
\put(151.25,36.10){\line(-1,-1){13.67}}
\put(137.58,22.44){\line(-3,-1){19.67}}
\put(98.58,26.77){\circle*{5.20}}
\put(111.92,30.44){\circle*{4.85}}
\put(133.92,29.44){\circle*{5.20}}
\put(96.25,27.44){\line(-1,2){5.67}}
\put(90.58,38.44){\line(0,1){12.00}}
\put(100.25,28.77){\line(2,5){4.00}}
\put(104.25,38.77){\line(0,1){11.33}}
\put(106.92,35.77){\line(0,1){14.33}}
\put(113.58,32.10){\line(1,1){6.00}}
\put(119.58,38.10){\line(0,1){11.67}}
\put(129.25,34.77){\line(0,1){15.67}}
\put(142.92,32.10){\line(0,1){18.33}}
\put(103.25,15.77){\circle*{1.89}}
\put(91.92,22.10){\circle*{2.00}}
\put(84.92,35.77){\circle*{1.33}}
\put(137.25,22.44){\circle*{1.33}}
\put(151.25,35.77){\circle*{1.33}}
\put(90.58,38.77){\circle*{1.49}}
\put(104.25,39.10){\circle*{1.49}}
\put(106.92,35.77){\circle*{1.33}}
\put(119.58,38.10){\circle*{1.33}}
\put(129.25,34.77){\circle*{1.33}}
\put(143.25,32.10){\circle*{1.33}}
\put(124.58,38.77){\makebox(0,0)[cc]{$\dots$}}
\put(122.25,14.77){\makebox(0,0)[cc]{$\Pi$}}
\put(47.58,17.20){\makebox(0,0)[cc]{$\Pi$}}
\put(39.92,33.20){\makebox(0,0)[cc]{$\Pi'$}}
\put(58.25,31.77){\makebox(0,0)[cc]{$O(p)$}}
\put(54.42,44.10){\makebox(0,0)[cc]{$p'$}}
\put(70.25,25.77){\makebox(0,0)[cc]{$p$}}
\put(18.25,26.10){\makebox(0,0)[cc]{$p$}}
\put(51.92,53.44){\makebox(0,0)[cc]{$W\equiv\Lab(\overline p)$}}
\put(109.25,17.44){\line(-1,1){8.67}}
\put(111.25,19.44){\line(0,1){8.67}}
\put(113.58,19.77){\line(0,1){8.67}}
\put(124.92,19.44){\line(1,1){7.67}}
\put(112.25,19.44){\line(0,1){3.67}}
\put(112.25,28.44){\line(1,-6){0.67}}
\put(99.25,24.44){\line(1,-2){1.67}}
\put(97.25,24.44){\line(-1,-1){1.67}}
\put(109.92,28.44){\line(-1,-1){1.67}}
\put(133.92,26.77){\line(0,-1){3.00}}
\put(98.25,29.10){\line(0,1){3.33}}
\put(111.92,32.44){\line(0,1){3.33}}
\put(135.25,31.44){\line(1,1){2.33}}
\put(116.25,19.44){\line(0,1){4.33}}
\put(119.58,19.77){\line(0,1){3.67}}
\put(122.58,19.77){\line(0,1){3.33}}
\put(110.25,18.44){\line(-3,4){2.00}}
\put(103.92,27.77){\makebox(0,0)[cc]{$\Pi'$}}
\put(117.25,29.77){\makebox(0,0)[cc]{$\Pi''$}}
\put(126.67,28.44){\makebox(0,0)[cc]{$\Pi^{(n)}$}}
\put(145.25,25.44){\makebox(0,0)[cc]{$p$}}
\put(85.92,25.44){\makebox(0,0)[cc]{$p$}}
\put(94.42,36.27){\makebox(0,0)[cc]{$p'$}}
\put(109.58,39.00){\makebox(0,0)[cc]{$p''$}}
\put(139.18,37.44){\makebox(0,0)[cc]{$p^{(n)}$}}
\put(97.25,44.10){\makebox(0,0)[cc]{$O(p')$}}
\put(112.92,44.10){\makebox(0,0)[cc]{$O(p'')$}}
\put(135.92,44.77){\makebox(0,0)[cc]{$O(p^{(n)})$}}
\put(97.25,52.77){\makebox(0,0)[cc]{$W'$}}
\put(112.92,52.77){\makebox(0,0)[cc]{$W''$}}
\put(135.92,52.77){\makebox(0,0)[cc]{$W^{(n)}$}}
\put(37.92,43.77){\makebox(0,0)[cc]{$O(p')$}}
\put(45.25,4.77){\makebox(0,0)[cc]{a}}
\put(117.25,4.77){\makebox(0,0)[cc]{b}}
\put(107.00,36.00){\line(1,-1){6.67}}
\put(129.34,34.67){\line(1,-1){6.00}}
\put(134.00,31.33){\line(0,1){5.67}}
\put(35.34,30.00){\line(1,-2){5.67}}
\put(143.34,32.00){\line(-2,-1){9.33}}
\end{picture}
\end{center}
\begin{center}
\nopagebreak[4]
Fig. \theppp.
\end{center}
\addtocounter{ppp}{1}

Let $p', ..., p^{(n)}$ be all ovals constructed for bands from the subset $\Omega_2$ of $\Omega_1$ consisting of all bands from $\Omega_1$ which end on discs $\Pi',...,\Pi^{(n)}$ (see Fig. \thederived b).
From Lemma \ref{3.8}, it follows that the shadows of these discs do not intersect and are inside the shadow of $\Pi$. The end edges of the bands from $\Omega_3=\Omega_1\backslash \Omega_2$ are
on the shadow $\overline p$ and cannot be on the shadows of the discs $\Pi',...,\Pi^{(n)}$. 
Therefore for the labels $W', ..., W^{(n)}$ of these shadows we have:

\begin{equation}\label{eq11.2}
|W'|_{y,\theta}+...+|W^{(n)}|_{y,\theta}+card(\Omega_3)\le |W|_{y,\theta}.
\end{equation}

By Lemma \ref{3.8} each of the regions $\ooo(p'),...,\ooo(p^{(n)})$ contains fewer 
discs than $\ooo(p)$. By the induction hypothesis, 

\begin{equation} \label{eq11.3}
2|u^{(i)}|_{y}\le |W^{(n})|_{y,\theta}
\end{equation}
for $i=1,...,n$ and word $u^{(i)}$ written on the boundary of $\Pi'$ between two consecutive $\kappa$-letters. At the same time

\begin{equation}\label{eq11.4}
5(|u'|_y+...+|u^{(n)}|_y)\ge card(\Omega_2)\ge 7ck|u|_y-card(\Omega_3)
\end{equation}
since different bands from $\Omega_2$ end on no more than 5 different $(\overline Y\cup \overline Q)$-parts of the boundary of a disc. From the inequalities \ref{eq11.2}-\ref{eq11.4}, it follows that

$$2|u|_y\le \sum_{i=1}^n\frac{10}{7}|u^{(i)}|_y+\frac{2}{7} card(\Omega_3)\le \sum_{i=1}^n|W^{(i)}|_{y,\theta}+card(\Omega_3)\le |W|_{y,\theta}.$$
$\Box$

\subsection{Bounded Distortion}

Let $T$ be the subgroup of $H_N(\sss)$ generated by all generators of $H_N(\sss)$ except $\kappa_1,...,\kappa_{4N}$ and generators from $\overline Q$. Then $G_b=\langle B\rangle$ is inside $T$. First we show that $G_b$ has bounded distortion in $T$.

\begin{lm}\label{lmm1} Let $\Delta$ be a minimal diagram whose  boundary label  is a word in generators of $T$. Let $\partial(\Delta)=pq$ where $u=\Lab(p)$ is a reduced word in the alphabet $B$. Then 
the length of the element $u$ in $G_b$ (that is the length of a geodesic which represents $u$ in the Cayley graph of $G_b$) does not exceed the length of $|\Lab(q)|_{y,b}$.
\end{lm}

\noindent {\bf Proof.} Let us assume that $\Delta$ has the smallest type among all minimal diagrams with $\partial(\Delta)=pq$ where $\Lab(p)$ is reduced and equal to $u$ in $G_b$. 

By Lemma \ref{lmm} the diagram $\Delta$ contains no discs since its boundary contains no $\kappa$-edges.
 
Consider any maximal $B$-band $\ttt$ which starts on $p$.  This band can end either on $p$ or on a $G_b$-cell or on $q$ or on a $d$-cell. 

In the first case we can get a contradiction as in the proof of Lemma \ref{B-bridge}: the $A$-word written on the top path of $\ttt$ is 1 in the free group, so either there are two cells in this band that cancel (which contradicts the minimality of $\Delta$) or the $B$-band is empty (which contradicts the assumption that $\Lab(p)$ is reduced). 

In the second case, again as in the proof of Lemma \ref{B-bridge}, we can replace the subdiagram consisting of $p$, the band $\ttt$ and the $G_b$-cell by a new subdiagram with the same boundary label and one $G_b$-cell which is attached to $p$ (this operation adds a number of $(A,B)$-commutation cells). Now we can remove the $G_b$-cell replacing $p$ by a new $B$-path $p'$ whose label is a reduced word which is equal to $u$ in $G_b$. The new diagram will have smaller type than $\Delta$ but $\partial(\Delta')=p'q$ where $\Lab(p')$ is reduced and equal to $u$ in $G_b$, a contradiction with the assumption that $\Delta$ has the smallest type among all such diagrams. 

Thus we can assume that every maximal $B$-band starting on $p$ ends either on $q$ or on the contour of a $d$-cell.

Let $\Omega$ be the set of all maximal $B$-bands starting on $p$ and let $\Omega_1$ be the subset of $\Omega$ consisting of all bands ending on a $d$-cell. Then by Lemma \ref{11.2y} at most two bands of $\Omega$ end on cells of the same maximal $\overline Y$-band. Since $\Delta$ does not contain $\overline Y$-annuli (Lemma \ref{Y-ann}), and $\overline Q$-cells (Lemma \ref{Q-ann}), every maximal $\overline Y$-band starts and ends on $\overline Y$-edges of the boundary of $\Delta$ (recall that $\Delta$ contains no discs) we have a one to one correspondence between bands from $\Omega_1$ and a subset of $\overline Y$-edges of $q$. 

Thus the number of edges of $p$ is smaller than the number of $\overline Y\cup B$-edges on $q$. $\Box$

\begin{lm} Let $u$ be an element of $G_b$ and let $v$ be a word in generators of $H_N(\sss)$. Suppose that $u=v$ in $H_N(\sss)$. Then the length of $u$ in $G_b$ does not exceed $|v|_{y,b,\theta}$.
\label{12.2}
\end{lm}

\noindent{\em Proof.} Let $\Delta$ be the diagram of the smallest type with the following property:
\medskip
(P) $\partial(\Delta)=pq\iv$ where $\Lab(p)$ is a word over the generators of $T$ and $\Lab(p)=u$ in $T$, $\Lab(q)\equiv v$.
\medskip
 We know by Lemma \ref{lmm1} that $|\Lab(p)|_{y,b}\ge |u|$. So it is enough to show that $|\Lab(p)|_{y,b}\le |v|_{y,b,\theta}$.

No edge in $p$ can belong to the contour of a $\overline Y$- or $G_b$-cell because otherwise by removing this cell we decrease the type of $\Delta$ preserving the property (P). Therefore every $B$-edge on $p$ also belongs to $q$. 
Thus if $p$ contains $B$-edges then by removing this edge we break our diagram into two parts;  these parts will have smaller ranks than $\Delta$, and we can consider them separately. Therefore we can assume that $p$ does not contain $B$-edges and in general $p$ and $q$ do not have any edges in common. 

Thus every $\overline Y$-edge on $p$ belongs to the contour of either a $(\overline Q,\Theta)$-cell or a 
disc. 

By Lemma \ref{sigma-ann2} $\Delta$ does not contain $\Theta$-annuli. Consider two cases.

(1) The number $N_\theta$ of maximal $\Theta$-bands in $\Delta$ is at least  
$\frac{1}{4c}|\Lab(p)|_y$ where $c$ is the constant from Lemma \ref{11.2}. 
Suppose that a $\Theta$-band $\ttt$ starts and ends on $p$. Then the same argument as in the proof of Lemma \ref{sigma-ann2} shows that the subdiagram $\Delta_1$ bounded by $\ttt$ and $p$ cannot contain $\kappa$-edges. It cannot contain $\overline Q$-edges also because $p$ does not contain $\overline Q$-edges and Lemma \ref{Q-Th-ann} rules out $(\overline Q, \Theta)$-annuli. Therefore the top and the bottom paths of $\ttt$ are words over the generators of $T$. Thus removing the subdiagram $\Delta_1$ from $\Delta$ would give us a diagram with smaller type still satisfying 
property (P). Therefore no $\Theta$-band in $\Delta$ starts and ends on $p$. Thus $q$ contains 
at least $\frac{1}{4c}|\Lab(p)|_y$ $\Theta$-edges. Since the weight of each of these edges is $4c$,
we get the desired  inequality $|\Lab(p)|_{y,b}\le |\Lab(q)|_{y,b,\theta}$.

(2) Now suppose that $N_\theta <\frac{1}{4c}|\Lab(p)|_y$.  By Lemma \ref{11.2} 
we deduce that at most $c$ $\overline Y$-edges of $p$ can belong to cells of the same $\Theta$-band. 
Therefore the number of $\overline Y$-edges of $p$ belonging to $\Theta$-cells is at most $\frac{1}{4}|\Lab(p)|_y$. 

Let $\Omega$ be the set of $\overline Y$-edges of $p$ which belong to discs. We have established that 
$|\Omega|\ge \frac{3}{4}|\Lab(p)|_y$. 

If some edges of $\Omega$ belong to the contour of the same disc then they cannot be separated
on this disc by $\kappa$-edges. To prove that, consider the subdiagram between the disc and $p$. 
The $\kappa$-bands which start on the boundary of this subdiagram would have nowhere to end (there are no discs in this subdiagram by Lemmas \ref{lmm} and \ref{normal}).

Therefore the $\overline Y$-edges of $p$ belong to one $(\overline Y\cup \overline Q)$-part of the boundary of any 
disc $\Pi$ which touches $p$. Therefore there exists an oval $z$ in the disc graph of $\Delta$ passing through $\Pi$ such that
no $\overline Y$-edge on $p$ belongs to the region $\ooo(z)$. In this case the shadow of the disc $\Pi$ is contained in $q$ because $p$ contains no $\kappa$-edges. By Lemma \ref{11.3} we have the following inequality for the label $Z$ of this shadow:

\begin{equation}\label{eq12.2}
|Z|_{y,\theta}\ge 2|w_\Pi|_{y}
\end{equation}
where $w_\Pi$ is the word written on the boundary of $\Pi$ between two consecutive $\kappa$-edges.

The oval $z$ can be chosen in such a way that every $\Gamma(\Delta)$-normal arc $x$ with end on an edge belonging to $\Pi$ and $p$ is 2-separated from $z$. Then ovals $z_1$ and $z_2$ constructed in this way for two such discs $\Pi_1$ and $\Pi_2$ will be diverging because the edges $e_1$ and $e_2$ where these discs are touching $p$ can be connected by suitable $\Gamma(\Delta)$-normal arc $x$. (Recall that there are no $\kappa$-edges between edges $e_1$ and $e_2$ on $p$.) 

By Lemma \ref{3.6} ovals $z_1$ and $z_2$ do not intersect in internal vertices of the disc graph $\Gamma(\Delta)$. Since the arc $x$ is outside regions $\ooo(z_1)$ and $\ooo(z_2)$, these regions also do not intersect. Therefore the shadows $q_1$ and $q_2$ of $z_1$ and $z_2$ do not intersect. 

Therefore we can sum up the inequalities \ref{eq12.2} for all discs touching $p$ and get 
$$|\Lab(q)|_{y,\theta} \ge 2\sum |w_\Pi|_{y} \ge 2\ \frac{3}{4}|\Lab(p)|_{y}\ge |\Lab(p)|_y$$
as desired. $\Box$.

Lemmas \ref{lmm1} and \ref{12.2} complete the proof that $G_b$ has bounded length distortion in $H_N(\sss)$.



\noindent Jean-Camille Birget\\
Department of Computer Science\\
University of Nebraska-Lincoln\\
birget@cse.unl.edu\\

\noindent Alexander Yu. Olshanskii\\
Department of Higher Algebra\\
MEHMAT. Moscow State University\\
olsh@nw.math.msu.su\\

\noindent Eliyahu Rips\\
Department of Mathematics\\
Hebrew University of Jerusalem\\
rips@math.huji.ac.il\\

\noindent Mark V. Sapir\\
Department of Mathematics\\
Vanderbilt University\\
http://www.math.vanderbilt.edu/$\sim$msapir\\

\end{document}